\documentclass[12pt]{amsart}
\usepackage{amsmath,amsfonts,amsthm,amssymb}
\def\dd{\delta}
\def\bb{\beta}
\def\e{\varepsilon}
\def\N{{\mathbb N}}
\def\R{{\mathbb R}}

\def\H{{\mathcal H}}
\def\TT{{\mathcal T}}

\def\id{\text{I} } 
\def\charp{\Sigma^{re}}
\def\charq{\Sigma^{im}}
\def\charxpq{\Sigma_x}
\def\charpq{\Sigma}
\def\ty{\tilde{y}}
\def\Sl{S_\l^2}

\def\dd{\delta}

\renewcommand{\l}{\lambda}
\renewcommand{\d}{\partial}
\renewcommand{\aa}{\alpha}

\newcommand{\etf}{e^{\tau \phi}}

\newcommand{\pim}{{p_{im}}}
\newcommand{\pr}{{p_{re}}}
\newcommand{\piw}{{p_{im}^w}}
\newcommand{\prw}{{p_{re}^w}}
\newcommand{\I}{\mathcal I}

 \newtheorem{theorem}{Theorem}
 \newtheorem{remark}{Remark}[section]
 \newtheorem{lemma}[remark]{Lemma}
 \newtheorem{cor}[remark]{Corollary}
 \newtheorem{prop}[remark]{Proposition}
 \newtheorem{definition}[remark]{Definition}

\begin{document}

\title{Dispersive estimates for principally normal pseudodifferential
  operators}

\author{Herbert Koch \and Daniel Tataru}

\address{Fachbereich Mathematik,Universit\"at Dortmund,  44221 Dortmund}
\email{koch@mathematik.uni-dortmund.de}

\address{Department of Mathematics, University of California at
  Berkeley, Berkeley, CA, 94596}
\email{tataru@math.berkeley.edu}
\thanks{The first author was supported in part by NSF Grant DMS-0226105}



\begin{abstract}
The aim of these notes is to describe some recent results concerning
dispersive estimates for principally normal pseudodifferential
operators. The main motivation for this comes from unique continuation
problems. Such estimates can be used to prove $L^q$ Carleman
inequalities, which in turn yield unique continuation results for 
various partial differential operators with rough potentials.
\end{abstract}

\maketitle

\section{Introduction}
Dispersive estimates are $L^q$ estimates for nonelliptic partial
differential operators which are a consequence of the decay properties
of their fundamental solutions. These decay properties follow from
spatial spreading of the singularities of the solutions. Since
solutions propagate in directions conormal to the characteristic set
of the operator, this spreading can be related to nonzero curvatures
of the characteristic set.  Dispersive estimates for constant
coefficient operators are closely related to the restriction theorem
in harmonic analysis.

Various types of dispersive estimates 
are known to be true for operators such as the wave
operator, the Schr\"odinger operator and the linear KdV, see 
Ginibre-Velo~\cite{GV}, Keel-Tao~\cite{KT}.
They have proved to be  useful in the study of nonlinear problems,
as well as of problems with unbounded potentials.

More recently, similar estimates have been obtained for wave operators
with variable coefficients, beginning with the smooth case in
Kapitanskii~\cite{759.35014}, Mockenhaupt, Seeger and
Sogge~\cite{MSS}, up to operators with $C^2$ coefficients in
Smith~\cite{Sm} and Tataru~\cite{MR1833146}, \cite{MR2003a:35120}.
Similar results were obtained for the Schr\"odinger equation in
Staffilani-Tataru~\cite{MR1924470} ($C^2$ coefficients) and in
Burq-Gerard-Tzvetkov~\cite{MR2003a:58055} (smooth coefficients).  In
the variable coefficient elliptic case one should also mention Sogge's
$L^q$ eigenfunction bounds on compact manifolds, see \cite{Sbook}.

All the above examples are operators with real symbols.
On the other hand, in unique continuation problems one is interested
in Carleman estimates. These are uniform weighted estimates with
respect to a family of exponential weights,
\[
\|\etf u\|_{L^q} \leq \|\etf P(x,D) u\|_{L^r}, \qquad \tau > \tau_0
\]
With the substitution $v = \etf u$ these can be rewritten as
\[
\|v\|_{L^q} \leq \| P(x,D+i\tau \nabla \phi) v\|_{L^r}, \qquad \tau > \tau_0
\]
Even if $P$ has constant coefficients and real symbol, the 
conjugated operator
\[
P_\phi = P(x,D+i\tau \nabla \phi)
\]
will have complex symbol and variable coefficients. In order for any
such estimates to hold $P_\phi$ must satisfy a so-called
pseudoconvexity condition. It is known that such a condition implies
$L^2$ estimates.  However, $L^r \to L^q$ estimates are considerably
more difficult to obtain. The first results in this direction were
obtained in Jerison-Kenig~\cite{MR87a:35058} for the Laplacian with a
polynomial weight. Later these were extended by
Sogge~\cite{MR91k:35068} to second order elliptic operators with
smooth coefficients. Further work of Wolff \cite{MR96c:35068} and of the
authors \cite{MR2001m:35075} addresses also the case of Lipschitz
coefficients, and $L^p$  gradient potentials.

To this one should add the work of Kenig-Ruiz-Sogge~\cite{MR88d:35037}
for second order constant coefficient operators, and the work of
Sogge~\cite{MR91m:35051} for parabolic operators with smooth
coefficients. The counterpart of the Jerison-Kenig estimates for
second order constant coefficient parabolic operators is proved in
Escauriaza~\cite{MR2001m:35135}.

All of the above mentioned results take advantage of the special form
of the operator in one way or another. On the other hand, it is clear
that only the geometry of the characteristic set should matter.

Motivated by problems in unique continuation and in local solvability,
in the present article we consider the problem of obtaining dispersive
estimates for operators which are principally normal.  However, of
independent interest is our parametrix construction for principally
normal operators, as well as the corresponding pointwise estimates for
the kernel of the parametrix. We only make assumptions on the geometry
of the characteristic set, and we also seek to use minimal regularity
for the symbols/coefficients.

An obstacle in applying our results to obtain Carleman estimates
for unique continuation problems is that the conjugated operator
$P_\phi$ introduced above does not satisfy the principal normality
condition. Fortunately the $L^2$ estimates which follow from
 the pseudoconvexity condition are strong enough so that they allow
spatial localization on a much smaller scale. This scale turns out to be
precisely the largest scale on which the principal normality 
survives.

To give the reader some idea of the results we obtain in this article,
we present some very simple examples. All these examples have constant
coefficients, however our results apply as well to operators with
variable coefficients.  In what follows $u$ is supported in the unit
ball in $\R^n$ and $\l > 1$. In many cases the support restriction is
easily removed by scaling.

 \begin{itemize} 
\item As a consequence of Theorem \ref{tp}  we have 
\[ 
\Vert u \Vert_{L^{\frac{2(n+1)}{n-1} }} \lesssim 
\l^{-\frac2{n+1}} \Vert (\Delta  + \l^2) u \Vert_{L^{\frac{2(n+1)}{n+3}}}.
\]
\item By Theorem \ref{tpq1} for any differential operator $Q(D)$ with
  constant coefficients and real symbol we have
\[ 
\Vert u \Vert_{L^{\frac{2(n+1)}{n-1} }} \lesssim 
\l^{-\frac2{n+1}}  \Vert (\Delta  + \l^2  + i \l^2 Q(D/\l)) u \Vert_{L^{\frac{2(n+1)}{n+3}}} 
\] 
\item Theorem \ref{tpq} implies  
\[ 
\Vert u \Vert_{L^{\frac{2(n+2)}{n}}} \lesssim \l^{-\frac4{n+2}}
\Vert (\Delta  + \l^2 + \l \d_1) u \Vert_{L^{\frac{2(n+2)}{n+4}}}. 
\]  
\item By Theorem \ref{tpdq}, for $|\delta|\le 1$
\[ 
\Vert u \Vert_{L^{\frac{2(n+2)}{n}}} \lesssim \delta^{-\frac{2}{n+2}
} \l^{-\frac4{n+2}} \Vert (\Delta +\l^2 + \delta \l \d_1) u
\Vert_{L^{\frac{2(n+2)}{n+4}}}.
\]  
\end{itemize} 

These results are slightly stronger than stated in our theorems.
However, there are obvious $L^2$ estimates for $\Delta + \l^2$
which allow to improve the estimates to the form stated above.

A new obstacle which we face in this analysis is that for principally
normal operators the characteristic set is a codimension $2$ manifold.
Its curvature properties are not as easy to describe as in the
codimension $1$ case. The propagation no longer occurs along rays, but
instead along two dimensional surfaces in the phase space. Even in the
case when the curvature of the codimension $2$ characteristic set is
``nondegenerate'', the spatial projections of these two dimensional
surfaces through a point must overlap substantially in some
directions. Thus, unlike in the case of operators with real symbols,
there is no hope to obtain uniformly strong kernel decay estimates in
all directions for the parametrix. Instead, there will always be a
lower dimensional set of directions with a weaker kernel decay, and
the geometry of this set can be quite intricate.  Two extreme examples
are $-\Delta_{\R^n} -1 + i D_1$, respectively $D_2 - D_1^2 + i(D_3 -
D_1^2)$.  In the first case the characteristic set is the intersection
of the sphere $|\xi|=1$ with the plane $\xi_1 = 0$; all two
dimensional planes which are normal to it intersect on a line, where
the kernel decays like $|x|^{-1}$ compared to $|x|^{-n/2}$ in the
other directions.  In the second case, the characteristic set is the
curve $\xi_1 \to \gamma(\xi_1) = (\xi_1,\xi_1^2,\xi_1^3)$, which has
nonzero curvature and torsion.  This time the directions of bad decay
of are those perpendicular to both $\dot \gamma$ and $\ddot{\gamma}$;
they are spread on a cone, but in exchange the kernel is not as large
in those directions as in the first case. The generic decay of the low
frequency part of the parametrix is like $|x|^{-3/2}$ compared to
$|x|^{-4/3}$ in the bad directions.

Fortunately we are able to produce a factorization of the parametrix
which allows us to establish the dispersive estimates without having
to study the kernel decay in the bad region. Instead, all we need is
to prove this decay in the good region, where focusing does not occur.
Some related results were independently proved by Dos Santos Ferreira~\cite{DS}.
However, his estimates are somewhat weaker as they are based on the
worst decay rate instead of the generic one.

Another feature of this work is that we seek to obtain 
our estimates with minimal regularity assumptions 
on the coefficients. Precisely, in the case of the 
principally normal operators it turns out that we need
to have $C^2$ coefficients. This gets even better for 
some of the applications to unique continuation, where
we have better spatial localization.

The structure of the article is as follows. In the next section we
state our assumptions and results in a dyadic setting.  The advantage
of doing it this way is that there is more than one interesting case
in which the dyadic results can be applied. In Section~\ref{norm} we
use elliptic arguments to reduce the problems to some cannonical
formulation.

Next we consider a set of increasingly complex problems.  We begin
in Section~\ref{four} with the case of operators with real symbols. First we construct a
wave packet type representation of the fundamental solution, then we
use curvature assumptions to prove pointwise bounds for it.  This in
turn yields the dispersive estimates.

In Section~\ref{five}  we turn to the construction of parametrices in
the symplectic case. This has an operator theory flavour but is also
borrows some ideas from the Littlewood-Paley theory.
Our construction and the $L^2$ type estimates for the parametrix 
apply in an abstract setup for operators 
\[ 
D_t - A(t )+ iB(t) 
\]
where $A(t)$, $B(t)$ are selfadjoint operators in a Hilbert space,
satisfying the fixed time commutator estimate
\[
\|[D_t -A,B] u\| \lesssim \|Bu\|+\|u\|
\] 

In Section~\ref{six} we combine these $L^2$ bounds with the dispersive
estimates for operators with real symbols.  As a first consequence we
obtain the same results for principally normal operators as for the
selfadjoint part.  The second case we consider is the involutive case,
when both the real and the imaginary part of the symbols are of
principal type, with transversal characteristic sets.  The dispersive
estimates follow from bounds for operators with real symbols and the
$L^2$ bounds for the parametrix.  Then we study the degenerate
involutive case, where the imaginary part is still of principal type
but small, say of size $\dd$. In this case we use the wave packet
representation of the fundamental solution in the real case to derive
a similar representation for the parametrix.  This gives pointwise
kernel decay in the good directions, while in the bad decay directions
we fall back on the approach for nondegenerate operators.  We obtain
estimates with a sharp dependence of the constants on $\dd$.

The first application we consider is to local solvability problems.
Principally normal operators are known to be locally solvable
with loss of one derivative. Here we consider instead principally
normal operators with unbounded potentials, and use dipersive
estimates in order to prove similar results.

The last part of the article is devoted to applications to unique
continuation problems; more precisely, we obtain $L^q$ Carleman
estimates for elliptic and parabolic operators with $C^1$
coefficients, and for (non-elliptic) principally normal operators with
$C^2$ coefficients. Our strategy is as follows. On the unit spatial
scale we use only on the $L^2$ Carleman estimates; these allow us to
localize the $L^q$ estimates on a much smaller spatial scale, see
Lemmas~\ref{las},\ref{lad}. On the smaller spatial scale we are able
to use the dispersive estimates for principally normal operators.

\section{The fixed frequency results}
For $\l > 1$ and $j = 0,1, \cdots$ we consider a class of symbols, 
denoted by $S_\lambda^j$, which
satisfy the conditions
\begin{equation}
\begin{array}{c}
|\d^{\aa}_x \d^\bb_\xi a(x,\xi)| \leq c_{\aa,\bb} \l^{-|\bb|}   \quad
|\aa|  \leq j \cr\cr
|\d^{\aa}_x \d^\bb_\xi a(x,\xi)| \leq c_{\aa,\bb}
\l^{\frac{|\aa|-j}2-|\bb|}   \quad |\aa|  \geq j \end{array}
\label{cs} \end{equation}
The parameter $\l$ plays the role of the frequency.  The symbols in
$S^j_\l$ are bounded, and the corresponding operators are $L^2$
bounded.  Both the space of symbols $S_\l^j$ and the corresponding
class of pseudodifferential operators are algebras.  One can see that
the first $j$ derivatives of the symbols behave as for $S_{1,0}$
symbols, while the rest are as for $S_{\frac12,1}$ symbols.  This
suffices in order for the usual calculus of pseudodifferential
operators to apply. In this article we use the Weyl calculus, but this
is not essential.

Later on we use symbols in $S_\lambda^j$ to describe be the frequency $\l$
part of order $0$ pseudodifferential operators. For operators of order
$k \neq 0$ we use the notation
\[
\l^k S_\lambda^j =\{ \l^k a;\ a \in S_\lambda^j\}
\]
and we denote the corresponding class of operators by $\lambda^k
OPS^j_\l$.  All estimates in this paper require only control of a
finite number of derivatives of the symbols. We do not keep track of
the number of derivatives needed. Nevertheless we can and do consider
the symbol spaces and the spaces of operators as Banach spaces. The
statement that a certain operator in $\lambda^k OPS^j_\l$ is bounded
as a linear map from $L^q$ to $L^r$ means that its operator norm is
controlled by the norm of its symbol and possibly by the constants
defined in the assumptions below.  It is however crucial to obtain the
correct dependence of the constants on $\lambda$ and possibly other
parameters.

Given real symbols $p_{re}, p_{im} \in \l \Sl$ and $0\leq \dd \leq 1$
we seek estimates from below for the operator 
\[
p^w(x,D)=\prw(x,D)+ i\delta \piw
\]
 of the form
\[
\|\chi^w u\|_{L^r} \leq c_1\| p^w u\|_{L^q} + c_2\Vert u \Vert_{L^2} 
\]
combining various pairs of exponents $q$ and $r$. The constants $c_1, c_2$
possibly depend on $\l$, $\delta$, finitely many of the constants
$c_{\alpha\beta}$ in \eqref{cs}, and on the geometry of the real
symbols $\pr$ and $\pim $ of the symbol in a quantitave way which will
be described below.  To keep the notation concise we introduce the
following
\begin{definition}
Given a Banach space $X$ and $\rho > 0$ we denote by $\rho X$
the same vector space equipped with the norm
\[
\|u\|_{\rho X} = \rho^{-1} \|u\|_{X}
\]
Thus the unit ball in $\rho X$ is $\rho$ times the unit ball in $X$.
\end{definition}

The main condition which connects $\pr$ and
$\pim$ is a principal normality condition:

{\bf (A1)} The operator $p^w$ 
is principally normal,  i.e.
\begin{equation}\label{bt} 
|\{\pr,\pim\}| \lesssim |p_{re}|+|p_{im}|+1
\end{equation}

As proved in \cite{MR1944027}, this condition guarantees that one has
good $L^2$ estimates from below for $p^w$
even for $p$ in a larger class than $\l \Sl$:

\begin{theorem}
  Let $p=p_{re}+ i p_{im}$ be a symbol which satisfies
\[
|\d^{\aa}_x \d^\bb_\xi p(x,\xi)| \leq c_{\aa,\bb}
\l^{\frac{|\aa|-|\bb|}2}   \quad |\aa|  \geq 2 
\]
 If  (A1) holds then
\begin{equation}
\|\prw u\|_{L^2}^2 + \|\piw u\|_{L^2}^2 \lesssim 
\| p^w u\|_{L^2}^2
+\|u\|_{L^2}^2
\end{equation}
\label{tl2}\end{theorem}

Here we are interested in obtaining $L^r$ estimates, and in order
to do this we need to restrict ourselves to a region of spatial size
$1$ and frequency size $\l$. After a translation we can assume this is
centered at the origin, therefore we set
\[
B_\l = \{(x,\xi);\ |x| < 1, |\xi| < \l\}
\] 
We let $\chi \in S_\l^0$ be  a symbol  which is compactly supported
in $ B_\l$.  For comparison we consider first the case when $p$ is
elliptic, i.e. $p \in \l S_\l^0$, $p^{-1} \in \l^{-1} S_\l^0$.
 \begin{prop} 
a) If $a \in S^0_\lambda$ is supported in $B_\lambda$ and 
$1\le r \le q \le \infty$ then 
\[ 
\Vert a^w u \Vert_{L^q} \lesssim \lambda^{\frac{n}r- \frac{n}q}
\Vert u \Vert_{L^r}.  
\]
b) If $p \in S^0_\l$ is elliptic then
\[
\|\chi^w u\|_{ L^q} \lesssim \lambda^{\frac{n}r- \frac{n}q-1} \|p^wu\|_{ L^{r}}
+\l^{-N} \|u\|_{L^2}, \qquad N > 0
\]
\label{elliptic} \end{prop}

For principally normal operators $p \in \l S^2_\l$ we seek to prove 
estimates of the form
\begin{equation}
\|\chi^w u\|_{\l^{\rho(q)} L^q} \lesssim \|p^wu\|_{L^2+\l^{-\rho(q)} L^{q'}}
+\|u\|_{L^2}
\label{estr}\end{equation}
The exponent $\rho$  is chosen as in the above elliptic estimates,
\[
\rho(q)= \frac{n-1}{2} - \frac{n}{q}
\]
In many examples this relation can also be derived from scaling
considerations.  If (\ref{estr}) holds for some $q$ then it is easy to
see that it must hold for all larger $q$.  However, if $p$ is not
elliptic then (\ref{estr}) cannot hold for $ q =2$.  Our goal is
to find the lowest value of $q$ for which (\ref{estr}) holds.

A dual form of (\ref{estr}) applies to the adjoint operator $ (p^w)^*=
{\bar p}^w$.  However, since the class of operators we work with is
invariant with respect to taking adjoints, it suffices to state it for
$p^w$. We seek a parametrix $K$ for the operator $P$ with
the following properties:
\begin{equation}
\|K f\|_{L^2}+ \|  K  \chi^w  f \|_{\l^{\rho(q)} L^q} \lesssim 
\|f\|_{L^2+\l^{-\rho(q)} L^{q'}}
\label{estrk} 
\end{equation}
\begin{equation}
\|  (p^w K-\id)\chi^w f\|_{L^2} 
\lesssim\|f\|_{L^2+\l^{-\rho(q)} L^{q'}}
\label{estrpk}
\end{equation}

We also want to obtain mixed norm estimates. For this
we split the coordinates $x = (x_1,x')$ and we use the notation
$L^qL^r = L^q_{x_1} L^r_{x'}$.  Then we seek estimates of the form
\begin{equation}
\|\chi^w u\|_{\l^{\rho(q_1,r_1)} L^{q_1}L^{r_1}} \lesssim 
\| p^w u\|_{L^2+\l^{-\rho(q_2,r_2)} L^{q'_2}L^{r'_2}} +\|u\|_{L^2}
\label{estrs}\end{equation}
Again the exponent $\rho$ is determined by comparison with elliptic 
operators,
\[
\rho(q,r)= \frac{n-1}{2} - \frac{1}{q} -\frac{n-1}{r}
\]
Correspondingly, we want a parametrix $K$ with the mapping 
properties
\begin{equation}
\|  K\chi^wf \|_{L^2 \cap\l^{\rho(q_1,r_1)} L^{q_1} L^{r_1}}
\lesssim \|f\|_{L^2+  \l^{-\rho(q_2,r_2)} L^{q'_2} L^{r'_2}}
\label{estrsk}\end{equation}
\begin{equation}
\|  (p^wK-\id)\chi^w f\|_{L^2} 
\lesssim \|f\|_{L^2+ \l^{-\rho(q_2,r_2)} L^{q'_2} L^{r'_2}}
\label{estrspk}\end{equation}
Most of the work in this article is devoted to the construction of 
parametrices for principally normal operators. The parametrix bounds
can be easily connected to the direct estimates using the following 
result.

\begin{prop}
  Assume that there is a parametrix $K$ for ${\bar p}^w$ which
  satisfies \eqref{estrsk} and \eqref{estrspk}. Then \eqref{estrs}
  holds with $(q_1,r_1)$ and $(q_2,r_2)$ interchanged. In particular,
  if there is a parametrix for ${\bar p}^w$ satisfying \eqref{estrk}
  and \eqref{estrpk} then \eqref{estr} holds for $p^w$.
\label{pe}
\end{prop}

\begin{proof}
  
  Let $K$ be the parametrix for ${\bar p}^w$. Given  $g \in
  \l^{-\rho(q_1,r_1)} L^{q'_1} L^{r'_1}$ we consider the decomposition
\[
\chi^w g = g_1 + {\bar p}^w v, \qquad v = K \chi^w g 
\]
The estimates  \eqref{estrsk} and \eqref{estrspk} for ${\bar p}^w$
show that
\[
\|g_1\|_{L^2} + \|v\|_{L^2 \cap  \l^{\rho(q_2,q_2)} L^{q_2}   L^{r_2}}
\lesssim \|g\|_{  \l^{-\rho(q_1,r_1)} L^{q'_1} L^{r'_1}}.
\]
Now we compute
\[
\langle \chi^w u,g\rangle = \langle  u,\chi^w g\rangle = \langle u,
g_1\rangle +  \langle  u,  {\bar p}^w v \rangle = 
\langle u, g_1\rangle +  \langle  p^w u,   v \rangle
\]
Using the previous inequality we estimate
\[
|\langle \chi^w u,g\rangle| \lesssim (\|u\|_{L^2} + \|p^w u\|_{L^2 +
  \l^{-\rho(q_2,q_2)}  L^{q'_2}   L^{r'_2}})
\|g\|_{  \l^{-\rho(q_1,r_1)} L^{q'_1} L^{r'_1}}
\]
This implies (\ref{estrs}) for $p^w$.
\end{proof}

To
obtain such estimates we need some geometric information on the
characteristic sets of $\pr$ and $\pim$. We denote
\[
\charp =  \text{char } \pr \cap B_\l, \qquad \charq =  
 \text{char } \pim \cap B_\l
\]
\[
\charpq =  \text{char } p \cap B_\l.
\]
Our first results use only the geometric information
for $\pr$, which is introduced next.

{\bf (A2)} $\pr$ is of real principal type, i.e.
\begin{equation}
|\nabla_\xi \pr(x,\xi)| \gtrsim 1 \qquad \text{ in } \charp
\end{equation} 

This implies that for each $x$ the fiber $\charp_x$ of the
characteristic set $\charp$ is a smooth codimension one hypersurface.
Since $p \in \Sl$, the second fundamental form of $\charp$ has size
$O(\l^{-1})$. Then we impose the following curvature condition on
$\charp$.

{\bf (A3)} The characteristic set $\charp$ has $n-1-k$ nonvanishing
curvatures, i.e. for each $x$ the second fundamental form of
$\charp_x$ has  rank at least $n-1-k$. More precisely there exists 
a $n-1-k$ minor $M$ of the second fundamental form of $\charp_x$ with
\begin{equation}
|\det M| \gtrsim \lambda^{k-n+1}.  
\end{equation}

Here $0 \leq k \leq n-2$, but the most interesting cases for
applications are $k=0,1$. The first corresponds to the Schr\"odinger
equation, while the second arises in the study of the wave equation.

Finally, we prove mixed norm estimates, and
this requires a choice of coordinates.

{\bf (A4)} The level sets of $x_1$ are noncharacteristic for $\pr$, i.e.
\begin{equation}
|\d_{\xi_1}\pr| \gtrsim 1 \qquad \text{in } \charp
\end{equation}

We begin with a result for operators with real symbols:

\begin{theorem}
a) Let $p  \in \l \Sl$ be a real symbol  satisfying  {\bf (A2-3)}. 
Let
\begin{equation}
 q = \frac{2(n+1-k)}{n-1-k}, \qquad \rho(q) = \frac{n-1+k}{2(n+1-k)}
\label{r}\end{equation}
Then \eqref{estr} holds for $p^w$, and there is a parametrix $K$ for $p^w$
which satisfies \eqref{estrk} and \eqref{estrpk}.

b) If in addition the coordinates are chosen so that {\bf (A4)} holds
 and  $(q_1,r_1)$, $(q_2,r_2)$  satisfy
\begin{equation}
\frac{2}{q}+\frac{n-1-k}{r} = \frac{n-1-k}{2}, 
\quad 2 \leq q,r \leq \infty, \quad (q,r) \neq (2,\infty)
\label{rs}\end{equation}
then \eqref{estrs} holds for $p^w$, and there is a parametrix $K$ for $p^w$
which satisfies \eqref{estrsk} and \eqref{estrspk}.
\label{tp}\end{theorem}

As it turns out, the conclusion of the previous theorem often remains valid
if we add a small imaginary component to the above real symbol.
We  subject such perturbations to a slightly strenghtened form of the
 principal normality condition {\bf (A1)}, namely
 
 {\bf (A1)'} There exist symbols $r_{1} \in \lambda S_\l^2 $, $r_2,
 r_3 \in S_\l^1$ and $r_4 \in S_\l^0$ such that
\begin{equation}
 \{ \pr,\pim\} =  r_1 + r_2 \pr + r_3 \pim + r_4, \qquad |r_1|
 \lesssim |p|
\end{equation}

\begin{theorem}
  a) Let $p \in \l \Sl$ be a symbol satisfying {\bf (A1)'}, {\bf (A2)}
  and {\bf (A3)}.  Let $q$ be as in \eqref{r}. Then \eqref{estr} holds
  for $p^w$, and there is a parametrix $K$ for $p^w$ which satisfies
  \eqref{estrk} and \eqref{estrpk}.
  
  b) Assume in addition that the coordinates are chosen so that {\bf
    (A4)} also holds and at least one of the following two conditions
  applies:
\[
{\bf (B1)}  \quad \qquad \qquad \qquad  |d_\xi p_{re} \wedge d_\xi p_{im}| \ll  1
\quad  \qquad  \qquad \text{in }\ \charpq
\]
\[
{\bf (B2)} \qquad \qquad |\d_{\xi_1} p_{im}| \ll 1,\quad
|\d_\xi \d_{\xi_1} p_{im}| \ll \l^{-1}\qquad  \text{in }\ \charpq
\]
If $(q_1,r_1)$, $(q_2,r_2)$ satisfy \eqref{rs} and $q_1,q_2 > 2$ then
\eqref{estrs} holds for $p^w$, and there is a parametrix $K$ for $p^w$
which satisfies \eqref{estrsk} and \eqref{estrspk}.
\label{tpq1}\end{theorem}

We had hoped to prove the result with {\bf (A1)'} replaced by {\bf
(A1)} but this seems to cause certain difficulties in the proof. 
It would be interesting to know if this can be done.

We believe that part (b) should hold in general, without any of the
assumptions {\bf (B1)} and {\bf (B2)}. Note that {\bf (B1)} is not so harmful, because
if it does not hold then {\bf (A2)'} below must hold so we can place
ourselves in the setup of the stronger Theorem~\ref{tpq} below.
Unfortunately the case (b) of Theorem~\ref{tpq} contains
the assumption\footnote{This is likely to be unnecessary, see the discussion after
  Theorem~\ref{tpq}.} {\bf (A5)'}, which rules out some of the
examples we wish to consider later on. To at least partially
compensate for that we have also added the case {\bf (B2)} to the theorem.

Another possible improvement to Theorem~\ref{tpq1} would be to also
prove it when $q_1=q_2=2$. It may be possible to modify our argument
to allow one of $q_1$, $q_2$ to equal $2$, but we do not know how to
allow both of them to be $2$.

Our next result uses geometric information for both $\pr$ and
$\pim$; this allows for an improved range of indices in the estimates.
Begin with

{\bf (A2)'} $\pr$ and $\pim$ are of real principal type and their
characteristic sets are transversal, i.e.
\begin{equation}
|d_\xi \pr \wedge d_\xi \pim | \gtrsim 1
\quad \text{ in }  \charpq
\end{equation} 

This implies that for each $x \in \R^n$ the $x$ section $\charxpq$ of
the characteristic set of  $p$ is a smooth codimension
two submanifold of $\R^n$.  At each $\xi \in \charxpq$ the two
dimensional normal space $N \charxpq$ is generated by $\d_\xi \pr (x,\xi)$
and $\d_\xi \pim (x,\xi)$.  Its second fundamental form $S_{x,\xi}$ maps $N
\charxpq \times T\charxpq$ into $T\charxpq$. If we consider it as a
quadratic form in $ T\charxpq$ depending on the parameter $\nu \in N
\charxpq$, its rank may well depend on $\nu$.  In particular is not
possible to have $S_{x,\xi}(\nu)$ nondegenerate for all $\nu \in N
\charxpq$.  However if it is nondegenerate for some $\nu$ then it must
be nondegenerate for all $\nu$ except for at most $n-2$ values.  In
this case the directions $\nu$ for which $S_{x,\xi}(\nu)$ is
degenerate are precisely the directions in which the kernel of the
parametrix for $p^w$ has less decay. Our contention is that for the
purpose of proving dispersive estimates we can neglect the bad
directions and instead use only an assumption on the generic
behaviour. Consequently, the curvature condition we impose on
$\charpq$ is as follows:

{\bf (A3)'} The characteristic set $\charpq$ has $n-2-k$ nonvanishing
curvatures, i.e. for each $(x,\xi) \in  \charpq$ there is $\nu  \in
N_\xi \charxpq$ so that the second fundamental form $S_{x,\xi}(\nu)$
has rank (at least) $n-2-k$. More precisely we assume that there is
a $n-2-k$ minor $M$ of $S_{x,\xi}(\nu)$ with 
\begin{equation}
 |\det M|\gtrsim \l^{k-n+2}. 
\label{n-2} 
\end{equation}

In order to obtain mixed norm estimates we need to impose some
restriction on how coordinates are chosen. Thus we add the replacement
of {\bf (A4)}, namely

{\bf (A4)'} The level sets of $x_1$ are noncharacteristic for $p^w$,
i.e.
\begin{equation}
|p_{\xi_1}| \gtrsim 1 \qquad \text{in } \charpq
\end{equation}

For technical resons we also invoke a last condition which essentially
says that there are no bad directions for the second fundamental
form which are tangent to the level sets of $x_1$. More precisely, the
condition {\bf (A4)'} guarantees that $N_\xi \charxpq$ is transversal to the
planes $x_1=const$. Hence for each $(x,\xi) \in \charpq$ there is an
unique direction $\nu_0 \in N_\xi \charxpq$ which is tangent
to $x_1=const$. 

{\bf (A5)'} There exists a $n-2-k$ minor $M$ of  $S_{x,\xi}(\nu_0)$
which satisfies \eqref{n-2}. 

\begin{theorem}
  a) Let $p \in \l \Sl$ be a  symbol satisfying {\bf (A1),
    (A2)'} and {\bf (A3)'}. Let $q$ satisfy
\begin{equation}
 q=  \frac{2(n+2-k)}{n-k} , \qquad \rho(q)= \frac{n-2+k}{2(n+2-k)}
\label{r1}\end{equation}
Then \eqref{estr} holds for $p^w$, and there is a parametrix $K$ for
$p^w$ which satisfies \eqref{estrk} and \eqref{estrpk}.
  
b) Assume that in addition the coordinates are chosen so that {\bf
  (A4)'}, {\bf (A5)'} hold and $(q_1,r_1)$, $(q_2,r_2)$ satisfy
\begin{equation}
\frac{2}{q}+\frac{n-k}{r} = \frac{n-k}{2}\quad 2 \leq q,r \leq \infty,
\quad (q,r) \neq (2,\infty)
\label{rs1}\end{equation}
with $q_1$ and $q_2$ not both equal to $2$.
Then \eqref{estrs} holds for $p^w$, and there is a parametrix $K$ for
$p^w$ which satisfies \eqref{estrsk} and \eqref{estrspk}. 
\label{tpq}\end{theorem}

We expect Theorem~\ref{tpq} to be valid without the condition {\bf
  (A5)'}, which insures good kernel decay for the parametrix near
$x_1$ slices.  The primary obstruction in this result should come from
low kernel decay in the $x_1$ direction, not along $x_1$ slices. 
We were able to prove partial results with logarithmic losses, but 
not fully remove {\bf(A5)'}.

Finally, we also consider a degenerate case. For $0 < \dd < 1$
we consider the operator $p=\pr+i \dd \pim$  and seek an estimate similar to
\eqref{tpq}, but with the correct control of the constants as a function of
$\dd$. 

As $\dd$ approaches $0$ one expects the parametrix for $p^w$
to concentrate closer to the Hamilton flow of $\pr$. Then it is
natural to assume that the direction of the  Hamilton flow of $\pr$
is not one of the bad directions:

{\bf (A6)'} At each point in $\charpq$  the second fundamental form of
$\charp$ restricted to $T\charpq$  has rank at least
$n-2-k$, i.e. it has a $n-2-k$ minor $M$  which satisfies \eqref{n-2}.

\begin{theorem}
  a) Let $\prw, \piw \in \Sl$ be real symbols satisfying {\bf (A1)},
  {\bf(A2)'}, {\bf (A3)} with $k$ replaced by $k+1$, {\bf (A3)'} and
  {\bf (A6)'}. Let $q$ satisfy \eqref{r1}.  Then \eqref{estr} holds
  for $p^w$ with $\l^{\rho(q)}$ replaced by
  $\dd^{-\frac{1}{n+2-k}}\l^{\rho(q)}$; also there is a parametrix $K$
  for $p^w$ which satisfies \eqref{estrk} and \eqref{estrpk} with a
  similar substitution.

    b) If in addition the coordinates are chosen so that {\bf
      (A4)},{\bf (A5)'} are satisfied and $(q_1,r_1)$, $(q_2,r_2)$ are as in
    Theorem~\ref{tpq}  then \eqref{estrs} holds for $p^w$ with
    $\l^{\rho(q,s))}$ replaced by
    $\dd^{-\frac{2}{q(n-k)}}\l^{\rho(q,s)}$; also there is a
    parametrix $K$ for $p^w$ which satisfies \eqref{estrsk} and
    \eqref{estrspk} with a similar substitution.
\label{tpdq}\end{theorem}

\section{Cannonical models}
\label{norm}

Our  aim  is to reduce the operator $p^w$ to a
cannonical form. We decompose the coordinates in the form 
$x = (x_1,x^\prime)$ and $\xi = (\xi_1,\xi^\prime)$ with $x_1,\xi_1 \in \R$ 
and $x^\prime, \xi^\prime \in \R^{n-1}$. 

\begin{definition}
We say that the pair of symbols $p_{re},p_{im} \in \l S_\l$ are in cannonical 
form if there are real symbols $a,b \in \l S_\l$ so that 
\[
p_{re}(x,\xi) = \xi_1 + a(x,\xi'), \qquad p_{im}(x,\xi) = b(x,\xi')
\]
\end{definition}

The main result of this section asserts that it suffices to prove 
the results in the paper when $p,q$ are in the cannonical form.
If $p,q$ are in cannonical form then we will prove the estimate 
for a  cutoff symbol $\chi$ of the form  
\[
\chi = \chi(x,\xi')
\]
supported  in 
\[ B^\prime_\lambda = \{ (x,\xi^\prime): |x| \le 1, |\xi'| \le \lambda\} \]
This strengthens the estimates.

\begin{prop}
  a) Assume that the results in Theorems~\ref{tp},\ref{tpq},\ref{tpdq}
  hold for $p_{re},p_{im}$ in cannonical form. Then they hold in
  general.

b) Assume that the result in Theorem~\ref{tpq1} 
holds  for $p_{re},p_{im}$ in cannonical form
with the weaker hypothesis that the curvature
assumption {\bf (A3)} applies to $\xi_1 + a(x,\xi')+ \alpha b(x,\xi)$
for some real $\alpha$ (instead of $\xi_1 + a(x,\xi')$).
Then it holds in general.
\label{cann}\end{prop}

\begin{proof} We prove a series of lemmas, which imply the assertion. 
First we note that operators localized in $B_\l$ are bounded in 
all $L^q L^r$ spaces.

\begin{lemma}
Let $\eta \in S_\l^0$ be supported in $B_\l$. Then $\eta^w$ is bounded
in  $L^q L^r$ for all $1 \leq q,r \leq \infty$.
\label{intk}\end{lemma}
This follows from the translation invariant kernel bound
\[  |k(x,y) | \le c_N \lambda^{n}(1 +\lambda |x-y|)^{-N}   \] 
where $k(x,y)$ is the kernel of $\eta^w$.  As an immediate application
we show that one can localize better the output of the parametrix $K$.

\begin{lemma} Suppose that $K$ is a parametrix satisfying
  \eqref{estrsk} and \eqref{estrspk}. Let $\tilde \chi$ be supported
  in $B_{\lambda}$ and identically $1$ in the support of $\chi$. Then
  $\tilde K = \tilde \chi^w K $ satisfies \eqref{estrsk} and
  \eqref{estrspk}.
\label{co} \end{lemma} 
\begin{proof}
Lemma~\ref{intk}  implies the estimate \eqref{estrsk} for $\tilde K$. 
For the error estimate we write
\[ 
(p^w \tilde \chi^w K - \id ) \chi^w = \tilde \chi^w ( p^w K - \id) \chi^w 
+ (\id- \tilde \chi^w) \chi^w + [p^w,\tilde \chi^w] K \chi^w 
\]
In the first term we use the $L^2$ boundedness of $\tilde \chi^w$.
The second is negligible because $\chi$ and $\id-\tilde \chi$ have
disjoint supports. The commutator in the third term is in $OPS^0_\l$
and therefore $L^2$ bounded.
\end{proof} 

Next we prove that the estimates \eqref{estrsk} and \eqref{estrspk}
do not change if we multiply $p$ by  a zero order elliptic symbol.  

 \begin{lemma}
   Let $e \in \Sl$ be an elliptic symbol. Then the conclusion of
   either of the Theorems~\ref{tp},\ref{tpq1},\ref{tpq},\ref{tpdq}
   holds for $p$ if and only if it holds for $\tilde p = e p$.
\label{ellmul}\end{lemma}
\begin{proof}
  Let $K$ be a parametrix of $p^w$ which satisfies \eqref{estrsk} and
  \eqref{estrspk}.  By Lemma~\ref{co} we may replace $K$ by $\tilde
  \chi^w K$. By Lemma~\ref{pe} is suffices to show that $ (e^{-1})^w
  \tilde \chi^w K $ is a parametrix for $\tilde p^w$ which again
  satisfies \eqref{estrsk} and \eqref{estrspk}. The inequality
  \eqref{estrsk} holds for $ (e^{-1})^w \tilde \chi^w K$ because, by
  Lemma~\ref{intk}, the operator $ (e^{-1})^w \tilde \chi$ is bounded
  on $L^q(L^r)$. The inequality \eqref{estrspk}  follows from
\[ 
e^w p^w (e^{-1})^w \tilde \chi^w K - \id = p^w \tilde \chi^w K - \id + 
 (e^w p^w (e^{-1})^w  - p^w)\tilde \chi^w K 
\]  
where $  (e^w p^w (e^{-1})^w  - p^w) \in OPS^0_\l $ is $L^2$ bounded.
 \end{proof} 
 
 We now consider the issue of localization. The next lemma asserts
 that the estimates in
 Theorems~\ref{tp},\ref{tpq1},\ref{tpq},\ref{tpdq} are
 microlocalizable. This allows us to carry out the reduction to the
 cannonical form locally in the phase space.

\begin{lemma}
  Let $0 < \e< 1$. Assume that the conclusion of either of Theorems
  ~\ref{tp},\ref{tpq1},\ref{tpq},\ref{tpdq} holds for $\chi$ compactly
  supported in any ball $\e B_{\l}$ contained in $B_\l$. Then the same
  holds for $\chi$ supported in $B_\l$.
\label{localize} \end{lemma}

\begin{proof}
We consider a finite covering
\[
\text{supp } \chi \subset  \bigcup_j \e B_{\l}^j
\]
and a corresponding partition of unity
\[
1 = \sum \chi_j \text{ in } \text{supp } \chi, \qquad \text{supp } \chi_j \subset
\e B_{\l}^j.
\]
We denote by $K_j$ the parametrix for $p^w$ associated to
$\e B_{\l}^j$.  We define a parametrix $K$ for $p^w$ in $B_\l$ by
\[
K = \sum  K_j \chi_j^w.
\]
Then (\ref{estrsk})  is verified directly, while for (\ref{estrspk}) we 
compute
\[
(p^w K - id)\chi^w = \sum  (p^w K_j - \id)\chi_j^w \chi^w +
(\id- \sum \chi_j^w) \chi^w.
\]
The first term is estimated using the hypothesis for $\e B_{\l}^j$,
and the second  is a smoothing operator since the supports 
of $\chi$ and $(1- \sum \chi_j)$ are separated. 
\end{proof}

 Our next concern is the choice of coordinates. These are 
uniquely determined in part (b) of the theorems, but we have 
a choice to make in part (a). The next lemma asserts
that we can always choose coordinates so that  part (a)
is a special case of part (b).

\begin{lemma}
Assume that Theorems~\ref{tp},\ref{tpq1},\ref{tpq},\ref{tpdq} are true
under the additional hypothesis {\bf (A4)}, and provided that all the
assumptions in part (b) also hold for part (a). Then
Theorems~\ref{tp},\ref{tpq1},\ref{tpq},\ref{tpdq} are true as stated.
\end{lemma}

\begin{proof}
  There is nothing to do for part (b) of
  Theorems~\ref{tp},\ref{tpq},\ref{tpdq}.  For Theorem~\ref{tpq} (b)
  we locally get {\bf (A4)} from {\bf (A4)'} if we multiply $p$ by a suitably
  chosen complex number.

Consider now part (a) of
Theorems~\ref{tp},\ref{tpq1},\ref{tpq},\ref{tpdq}.
In Theorem~\ref{tp}(a) we take some $(x_0,\xi_0) \in \charp$. By {\bf (A1)}
$\d_{\xi} p_{re}(x_0,\xi_0) \neq 0$, therefore we can choose
coordinates so that $\d_{\xi_1} p_{re}(x_0,\xi_0) \neq 0$. Then
the same must be true in an $\e B_\l$  neighbourhood of $(x_0,\xi_0)$,
i.e. {\bf (A4)} holds there.
 
In Theorem~\ref{tpq1}(a) we we take $(x_0,\xi_0) \in \charp$ and
consider three cases.

\begin{itemize}
\item  If $|p_{im}(x_0,\xi_0)| \gtrsim \l$ then we are
in the elliptic region and all our estimates are straightforward. 
\item If $|p_{im}(x_0,\xi_0)| \ll \l$ and $|\d_\xi p_{im}(x_0,\xi_0)|
  \ll 1$ then {\bf (B1)} holds in a neighbourhood, therefore we can
  proceed as in the case of Theorem~\ref{tp}.
\item If $|p_{im}(x_0,\xi_0)| \ll 1$ and $|\d_\xi p_{im}(x_0,\xi_0)|
  \gtrsim 1$ then $\xi_0$ is close to $\charpq_{x_0}$, therefore
we can move it slightly and assume that $(x_0,\xi_0) \in \charpq$.
 \begin{itemize} \item If {\bf (B1)} holds there then we proceed
as in the case of Theorem~\ref{tp}. 
\item If {\bf (B1)} does not hold at $(x_0,\xi_0)$ then {\bf (A2)'}
  holds there.  If $\charp_{x_0}$ has $n-k-1$ nonvanishing curvatures
  at $\xi_0$ then $\charpq_{x_0}$ has at least $n-k-3$ nonvanishing
  curvatures at $\xi_0$.  Hence we can use Theorem~\ref{tpq}(a) to get
  a result which is at least as good as Theorem~\ref{tpq1}(a).
\end{itemize}
\end{itemize}

In the case of  Theorem~\ref{tpq}(a) we note that {\bf (A5)'} includes {\bf
  (A4)'}.  Given $(x_0,\xi_0) \in \charpq$ we need to choose
coordinates near it so that {\bf (A5)'} holds.  This only
requires that $dx^1$ is not orthogonal to the set of finitely many
directions of bad decay associated to $(x_0,\xi_0)$, so it can always
be achieved.

For Theorem~\ref{tpdq}(a) we argue as above, with the only difference
that now {\bf (A6)'} requires {\bf (A4)}, namely that $dx^1$ not be
orthogonal to $\d_\xi p_{re}(x_0,\xi_0)$.

\end{proof}

We now proceed to the main result of the section, which asserts that
we can always choose an elliptic symbol $e$  so that
$e p$ is in the cannonical form.

\begin{lemma}
(a) Let $p \in \l\Sl$ be a real symbol which satisfies {\bf (A4)}. 
Then near each $(x,\xi) \in \charp$ there is
a real elliptic symbol $e \in \Sl$ and a real symbol
$a \in \l \Sl$
 so that
\[
e(x,\xi) p(x,\xi) =  \xi_1+a(x,\xi')
\]
(bc) Let $p=\pr+i\pim \in \l\Sl$ be a symbol which satisfies {\bf (A4)}.  Then
near each $(x,\xi) \in \charp$ there is an elliptic symbol $e \in \Sl$
and real symbols $a,b \in \l \Sl$ so that 
\[
e(x,\xi) p(x,\xi) =  \xi_1+a(x,\xi')+ i b(x,\xi').
\]
(d) Let $p = \pr+i\dd \pim$, with $\pr,\pim \in \l \Sl$ symbols which
satisfy {\bf (A4)} and $\dd \in [0,1]$.  Then near each $(x,\xi) \in
\charp$ there is an elliptic symbol $e \in \Sl$ with $\Im e \in \dd
\Sl$ and real symbols $a,b \in \l \Sl$ so that
\[
e(x,\xi) p(x,\xi) = \xi_1+a(x,\xi') + i \dd b(x,\xi') 
\]
with the bounds for $e,a,b$ independent of $\dd$.
\label{fsn} \end{lemma}

\begin{proof}
(a) Since $p$ is real and $\d_{\xi_1} p \approx 1$ it follows that the zero set of $p$
can be expressed as
\[
\{ p(x,\xi) = 0 \} = \{ \xi_1 + a_0(x,\xi') = 0 \}   \qquad a \in \Sl
\]
We make a change of variable $\xi_1 \to \xi_1+a_0(x,\xi')$
which leaves unchanged the classes of symbols we work with.
This reduces the problem to the case when 
\[
\{ p(x,\xi) = 0 \} = \{ \xi_1 = 0\}
\]
and  we can take
\[
e(x,\xi) = \frac{\xi_1}{p(x,\xi)}.
\]
It is easy to verify that $e$ has the desired regularity.

(bc) If $\pim \neq 0$ then we begin by reducing the problem as before to the
case when $\pr = \xi_1$.  We want to find an elliptic symbol
$e \in \Sl$ and real symbols $a,b \in \l \Sl$  such that
\[ 
e(\xi_1+i \pim) = \xi_1 + a(x,\xi^\prime) + i b (x,\xi^\prime) 
\]
First we produce a formal series with this property:

\begin{lemma}
Let 
\[
 \pim(x,\xi_1,\xi') \approx \sum_{k \geq 0} q_k(x,\xi') \xi_1^k, \qquad q_k \in
 \l^{1-k} \Sl
\]
be the formal Taylor series for $\pim$ at $\xi_1=0$.  Then there are
formal series
\[
\sum_{k,l \geq 0} e_{k,l}(x,\xi') \xi_1^k q_0^l, \qquad  
\sum_{l\geq 1} a_{l}(x,\xi')  q_0^l,
\qquad \sum_{l \geq 1} b_{l}(x,\xi')  q_0^l
\]
with coefficients $e_{k,l} \in \l^{-k-l}\Sl$, $a_l,b_l \in
\l^{1-l}\Sl$ whose partial sums
\[
e^N = \sum_{k+l \leq N} e_{k,l}(x,\xi') \xi_1^k q_0^l, \quad a^N =
\sum_{l \leq N} a_{l}(x,\xi') q_0^l,
\quad b^N=\sum_{l \leq N} b_{l}(x,\xi')  q_0^l
\]
satisfy
\[ 
e^{N-1}(\xi_1+i \pim(x,\xi)) = \xi_1 + a^N(x,\xi) + i b^N(x,\xi) 
+ O(\l^{-N}(|\xi_1|+|q_0|)^{N+1})
\]
\label{ssn} \end{lemma}
\begin{proof}
The coefficients $e_{k,l}$, $a_l$, $b_l$ are uniquely determined
inductively. We begin with $N=1$, where we must have 
\[
e_{00} (\xi_1+ iq_0+ iq_1 \xi_1) = \xi_1 + a_1 q_0  + i b_1 q_0. 
\]
We check this first at points where $q_0=0$ and second  where  $\xi_1=0$. This 
yields
\[
e_{00} = (1+i q_1)^{-1}, \qquad a_1 = - \Im (1+i q_1)^{-1}, \qquad
b_1 = \Re  (1+i q_1)^{-1}
\]
For the induction step we must have
\[
\sum_{k+l=N-1}\! e_{k,l} \xi_1^k q_0^l (\xi_1(1+iq_1)+i q_0) =
 (a_N+ib_N) q_0^N - \!\!\! \sum_{k+l < N-1}\! e_{k,l} q_{N-k-l} \xi_1^{N-l} q_0^l 
\]
Then the coefficients $e_{k,l}$ are obtained by polynomially dividing
the sum on the right by $\xi_1(1+iq_1) +i q_0$ as polynomials in
$\xi_1$; finally, and $(a_N+ib_N) q_0^N$ is the remainder.
\end{proof}

The second step is to find smooth functions which match the formal
series up to any order. This is a classical argument, see e.g.
H\"ormander~\cite{MR87d:35002a}, Proposition 18.1.3.

\begin{lemma}
There are symbols $e\in  \Sl$, $a,b \in \l \Sl$ so that
\[
e-e^N = O(\l^{-N-1}(|\xi_1|+|q_0|)^{N+1}), \quad a - a^N, b-b^N = 
O(\l^{-N} |q_0|^{N+1})
\]
\label{ellnorm} \end{lemma}

Now we continue the proof of Lemma~\ref{fsn}.  Combining
Lemma~\ref{ssn} and Lemma~\ref{ellnorm}, we obtain
\[
e(\xi_1+ i q(x,\xi)) = \xi_1 + a(x,\xi') + i b(x,\xi') + r(x,
\xi',\xi_1,q_0)
\]
where the remainder term $r$ vanishes of infinite order
at $\xi_1 = 0, q_0 = 0$. Finally we elliminate $r$ with the
substitution 
\[
e:= e+ r(x, \xi',\xi_1,q_0)(\xi_1+ i q(x,\xi))^{-1}
\]
It is clear that the second right hand side term 
is a smooth symbol.

(d) This follows simply by replacing $q_0$ by  $\delta q_0$ in the 
argument above. 

\end{proof}

To conclude the proof of Proposition~\ref{cann} it remain to study 
how  the hypothesis of our theorems is modified by the
multiplication with the elliptic symbol constructed in
Lemma~\ref{fsn}.  We begin with

\begin{lemma}
a) The hypotheses of Theorem~\ref{tp} remain unchanged if we multiply
$p$ by any real elliptic symbol $e\in \Sl$. 

(b) The hypotheses of Theorem~\ref{tpq1} remain unchanged if we
multiply $p$ by any real elliptic symbol $e \in \Sl $.

c) The hypotheses of Theorem~\ref{tpq} remain unchanged if we multiply
$p$ by any elliptic symbol $e \in \Sl$.

d) The hypotheses of Theorem~\ref{tpdq} remain unchanged if we
multiply $p$ by any  symbol $e \in \Sl$ with $\Re\, e$ elliptic
and $\Im\, e \in \dd \Sl$.
\label{reduction}\end{lemma}

The proof of the lemma is straightforward. It completes the proof of
Proposition~\ref{cann} (a).  However, Proposition~\ref{cann} (b),
which refers to Theorem~\ref{tpq1}, requires a more detailed
discussion because the symbol $e$ does not necessarily satisfy the
condition in part (b) of Lemma~\ref{reduction}.

To understand what happens we consider each of the steps
in the proof of Lemma~\ref{fsn}(bc). By Lemma~\ref{reduction}(b)
the initial multiplication by an elliptic symbol changes nothing.

We fix some $(x_0,\xi_0)$ on $\charp$.  The change of variable
$\xi_1+a_0(x,\xi') \to \xi_1$ does not affect either of {\bf (B1)} or
{\bf (B2)}.  If we have $|q_0(x_0,\xi'_0)| \gtrsim \l$ then $p$ is
elliptic at $x_0,\xi_0$ and all our estimates are straightforward.
Hence in what follows we assume that $|q_0(x_0,\xi'_0)| \ll \l $.
We consider two cases:

{\em Case 1.} Suppose {\bf (B1)} holds near $(x_0,\xi_0)$. Then we
must have
\begin{equation}
|q_0(x_0,\xi'_0)| \ll \l, \qquad 
|\d_\xi q_0(x_0,\xi'_0)| \ll 1.
\label{qdq}\end{equation}
The symbols $a(x,\xi')$, $b(x,\xi')$ obtained in Lemma ~\ref{fsn}(bc)
can be written in the form
\begin{equation}
a(x,\xi') = a_0(x,\xi') + a_1( x,\xi') q_0( x,\xi') + f(x,\xi')
q_0^2( x,\xi'), \quad f \in \l^{-1} \Sl
\label{aa1}\end{equation}
\begin{equation}
\label{bb1} b(x,\xi) = b_1( x,\xi') q_0( x,\xi')  + g(x,\xi')
q_0^2( x,\xi'), \qquad g \in \l^{-1} \Sl
\end{equation}
Here $b_1 \in \Sl$ is elliptic. Since the symbol $a_1( x,\xi') \in \Sl$
is not necessarily small, near $(x_0,\xi_0)$ we write
\[
a(x,\xi') = a_0(x,\xi') + \frac{a_1( x_0,\xi'_0)} {b_1( x_0,\xi'_0)} b( x,\xi') + \tilde a(x,\xi')
\]
where, due to \eqref{qdq}, the remainder $\tilde a(x,\xi')$ satisfies
\begin{equation}
|\d_\xi \tilde a(x_0,\xi'_0)|\ll 1, \qquad |\d_\xi^2 \tilde a(x_0,\xi'_0)| \ll \l^{-1}  
\label{tta}\end{equation}
Since $\xi_1+a_0(x,\xi')$ satisfies the curvature condition {\bf (A2)} this
implies that so must $\xi_1+ a(x,\xi') - \alpha  b(
x,\xi')$ near $(x_0,\xi_0)$, where $\alpha = \frac{a_1( x_0,\xi'_0)}
{b_1( x_0,\xi'_0)}$.

{\em Case 2.}  Suppose {\bf (B1)} does not hold near
$ (x_0,\xi_0)$.  Then we must have
\[
|q_0(x_0,\xi_0)| \ll \l \qquad 
|\d_\xi q_0(x_0,\xi_0)| \gtrsim  1.
\]
Hence we can shift $\xi_0$ slightly to arrive at the case when
\[
|q_0(x_0,\xi_0)| =0  \qquad 
|\d_\xi q_0(x_0,\xi_0)| \gtrsim  1.
\]
Since {\bf (B1)} does not hold, {\bf (B2)} must hold at $(x_0,\xi_0)$.
This implies that
\[
|q_1(x_0,\xi_0) \ll 1, \qquad |\d_\xi q_1(x_0,\xi_0)| \ll \l^{-1}
\]
Hence the same must hold for the symbol $a_1 = \Im (1+iq_1)^{-1}$.
Going back to \eqref{aa1}, we obtain a relation of the form
\[
a(x,\xi') = a_0(x,\xi') + \tilde a(x,\xi')
\]
where $\tilde a$ satisfies \eqref{tta}. Then $\xi_1+a(x,\xi')$
 satisfies the curvature condition {\bf (A2)} near $(x_0,\xi_0)$.
This concludes the proof of Proposition~\ref{cann}.
\end{proof}

\section{The real case }
\label{four}

\subsection{The parametrix construction}

In this section we denote $x_1$ by $t$, set $d=n-1$ and redenote
$x^\prime$ by $x\in \R^d$.  We shall construct a phase space
representation of the fundamental solution for the initial value
problem
\begin{equation}
(D_t + a^w(t,x,D)) u = f  \qquad u(0) = u_0
\label{ivp}\end{equation}
where $t \in [0,1]$ and $x\in \mathbb{R}^d$. If the symbol $a$ is real
then $a^w$ is selfadjoint, therefore it generates an isometric evolution
operator $S(t,s)_{t,s \in [0,1]}$ in $L^2(\R^d)$.

To keep the argument simple we work in a normalized setup where all
scales are of order $1$.  Thus we assume that the symbol $a$ is
measurable in $t$ and that it satisfies the bounds
\begin{equation}
  |\d_x^\alpha \d_\xi^\beta a (t,x,\xi)|  \leq c_{\aa,\bb} \qquad
|\aa|+|\bb| \geq 2
\label{a11}\end{equation}

We first outline a simple parametrix construction based on the FBI
transform, following ideas in \cite{MR1833146}.  There the equation
is conjugated with respect to the FBI transform and replaced by a
simpler transport equation in the phase space.

The FBI transform\footnote{Often one adds the factor $e^{\frac12
    \xi^2}$ in the formula; this would generate some obvious changes
  in what follows.} is an isometry from $L^2(\R^d)$ to $L^2(\R^{2d})$
which is defined by
\[ 
Tu(x,\xi) =2^{-\frac{d}2} 
\pi^{- \frac{3d}4} \int e^{- i \xi(x-y)-\frac12(x-y)^2} f(y) dy
\]
We approximate the  conjugated operator by
\[
\tilde{A} = a(x,\xi) +    a_\xi ( \frac1i\d_x  - \xi) - \frac1i a_x \d_\xi
\]
and we use the error estimate in  \cite{MR1944027}, Theorem 6:

\begin{lemma} \label{error}
Assume that the symbol $a$ satisfies \eqref{a11}. Then we have
\[
\|T a^w - \tilde A T\|_{L^2 \to L^2} \lesssim 1
\]
and the dual estimate
\[
\|T^* \tilde A - a^w T^*\|_{L^2 \to L^2} \lesssim 1
\]
\end{lemma}    
The operator $\tilde A$ is selfadjoint, therefore it generates
an isometric evolution operator  $\tilde S(t,s)$ in $L^2(\R^{2n})$.
Then a natural choice for a forward  parametrix is the operator
\[
K(t,s) = 1_{t \geq s}  T^* \tilde S(t,s) T
\]
Given the above error estimates, it is straightforward to prove that
this provides a good approximate solution in the $L^2$ sense:

\begin{prop}
Assume that the symbol $a$ satisfies \eqref{a11}. Then the operator
$K(t,s)$ satisfies
\[
\|K(t,s)\|_{L^2 \to L^2} \leq 1
\]
\[
\|(D_t + a^w) K(t,s)\|_{L^2 \to L^2} \lesssim 1
\]
\[ 
\lim_{t \to s+}  K(t,s)  = 1_{L^2}  
\]
\end{prop}

This result is strong enough in order to solve the original equation
(\ref{ivp}) iteratively for $f \in L^1L^2$ and $u_0 \in L^2$.  The
kernel of the parametrix $K$ is easy to describe explicitely. To do
this we begin with the evolution operator $\tilde S(t,s)$ in the phase
space.  It corresponds to the transport type operator
\[
D_t + \tilde A = -i (\d_t + a_\xi  \d_x - a_x \d_\xi)  + a(x,\xi) -
\xi  a_\xi  
\]
Solutions are  transported along the Hamilton flow for $D_t + a^w$,
namely
\[
 \dot x = a_\xi \qquad \dot \xi = - a_x 
\]
We denote its solution by $ x^t(x,\xi)$ and $\xi^t(x,\xi)$ where $x$
and $\xi$ are the initial data at time $t=0$.  There is also a
phase shift. We define the  real phase function $\psi$ by
\begin{equation} \label{psi} 
 \dot \psi = -a + \xi a_\xi, \quad  \psi(0,\bar x,\bar \xi) = 0    
\end{equation} 
where $\dot{\psi}$ denotes the differentiation along the flow. 
Then $\tilde S(t,s)$ is given by 
\[
(\tilde S(t,s) u)(x^t,\xi^t) = u(x^s,\xi^s) e^{i (\psi(t,x,\xi) 
- \psi(s,x,\xi))} 
\]
and  the parametrix $K$ has the kernel
 \begin{eqnarray*}
 K(t, y,s,\tilde y)  = 
2^{-d} \pi^{-\frac{3d}2}  \int_{\mathbb{R}^{2n}}  \!\!\!\!\!\!\!\!\!&& 
 e^{- \frac12 ( y - x^t)^2- \frac12 (\tilde y-x^s)^2 + i \xi^t(y-x^t) 
- i \xi^s (\tilde y-x^s)}  
\\  && e^{i (\psi(t,x,\xi)- \psi(s,x,\xi))}  dx \, d\xi. 
\end{eqnarray*}
One can use this directly to prove the dispersive estimates, as in
\cite{MR1833146}.

Here we prefer to use a different approach and work with exact instead
of approximate solutions.  Inspired by the above parametrix,
we seek to obtain a similar representation for the solution.

\begin{prop}
The kernel $K$ of the fundamental solution operator $D_t + a^w$ 
can be represented in the form
\begin{eqnarray}
K(t, y,s,\tilde y)  = 
  \int_{\mathbb{R}^{2n}} \!\!\!\!\!\!\!\!\!&& 
e^{- \frac12 (\tilde y-x^s)^2} e^{ 
- i \xi^s (\tilde y-x^s)}    e^{i (\psi(t,x,\xi)- \psi(s,x,\xi))} \nonumber 
\\ &&   e^{i \xi^t( y-x^t)} G(t,s,x,\xi,y)   dx \, d\xi   
\label{fs} \end{eqnarray}
where the function $G$ satisfies  
\begin{equation}
|(x^t-y)^\gamma \d_x^\aa \d_\xi^\bb \d_{y}^\nu G(t,s,x,\xi,y)| \lesssim
c_{\gamma,\aa,\bb,\nu}
\label{Gbd}\end{equation}
\label{exact}\end{prop}
\begin{proof}
  Without any restriction in generality take $s=0$ and drop the
  argument $s$ in the notations.  Given $u_0 \in L^2$ we need to find
  the solution $u=S(t,0)u_0$ for the equation
\[
(D_t + a^w)u = 0, \qquad u(0) = u_0
\]
  We use the FBI transform to
  decompose $u_0$ into coherent states, and write
\[
u = S(t,0) T^* T u_0 = \int S(t,s) \phi_{x,\xi} Tu(x,\xi) dx d\xi
\]
where the coherent states $\phi_{x,\xi}$ are given by
\[
\phi_{x,\xi}(y) = 2^{-\frac{d}2} 
\pi^{- \frac{3d}4}  e^{ i \xi(x-y)-\frac12(x-y)^2}. 
\]
Then we can define the function $G$  by
\[
G(t,x,\xi,y)   = 2^{-\frac{d}2} \pi^{-\frac{3d}4}  e^{ -i \xi^t(y-x^t)}  
 e^{- i \psi(t,x,\xi)} (S(t,0) \phi_{x,\xi})(y)
\]
so that (\ref{fs}) holds. It remains to prove that
$G$ satisfies the bounds (\ref{Gbd}). For this we need to 
study the regularity of the Hamilton flow and of the phase
function.

\begin{lemma}
Assume that the symbol $A$ satisfies \eqref{a11}. Then 
\begin{equation}
| \d^\aa_x   \d^\bb_\xi x^t | +  |  \d^\aa_x   \d^\bb_\xi   \xi^t| \leq c_{\aa,\bb} \qquad
|\aa|+|\bb| \geq 1
\label{xtxit}\end{equation}
and, for $|\aa|+|\bb| \geq 1$ and $x_0,\xi_0 \in \R^d$,
\begin{equation}
\left|\d_x^\aa \d_\xi^\bb [\xi^t( y-x^t)+
  \psi(t,x,\xi)+\xi_0(x-x_0)]\right|_{x=x_0, \xi=\xi_0} 
\leq c_{\aa,\bb}(1+| y-x^t_0|) 
\label{phase}\end{equation}
\end{lemma}
\begin{proof}
The bound (\ref{xtxit})  is a consequence of the structure of the Hamiltonian 
equations. 

For (\ref{phase}) we note that if any  derivatives fall on $\xi^t$ then
the corresponding term can be easily estimated using (\ref{xtxit}).
It remains to consider the quantity
\[
e = - \xi^t \d_x^\aa \d_\xi^\bb x^t +
  \d_x^\aa \d_\xi^\bb  \psi(t,x,\xi)+ \xi_0 \d_x^\aa \d_\xi^\bb(x-x_0)
\]
At $t=0$ and $(x,\xi) = (x_0,\xi_0)$ we trivially have $e=0$ therefore
it suffices to bound its derivative along the flow.  But
\[
\dot e = a_x(x_t,\xi_t)  \d_x^\aa \d_\xi^\bb x^t - \xi^t \d_x^\aa \d_\xi^\bb 
a_\xi(x_t,\xi_t) - \d_x^\aa \d_\xi^\bb a(x_t,\xi_t) + \d_x^\aa
\d_\xi^\bb (\xi^t a_\xi(x_t,\xi_t))
\]
We use Leibnitz's rule for the last expression. All terms are bounded
except the ones where all derivatives fall on $\xi^t$,
respectively $a_\xi(x_t,\xi_t))$. The latter cancels the second
expression in the above formula, and we are left with 
\[
\dot e = a_x(x_t,\xi_t)  \d_x^\aa \d_\xi^\bb x^t  - \d_x^\aa
\d_\xi^\bb a(x_t,\xi_t) +a_\xi(x_t,\xi_t)  \d_x^\aa \d_\xi^\bb \xi^t +
O(1)
\]
For the middle expression we use the chain rule. The terms
where at least two derivatives fall on the symbol $a$ are bounded.
But this leaves us with only two terms which cancel the first and the
last in the above formula.
\end{proof}

We now continue the proof of Proposition~\ref{exact}.
Using (\ref{phase}) one easily obtains at $x=x_0, \xi=\xi_0$
\begin{eqnarray}
&\!\!\!\!\! \left|\d_{y}^\nu \d_x^\aa \d_\xi^\bb [\xi^t(y - x^t)+ \psi(t,x,\xi) +
\xi_0(x - x_0) -
\xi^t_0(y - x^t_0) - \psi(t,x_0,\xi_0)]
\right| \nonumber\!\!\!\!\! &\\ &\leq c_{\aa,\bb,\nu}
(1+ | y - x^t_0|)  &
\label{phase2}\end{eqnarray}
Then it suffices to prove the estimates (\ref{Gbd}) at $(x_0,\xi_0)$
for the modified function
\[
G_1(t,x,\xi,y)   =  e^{ -i \xi_0^t( y-x_0^t)}  e^{-i\psi(t,x_0,\xi_0)}  
\left(S(t,0) (e^{i\xi_0(x-x_0)} \phi_{x,\xi})\right)(y)
\]
We translate $G_1$ to the origin by setting
\[
G_2(t,x,\xi,y) = G_1(t,x_0+x,\xi_0+\xi,x_0^t+y)
\]
The $x$ and $\xi$ variables are translated so that they are now
centered at the origin. We do not change the notations since at this
point $x$ and $\xi$ appear only in the initial data for $G_1$ and
$G_2$.

The relation (\ref{Gbd}) at $(x_0,\xi_0)$ is replaced by a similar
relation at $(0,0)$,
\begin{equation}
\left| y^\gamma \d_x^\aa \d_\xi^\bb \d_{y}^\nu 
G_2(t,x,\xi,y)\right|_{x=0,\xi=0}
 \lesssim c_{\gamma,\aa,\bb,\nu}
\label{Gbd1}\end{equation}
A routine computation shows that the function $G_2$ solves the modified equation
\[
(D_t + a_2^w(t,y,D_y)) G_2 = 0, \qquad G_2(0) = \phi_{x,\xi}
\]
where
\[
a_2(t,y,\eta) = a(t,x_0^t+y,\xi_0^t+\eta) - a(t,x_0^t,\xi_0^t) - y  a_x(t,x_0^t,\xi_0^t) -
\eta  a_\xi(t,x_0^t,\xi_0^t)
\]
still satisfies (\ref{a11}) but in addition vanishes of second order 
at $0 \in \R^{2d}$. To differentiate it with respect to $x,\xi$ it suffices to
differentiate the initial data.
But the functions 
\[
\d^\aa_x \d^\bb_\xi \phi_{x,\xi}(y)|_{x=0,\xi=0}
\]
are Schwartz functions in $y$. Hence it suffices to consider the
problem
\[
(D_t + a_2^w(t,y,D)) v = 0, \qquad v(0)= v_0
\]
where the initial data $v_0$ is a Schwartz function, and prove that
the solution $v(t)$ is also a Schwartz function. This follows if we
can prove energy estimates for the functions $y^\aa \d^\bb v$, which
we do by induction over $k=|\aa|+|\bb|$. If $k = 0$ then we trivially
have
\[
\|v(t)\|_{L^2} = \|v(0)\|_{L^2}
\]
For $k=1$ we compute the equations for $y v$ and $\d v$:
\[
(D_t + a_2^w(t,y,D)) (y v) = - i (\d_\eta a)^w(t,y,D) v
\]
\[
(D_t + a_2^w(t,y,D)) (\d_y  v) =  i (\d_y a)^w(t,y,D) v
\]
To bound the right hand side we need the next lemma
for the symbol $b = \d_y a_2$ and $b = \d_\xi a_2$.
This is a special case of Theorem 3  in \cite{MR1944027}.

\begin{lemma}
Suppose that the symbol $b(x,\xi)$ satisfies
\[
\|\d^\aa_y \d^\bb_\eta b(y,\eta)\| \leq c_{\aa,\bb} \qquad |\aa|+|\bb|
\geq 1 
\]
and also $b(0,0)= 0$.
Then
\[
\|b^w(y,D) u\|_{L^2} \lesssim \|y u\|_{L^2}+\|\d u\|_{L^2} +\|u\|_{L^2}
\]
\label{b9}\end{lemma}

Using Lemma \ref{b9} and the Gronwall inequality we conclude that
\[
\|y v(t)\|_{L^2} + \|\d v(t)\|_{L^2} \lesssim
\|yv_0\|_{L^2} + \|\d v_0\|_{L^2} + \|v_0\|_{L^2}
\]
It remains to do the induction step.  We denote by $L_k$ all operators
of the form $x^\aa \d^\bb $ with $|\aa|+|\bb|= k$. Suppose that
\[
\sum_{j \leq k} \|L_j v(t)\|_{L^2}
\leq c_k \sum_{j \leq k} \|L_j v_0\|_{L^2}
\]
The functions $L_{k+1} v$ solve a weakly coupled system of the
form
\[
(D_t \!-\! a_2^w) L_{k+1}  v = ( \d_{y,\eta} a_2)^w
L_k  v + \sum_{i \geq 2}^{i+j \leq k+1}  ( \d_{y,\eta}^i a_2)^w
L_j v
\]
For this we use energy estimates and Gronwall's inequality.
The first right hand side term is estimated using Lemma~\ref{b9}
and the second using the induction hypothesis.
\end{proof}

We conclude the section with a supplementary result where we establish
some time regularity of the function $G$ in Proposition~\ref{exact}.
This will not be needed until Section~\ref{degenerate}.

\begin{prop}
Let $a$ be a symbol satisfying (\ref{a11}), and $G$ be as in
Proposition~\ref{exact}. Suppose that
$a$  is smooth in $t$ near some $t_0 \in [0,1]$ and satisfies the additional relations
\begin{equation}
|a_{x}(t_0,x,\xi)| + |a_\xi(t_0,x,\xi)| \leq \mu
\end{equation}
\begin{equation}
 |\d_t^\sigma \d_x^\alpha \d_\xi^\beta a (t_0,x,\xi)|  \leq c_{\aa,\bb,\sigma}
 \mu^{\sigma+1}, \qquad
\sigma \geq 1, \quad |\aa|+|\bb| \geq 1
\label{a11b}\end{equation} 
where $\mu > 1$ is a large parameter. Then $G$ satisfies
the additional bound
\begin{equation}
|(x^t-y)^\gamma \d_t^{\sigma}   \d_x^\aa \d_\xi^\bb \d_{y}^\nu G(t_0,s,x,\xi,y)| \lesssim
c_{\gamma,\aa,\bb,\nu,\sigma} \mu^{\sigma}
\label{Gbdt}\end{equation}
\label{tdif}\end{prop}
\begin{proof}
We begin with an analysis of the time derivatives of the Hamilton
flow.  The counterpart of (\ref{xtxit}) is
\begin{equation}
| \d_t^\sigma \d^\aa_x   \d^\bb_\xi x^t | +  | \d_t^\sigma \d^\aa_x
\d^\bb_\xi   \xi^t| \leq c_{\aa,\bb, \sigma} \mu^\sigma,  \quad
\sigma+ |\aa|+|\bb| \geq 1, \quad t=t_0
\label{xtxit1}\end{equation}
This can be proved by induction with respect to $\sigma$. The details are left for
the reader. An immediate consequence of it is
\begin{equation}
| \d_t^\sigma \d^\aa_x   \d^\bb_\xi a_{x}(t,x^t,\xi^t)|+
 | \d_t^\sigma \d^\aa_x   \d^\bb_\xi a_{\xi}(t,x^t,\xi^t)| \leq
c_{\aa,\bb, \sigma} 
\mu^{\sigma+1}, \quad t=t_0  
\label{axaxi}\end{equation}

Next we examine the steps in the proof of Proposition~\ref{exact}
and track the time derivatives of $G$.

{ \em (i) The phase correction.} Here we need to strenghten 
(\ref{phase2}) to a form which also includes time derivatives, 
\begin{eqnarray}
&\!\!\!\!\!\left|\d_t^\sigma \d_{y}^\nu \d_x^\aa \d_\xi^\bb [\xi^t(y -
  x^t) \!+\! \psi(t,x,\xi) \!+\!
\xi_0(x - x_0)\! -\!
\xi^t_0(y - x^t_0) \!-\! \psi(t,x_0,\xi_0)]
\right| \nonumber \!\!\!\!\!&\\ &\leq c_{\aa,\bb,\nu,\sigma}
(1+ | y - x^t_0|) \mu^{\sigma}, \qquad t=t_0. & \!\!\!\!\! 
\label{phase3}\end{eqnarray}
We know this for $\sigma = 0$, and it is trivial if $|\aa|+|\bb|=0$.
Then we denote by $e$ the expression which is differentiated in the
above formula, and we first compute
\begin{eqnarray*}
\d_{x,\xi} \d_t e &=& \d_{x,\xi} \d_t (\xi^t(y - x^t)+ \psi(t,x,\xi)) \\
&=& \d_{x,\xi} (-a_x(t,x^t,\xi^t)(y - x^t)  - a(t,x^t,\xi^t)) \\
&=& - (y - x^t) \d_{x,\xi} a_x(t,x^t,\xi^t) -
a_\xi(t,x^t,\xi^t)\d_{x,\xi} \xi^t
\end{eqnarray*}
The proof of (\ref{phase3}) is completed using Leibnitz's rule and
(\ref{xtxit1}), (\ref{axaxi}).

{\em (ii) The coordinate change.} The function $G_2$ is obtained from
$G_1$ after a time dependent translation of $x^t_0$. This translation
preserves the bounds (\ref{Gbdt}) on $G$ since by (\ref{xtxit1}) we
know that 
\[
|\d_t^\sigma x_t^0| \leq c_\sigma \mu^\sigma \qquad t=t_0
\]

{\em (iii) The localized evolution.}
The bounds for $\d_t v$ are obtained directly from the equation.
To obtain bounds for higher order time derivatives of $v$ we 
repeatedly differentiate the equation. For this we need to control
the derivatives of the symbol $a_2$. We know that
\[
a_2(t,0,0) = 0, \qquad \d_{x,\xi} a_2(t,0,0) =0
\]
Then the same will apply to the time derivatives of $a_2$.
It remains to show that
\[
|\d_t^\sigma  \d_x^\aa \d_\xi^\bb a(t,x+x_0^t,\xi+\xi_0^t)|\leq
c_{\aa,\bb,\sigma} \mu^{1+\sigma}, \quad |\aa|+|\bb| \geq 2, \quad
\sigma \geq 1, \quad t=t_0
\]
which is easily done using the chain rule and (\ref{xtxit1}), (\ref{axaxi}).
\end{proof}

\subsection{Fixed time estimates}

Here we combine the above representation of the fundamental solution
for $D_t + a^w$ with the curvature condition in order to obtain
pointwise bounds for its kernel.

\begin{prop}
  Let $a \in \Sl$ and $0 \leq k \leq d$ so that for each $(t,x,\xi)
  \in B_\l$ there exists an $d-k$ nondegenerate minor $M$ of
  $\d_\xi^2 a(t,x,\xi)$ satisfying
\[
|\det M| \gtrsim \l^{-(d-k)}.
\]
Then there exists $T > 0$ so that for all $|t-s| < T$  we have
\begin{equation}
\|S(t,s) \chi^w u_0\|_{ L^\infty} \lesssim \l^{\frac{d+k}{2}}
|t-s|^{\frac{-d+k}2}  \lesssim \|u_0\|_{L^1}
\label{pure}\end{equation}
\label{fixedtime}\end{prop}
\begin{proof}
  Without any restriction in generality we take $s=0$.  We fix $t_0$
  in $[0,1]$ and seek to prove the above estimate when $t = t_0$.  
 The result is trivial  if $t_0 < \l^{-1}$. Hence in the sequel we
 assume that $t_0 \geq \l^{-1}$.

We rescale the problem to reduce it to an estimate for  $t=1$. 
If $u = S(t,0) \chi^w u_0$ then we set  
\[
v(t,x) = u\left(\frac t t_0,\frac{x\sqrt{t_0}}{\sqrt\l}\right)
\]
The function $v$ solves the equation 
\[
(D_t + \tilde a^w(t,x,D) )v =0,   \qquad v(0) = \tilde \chi(t,x,D) v_0
\]
where 
\[
\tilde a (t, x,\xi) = t_0 a\left(\frac{t}{t_0}, \frac{x\sqrt{t_0}}{\sqrt\l},
  \frac{\xi\sqrt{\l}}{\sqrt{t_0}}\right), \qquad \tilde \chi(t,x,\xi)
= \chi^w\left(\frac{x\sqrt{t_0}}{\sqrt\l},\frac{\xi
    \sqrt{\l}}{\sqrt{t_0}} \right)
\]
The new frequency scale is $\mu = \sqrt{t_0 \lambda}$.
The rescaled version of \eqref{pure} has the form
\begin{equation}
\|v(1)\|_{L^\infty} \lesssim   \mu^k  \|v_0\|_{L^1}
\label{purer}\end{equation}
It is easy to verify that $\tilde a$ satisfies (\ref{a11}), therefore we can
use the parametrix in Proposition~\ref{exact},
\[
v(t,y)  \!  = \!\!
  \int_{\R^{3d}} \!\! 
 G(t,x,\xi,y) e^{- \frac12 (\ty-x)^2 + i \xi^t(y-x^t) 
- i \xi (y-x)}  
 e^{i \psi(t,x,\xi)} (\tilde \chi v_0)(\ty) dx \, d\xi   d\tilde y
\]
where the symbol $\tilde \chi$ is compactly supported in
\[
\tilde B = \{ |x| \leq \mu t_0^{-1}, \ |\xi| \leq \mu\}
\]
and it is smooth on the scale of $\tilde B$.  The contribution of the
complement of $\tilde B$ to the above integral is negligible. More
precisely we can write
\begin{eqnarray*}
v(t,y)  \!\!\!\! &=& \!\!\!\!\!
  \int_{\tilde B} \!\!
 G(t,x,\xi,y) e^{- \frac12 (\ty-x)^2 + i \xi^t(y-x^t) 
- i \xi (y-x)}  
 e^{i \psi(t,x,\xi)} (\tilde \chi^w v_0)(\ty) dx \, d\xi   d\tilde y\\ & +& O(\mu^{-\infty})
\end{eqnarray*}
Using an $L^1$ bound for $\tilde \chi v_0$ and a trivial estimate for
the kernel, the inequality \eqref{purer} would follow from
\[
\int_{\tilde B} |G(1,x,\xi,y) | d \xi \lesssim \mu^{k}
\]
Given the bounds \eqref{Gbd} for $G$, this reduces to
\begin{equation}
\int_{\tilde B} (1+ |x^1-y |)^{-N} d \xi \lesssim \mu^{k}, \qquad N \ \text{large}
\label{mg}\end{equation}
The key factor here is the dependence of $x^1$ on $\xi$.  We study this
using the linearization of the Hamilton flow. The functions 
\[
X = \frac{\d  {x^t}}{\d {\xi}}, \qquad \Xi = \frac{\d {\xi^t}}{\d {\xi}}
\]
solve the ordinary differential equation along the Hamilton flow
\[
\left\{ \begin{array}{l} \dot{X} = \ \ {\tilde a}_{\xi x} X + {\tilde a}_{\xi \xi} \Xi \cr
\dot{\Xi}=- {\tilde a}_{xx} X - {\tilde a}_{x \xi} \Xi \end{array}  \qquad \right\{ 
\begin{array}{l} X(0) = 0 \cr \Xi(0) = I \end{array}
\]
Since ${\tilde a}_{\xi x},{\tilde a}_{xx} ,{\tilde a}_{x \xi}= O(\sqrt{t_0})$ we obtain
\[
\dot X = {\tilde a}_{\xi \xi} + O(\sqrt {t_0} )
\]
We can also compute
\[
\dot{\tilde a}_{\xi \xi} = {\tilde a}_{t \xi \xi} + {\tilde a}_{\xi
  \xi x} {\tilde a}_\xi - {\tilde a}_{\xi \xi \xi } {\tilde a}_x =
O(\sqrt {t_0} )
\]
Hence we obtain 
\[
X(t) = t( {\tilde a}_{\xi \xi}(0,x,\xi) + O(\sqrt {t_0} ))
\]
which at time $1$ gives
\begin{equation}
\frac{\d {x^1}}{\d {\xi}} = {\tilde a}_{\xi \xi}(0,x,\xi) + O(\sqrt {t_0} )
\label{cxxxi}\end{equation}
 Given $\xi_0 \in \tilde B$ we choose coordinates 
\[
\xi =(\xi',\xi''), \qquad \xi'=(\xi_1, \cdots, \xi_{d-k})
\] 
so that the matrix $(\d^2_{\xi'} \tilde a(0,x,\xi_0) )$ is nondegenerate.
Since $|\d_\xi^3 \tilde a(x,\xi)| \lesssim \mu^{-1}$, it follows that the same
must hold for $\xi \in B(\xi_0,\dd \mu)$ for small fixed $\dd$.  To
prove (\ref{mg}) we split $\tilde B$ into balls of radius $\dd \mu$.
Since
 \[
\int_{B(\xi_0,\dd \mu)} (1+ |x^1-y |)^{-N} d \xi \lesssim \mu^{k}  \sup_{\xi''}
\int_{B'(x_0,\dd \mu)} 
 (1+ |x^1-y |)^{-N} d \xi'
\]
it suffices to show that
\[
\int_{B(\xi_0,\dd \mu)}   (1+ |x^1-y |)^{-N} d \xi' \lesssim 1
\]
But in this region (\ref{cxxxi}) shows that $\d_{\xi'} x^1$ is a small
perturbation of the nondegenerate matrix $\d_{\xi'}^2 \tilde a(x,\xi_0)$.
Hence the above estimate follows.
\end{proof}

\subsection{ Mixed norm estimates }

Here we use the fixed time bounds obtained before in order to derive
space-time estimates in mixed norm spaces.

\begin{prop}\label{oireal2}
Assume that $D_t+a^w$ satisfies {\bf (A2)}, {\bf (A3)} in $B_\l$. 
Let $\chi,\tilde\chi \in \Sl$ be  symbols which are supported in $B_\l$.
Let $u$ solve
\[
(D_t + a^w) u = \tilde \chi^w(x,D) f_1 + f_2  \qquad u(0) = u_0
\]
in $[0,1]$. Then for $(r,s)$ as in \eqref{rs} we have
\begin{eqnarray}  
\|u\|_{L^\infty L^2}+\Vert \chi^w(x,D)u \Vert_{  
\lambda^{\rho(r,s)} L^r L^s}   
&\lesssim &
\Vert f_1  \Vert_{ \lambda^{-\rho(r,s)}  
L^{r'} L^{s'}  } \nonumber \\ &&+ \| f_2\|_{L^1 L^2}
+ \Vert u_0  \Vert_{L^2} \hspace{.5in}
\label{str}\end{eqnarray} 
\end{prop}
\begin{proof}
It suffices to prove this in a sufficiently small time interval,
as we can iterate it and obtain it in the full interval $[0,1]$.

Besides the trivial energy estimates need to show that
\[
\chi^w S(t,s): L^2 \to\lambda^{\rho(r,s)} L^r L^s 
\]
\[
S(t,s) \chi^w:  \lambda^{-\rho(r,s)}  
L^{r'} L^{s'}  \to L^2
\]
\begin{equation}
1_{t > s} \chi^w S(t,s)\tilde \chi^w :   \lambda^{-\rho(r,s)}  
L^{r'} L^{s'}  \to  \lambda^{\rho(r,s)} L^r L^s
\label{tsi}\end{equation}
The first two statements are dual. Using a $TT^*$ argument
they reduce to the third without $1_{t > s}$. But the third with $1_{t >
  s}$ is dual to the third with $1_{t < s}$ instead, therefore it
implies the first two. It remains to prove the third.

 The energy estimates yield the trivial bound 
\[
\|S(t,s)\|_{L^2 \to L^2} \leq 1
\]
On the other hand the decay estimates in Proposition~\ref{fixedtime}
show that 
\[
\|\chi^w S(t,s) \tilde \chi^w \|_{L^1 \to L^\infty} \lesssim |t-s|^{-\frac{d-k}2} 
\l^{\frac{d+k}2}
\]
If $r > 2$ then (\ref{tsi}) follows from the two estimates above by
interpolation and the Hardy-Littlewood-Sobolev inequality.
The case $r=2$ needs some extra work, and it can be obtained as in  
 Keel-Tao~\cite{KT}.
\end{proof}

\begin{cor} Assume that $D_t+\tilde a$ satisfies {\bf (A2)}, {\bf (A3)} in $B_\l$. 
Let $\chi \in \Sl$ be  a symbol which is supported in
$B_\l$. Then\footnote{ If $k = n-2$ then the $L^2L^\infty$ bound is not
  quite true, but all the intermediate bounds are still valid.}
\begin{equation}
\|\chi^w(x,D) u\|_{L^\infty(L^2)\cap  
\lambda^{\frac{n+k-2}{2(n-k)}} L^2(L^{\frac{2(n-k)}{n-k-2}})}   
\lesssim \|u\|_{L^2} + \|(D_t+A) u\|_{L^2}.
\label{l2l2e}\end{equation}
\label{l2l2}\end{cor}
This follows easily from the previous proposition applied to $\tilde
\chi^w u$, where the symbol $\tilde \chi$ is chosen to equal $1$ in a
neighbourhood of the support of $\chi$.

\section{The parametrix in the general case}
\label{five}
In this section we construct two parametrices for operators in
cannonical form. These constructions do not use much of the 
structure of our operators, so we prefer to write it in a more
abstract setup.

Thus, given selfadjoint operators $A(t)$, $B(t)$ in a Hilbert space
$X$ for $t \in [0,1]$, we seek to construct a parametrix for the
operator
\[
D_t + A + iB
\]
For simplicity we assume that both $A(t)$ and $B(t)$ are bounded
and smooth as functions of $t$, but the constants in our estimates
are independent of any such bounds. 

Making a slight abuse of notation we use simply $\|\cdot\|$ both for
the norm of $X$ and for the operator norm in $L(X)$ through this
section. All other norms will be indicated with a subscript. We also
abbreviate the notation $\|\cdot\|_{L^q}:= \|\cdot\|_{L^q(0,1;X)}$.

The main relation connecting $A$ and $B$ is the fixed time commutator
estimate
\begin{equation}
\| [D_t + A,B] u\| \lesssim \|Bu\|+\|u\|
\label{comut}\end{equation}

We do not use this directly, instead we first obtain a simple
consequence of it. Denote by $S(t,s)$ the unitary evolution generated
by $D_t + A$ in $X$.

\begin{lemma} Assume that \eqref{comut} holds in $[0,1]$. Then
\begin{equation}
\| (B(t) S(t,s)- S(t,s) B(s)) u\| \lesssim |t-s| (\| B(s)
u\| + \|u\|)
\label{bdiff}\end{equation}
\label{l0}\end{lemma}
\begin{proof} 
First we compute the equation for $ B(t) S(t,s) u$,
\[
(D_t + A)  B(t) S(t,s) u=  [D_t + A, B] S(t,s) u
\]
Using energy estimates and \eqref{comut} we obtain
\[
\| B(t) S(t,s) u\| \lesssim \|B(s) u\| + \int_{s}^t \|B(r)
S(r,s) u\| + \|S(r,s) u\| dr
\]
Applying Gronwall's lemma this yields
\begin{equation} 
\| B(t) S(t,s) u\| \lesssim \|B(s) u\|+\|u\|
\label{gron2}\end{equation}
Then we write
\[
(D_t + A) (B(t) S(t,s)- S(t,s) B(s)) = [D_t + A, B] S(t,s)
\]
which implies that
\begin{equation}\label{gron1} 
\| (B(t) S(t,s)- S(t,s) B(s)) u\| \lesssim \int_{s}^t \|B(r)
S(r,s) u\| + \|S(r,s)u\| dr
\end{equation} 
To obtain \eqref{bdiff} it suffices to estimate the right hand side
using (\ref{gron2}).
 \end{proof}

\subsection{ A simple parametrix }
Assume first that $D_t + A$ commutes with $B$. Then we can produce an
exact parametrix for $D_t + A + iB$, namely
\[
H(t,s) = 1_{(t-s)B(t) < 0} e^{(t-s)B(t)}  S(t,s)   = 
 S(t,s) 1_{(t-s)B(s) < 0} e^{(t-s)B(s)}  
\] 
where the $B$ dependent part is interpreted in the sense of operator
calculus for selfadjoint operators.

In the noncommuting case we can use either of these two expressions
as a parametrix, but they are no longer equal. Set
\begin{equation}
H(t,s) =  S(t,s) 1_{(t-s)B(s) < 0} e^{(t-s)B(s)}
\label{h}\end{equation}

\begin{prop}
Suppose that \eqref{bdiff} holds for $t,s \in [0,1]$. Then
for $t,s$ in a bounded interval  we have the fixed time estimates
\begin{equation}
\|H(t,s)\|_{X \to X} \leq 1, \quad
\|H(t,s)B(s)\|_{X \to X} \leq |t-s|^{-1}, 
\label{ks1}\end{equation}
\begin{equation}
\|B(t) H(t,s)\|_{X \to X} \lesssim |t-s|^{-1}, \quad
\|B(t)H(t,s)B(s)\|_{X \to X} \lesssim |t-s|^{-2}
\label{ks2}\end{equation}
In addition, the following space time error estimate holds:
\begin{equation}
\| (D_t+A+iB) H -I  \|_{L^1 \to L^\infty} \lesssim 1
\label{ks3}\end{equation} 
\label{pari}\end{prop}
\begin{proof}
The bounds in (\ref{ks1}) are trivial. So are the ones in (\ref{ks2})
provided that $D_t + A$ and $B$ commute. Otherwise,
they follow from \eqref{bdiff}.

It remains to prove (\ref{ks3}). A simple computation shows that
\[ 
[(\! D_t+A+iB) H -I ](t,s) \!=\! i (\! B(t) S(t,s) - S(t,s) B(s))
1_{(t-s)B(s) < 0} e^{(t-s)B(s)}
\]
so (\ref{ks3}) follows also from \eqref{bdiff}.
\end{proof}

\subsection{A more robust parametrix.}

While the above parametrix is quite simple, it is not clear whether
one can use it to show that the operator $D_t+A+iB$ inherits most of
the dispersive estimates from $D_t+A$.  To do this we use a modified
version of the above parametrix, which is somewhat reminiscent of the
Littlewood-Paley theory.  We consider a dyadic partition of the
unity
\[
1 = \sum_{j=0}^\infty \kappa_j^2
\]
where the functions $\kappa_j$ are supported in $\{2^{j} \leq \max\{|\xi|,1\}
\leq 2^{j+2} \}$ and are smooth on the scale of their support.  For $j
> 0$ we denote by $\kappa_j^+$ respectively $\kappa_j^-$ the parts of
$\kappa_j$ supported in $[0,\infty)$, respectively $(-\infty, 0]$. For
$j=0$ we set $\kappa_0^+ = \kappa_0$ and $ \kappa_0^-= 0$.  Using the functional
calculus for selfadjoint operators we define the dyadic operators
$\kappa_j(B(t))$.  Then the modified parametrix $H$ has the form
\begin{eqnarray}
H(t,s) &=& 1_{t > s} \sum_j \kappa_j^-(B(t))S(t,s) \kappa_j^-(B(s))
e^{(t-s)B(t)} \label{param1}
 \\&&-  1_{t < s} \sum_j \kappa_j^+(B(t))S(t,s) \kappa_j^+(B(s))
e^{(t-s)B(t)} \nonumber
\end{eqnarray}

To measure the regularity of this parametrix we introduce some
function spaces which depend only on $A$ and not on $B$.  This will
allow us later to transfer the dispersive estimates from $D_t - A$
to $D_t - A + iB$.  We begin with the energy space $L^\infty$, 
and the ``classical'' solutions for $D_t + A$ which are in the
space $W^{1,1}_A$ of $X$ valued functions with norm
\[
\|u\|_{W^{1,1}_A} = \|u\|_{L^\infty} + \|(D_t+A) u\|_{L^1}.
\]

In between these spaces we define the space $V^2_A$ of functions with
bounded $2$-variation along the $D_t +A$ flow, with norm
\[
\|u\|_{V^2_A}^2 = \|u(0)\|^2+ \sup_{ (t_j) \in \TT} \sum_j \|u(t_{j+1}) -
S(t_{j+1},t_j)u(t_j)\|^2
\] 
where $\TT$ is the set of finite increasing sequences in $[0,1]$.
Functions in $V^2_A$ have at most countably many discontinuities.  To
elliminate functions in $V^2_A$ which are zero a.e. we assume that all
functions in $V^2_A$ are right continuous.  The $V_A^2$ space satisfies
\begin{equation}
W^{1,1}_A \subset V^2_A \subset L^\infty
\label{wvl}\end{equation}
The closure of the space of smooth $X$ valued functions in $V^2_A$ is
$V^2_A \cap C$, i.e. the subspace of continuous functions in $V^2_A$.

Making a slight abuse of notation we denote 
the dual space  of $V^2_A \cap C$ by $(V^2_A)^{*}$.
This is a space of distributions which has an atomic structure.
There are two kinds of atoms in $(V^2_A)^{*}$:
\begin{eqnarray*}
\text{ ($L^1$ type)}\qquad & f = f_0(x) \dd_{t_0}, & \quad
  \|f_0\| =1
\\[2mm]
\text{($2$-variation type)}  & \displaystyle f  \!=\! \sum_j  
 S(t_{j+1},t_j) f_j \dd_{t_{j+1}}  -  f_j \dd_{t_j}, &
\sum_j \|f_j\|^2 = 1
\end{eqnarray*}
where $(t_j) \in \TT$. This has to be understood in the sense that
the atomic space generated by these atoms is a weakly* dense
subspace of  $(V^2_A)^{*}$ with an equivalent norm.
It follows from \eqref{wvl}  that
\[
L^1 \subset  (V^2_A)^{*} \subset (D_t +A) L^\infty 
\]

Following is the main result of this section, which describes the
mapping properties of the parametrix $H$ in  (\ref{param1}) in terms
of the $V^2_A$ and the  $(V^2_A)^{*}$ spaces.

\begin{prop}
  Assume that $A$, $B$ are selfadjoint operators which satisfy
  (\ref{bdiff}). Then the parametrix $H$ for $D_t+A+iB$ in \eqref{param1}
  satisfies the estimates
\begin{equation}
H: (V^2_A)^{*} \to V^2_A
\label{k}\end{equation}
\begin{equation}
(D_t+A+iB)H-\id : (V^2_A)^{*} \to L^\infty 
\label{pk}\end{equation}
\label{precise}
\end{prop}

\begin{remark}
  The same result and proof apply if we construct a para\-me\-trix for the
  operator $D_t+A+\alpha B + iB$, $\alpha \in \R$ using the same
  formula but with $e^{(t-s)B(s)}$ replaced with
  $e^{(t-s)(1-i\alpha)B(s)}$.
\label{calpha} \end{remark}

\begin{proof}
  We first observe that we can conjugate the result with respect to
  the group of isometries generated by $A$ and reduce the problem to
  the case when $A=0$. We omit $A$ in the notation $V^2 := V^2_0$ and $
  (V^2)^* := (V^2_0)^*$.  The condition \eqref{bdiff} implies that
\begin{equation}
\| (B(t) - B(s)) u\| \lesssim |t-s| (\| B(s)
u\| + \|u\|)
\label{bt1}\end{equation} 

Next we consider the operators $\kappa_j(B(t))$. They depend on $t$,
but the next lemma shows that this dependence is mild.

\begin{lemma}
  Suppose that $B$ satisfies (\ref{bt1}). Then for the operators $\kappa_i(B(t))$
  defined above we have

a) (bound for low modes of $B$)
\begin{equation}
\| [\kappa_j(B(t)) - \kappa_j(B(s)) ] f\|  \lesssim |t-s|
2^{-j} \|B(s) f\|  \label{chidiff}
\end{equation}

b) (almost orthogonality)
\begin{equation}
\| \kappa_i(B(t))  \kappa_j(B(s)) \| \lesssim |t-s| 2^{-|i-j|}
\qquad  |i-j| \geq 3
\label{ij}\end{equation}

c) (Lipschitz bound)
\begin{equation}
\| [\kappa_j(B(t)) - \kappa_j(B(s)) ] f\|  \lesssim |t-s|
\| f\|  \label{chidiffh}
\end{equation}

d)  (bound for high modes of $B$)
\begin{equation}
\|B(t) [\kappa_j(B(t)) - \kappa_j(B(s)) ] f\|  \lesssim |t-s| 2^j
\| f\|  \label{chidiffhh}
\end{equation}
\label{ort}
The same estimates hold if we replace $\kappa_j$ by $\kappa_j^\pm$
or any other bump functions on the same scale and with similar supports.
\end{lemma}

\begin{proof}
a) We use the representation 
\[
\kappa_j(B(t)) = \int e^{i \tau B(t)} \widehat{\kappa}_j(\tau) \ d \tau
\]
Since $\kappa_j$ is an integrable bump function on
the $2^{-j}$ scale, the estimate (\ref{chidiff}) follows
if we prove that
\[
\|(e^{i \tau B(t)} -e^{i \tau B(s)}) f\|  \lesssim |\tau| |t-s|
 \|B(s) f\|  
\]
For this we compute
\[
\frac{d}{d\tau}\left(e^{i \tau B(t)} e^{- i \tau B(s)}\right) = e^{i
  \tau B(t)} (B(t)-B(s)) e^{-i \tau B(s)}
\]
which by (\ref{bt1}) gives
\[
\|(e^{i \tau B(t)} -e^{i \tau B(s)}) f\|  \lesssim \int_0^\tau
\| (B(t)-B(s)) e^{i \theta B(s)} f\|  d\theta \lesssim 
|\tau| |t-s| \|B(s) f\| 
\]

b) By duality we can assume without any restriction in generality that
$i - j \geq 3$. Then by (\ref{chidiff}) we get 
\begin{eqnarray*}
\|\kappa_i(B(t))  \kappa_j(B(s)) f \|  &=& 
\|(\kappa_i(B(t)) - \kappa_i(B(s))) \kappa_j(B(s)) f \|  
\\ &\leq&
2^{-i} |t-s| \|B(s)\kappa_j(B(s)) f \| 
\\
&\lesssim&
 2^{j-i} |t-s| \|f\| 
\end{eqnarray*}

c) We write
\begin{eqnarray*}
\kappa_j(B(t)) - \kappa_j(B(s)) &= &  (\kappa_j(B(t)) - \kappa_j(B(s))) 
\sum_{i<j+3} \kappa_{i}(B(s)) \\ &  & + \sum_{ i \geq j+3} \kappa_j(B(t)) \kappa_i(B(s))
\end{eqnarray*}
For the first term we use \eqref{chidiff}, while for the second
we use \eqref{ij}.

d) Compute
\begin{eqnarray*}
B(t)(\kappa_j(B(t)) - \kappa_j(B(s))) & = & [B(t) \kappa_j(B(t)) - B(s)\kappa_j(B(s))] 
\\ & & + (B(t)-B(s)) \kappa_j(B(s))
\end{eqnarray*}
For the first term we use \eqref{chidiffh} while for the second we use
\eqref{bt1}.
\end{proof}

We want to replace the estimates (\ref{k}) and (\ref{pk}) by their
dyadic counterparts. For this we need the following

\begin{lemma}
Suppose $B$ satisfies (\ref{bt1}). Then
we have the estimate
\[
\| \sum_{j \in \N} \kappa_j (B) u_j \|_{V^2}^2 \lesssim \sum_{j \in \N} \| u_j\|_{V^2}^2
\]
and its dual
\[
\sum_{j \in \N}  \| \kappa_j(B) f\|_{(V^2)^*}^2 \lesssim \| f\|_{(V^2)^*}^2. 
\]
The same holds if we replace $\kappa_j$ by any other symbols
with similar support, size and regularity.
\label{orth}\end{lemma}

\begin{proof}
Denote 
\[
u =  \sum_{j \in \N} \kappa_j (B) u_j.
\]
Then the first estimate follows from the definition of the $V^2$ norm
and the inequality
\begin{equation}
\|u(t)-u(s)\| ^2 \lesssim  \sum_{j \in \N} \| u_j(t) -  u_j(s)\| ^2
+ |t-s|(\|u_j\|_{L^\infty}^2 +\|v_j\|_{L^\infty}^2 )
\label{ciu-v}\end{equation}
To prove this we write
\begin{eqnarray*}
u(t) - u(s) &=&  \sum_{j \in \N}  \kappa_j(B(t))u_j(t) -  \kappa_j(B(s))u_j(s)
\\ &=&  \sum_{j \in \N}  \kappa_j(B(t))(u_j(t) - u_j(s))  \\ &&  + 
 \sum_{j \in \N}    (\kappa_j(B(t))- \kappa_j(B(s)))u_j(s)
\end{eqnarray*}
The terms in the first sum are clearly almost orthogonal, and are estimated
by the first right hand side term in \eqref{ciu-v}. To estimate the
second sum by the second right hand side term in \eqref{ciu-v} we need
to prove that its terms are almost orthogonal as well. Precisely, it
would suffice to show that we have the off-diagonal decay
\[
\| (\kappa_j(B(t))- \kappa_j(B(s)) (\kappa_k(B(t))- \kappa_k(B(s))\|
\lesssim 2^{-|k-j|} |t-s|
\]
For $|i-j| \leq 3$ this follows from (\ref{chidiffh}). For $|i-j| > 3$
we have
\begin{eqnarray*}
& (\kappa_j(B(t))- \kappa_j(B(s)) (\kappa_k(B(t))- \kappa_k(B(s)) = &\\[2mm]
&-\kappa_j(B(t))\kappa_k(B(s)) -\kappa_j(B(s)\kappa_k(B(t))&
\end{eqnarray*}
and for each of the two terms we use (\ref{ij}).

\end{proof}

We now return to the proof of the theorem. Without any restriction in
generality we consider only the forward part of the parametrix $H$.
Due to Lemma~\ref{orth}, the bound (\ref{k}) for $H$ will follow from the dyadic
estimates
\begin{equation}
\| H_j f\|_{V^2} \lesssim  \|f\|_{(V^2)^*}
\label{dy} \end{equation}
where
\[
H_j(t,s) = 1_{t > s} \kappa_j^-(B(t)) \kappa_j^-(B(s))
e^{(t-s)B(t)} 
\]
Since $(V^2)^*$ is an atomic space it suffices to prove \eqref{dy}
when $f$ is a $(V^2)^*$ atom. We begin with an $L^1$ type atom
$f = f_0 \dd_{t_0}$ for which (\ref{dy}) is a consequence of the next
simple lemma:

\begin{lemma}
For a fixed $t_0 \in [0,1]$ set 
\[
v(t) = H_j(t,t_0) f_0
\]
Then
\[
\|  e^{ 2^j (t-t_0)} v\|_{V^2} \lesssim   \|f_0\| 
\]
\label{longt}\end{lemma}

The additional gain provided by the exponential factor is not needed 
for $L^1$ type atoms, but we will need it later for the $2$-variation
type atoms.

It remains to prove \eqref{dy} for an atom $f \in (V^2)^*$ of the form
\begin{equation}
f = \sum_k  (\dd_{t_{k+1}} -  \dd_{t_k}) f_k 
\label{atom}\end{equation}
where $(t_k) \in \TT$ is an increasing  finite sequence.
Denote
\[
u =  H_j f, \qquad u_k =  H_j   (f_k \dd_{t_{k+1}} -  f_k \dd_{t_k}) 
\]
or, more explicitely,
\begin{equation}
u_k(t)  = H_j(t,t_{k+1}) f_k -  H_j(t,t_{k}) f_k
\label{uk}\end{equation}
For each $k$ the function $u_k$ is supported in $[t_k,1]$
and decays exponentially in time on the $2^{-j}$ time scale.
We decompose $u$ into three parts,
\[
u = v_1 + v_2 + v_3
\]
where 
\begin{eqnarray*}
&\displaystyle v_1 = \sum_k 1_{[t_k,t_{k+1}]} u_k, \qquad v_2 =  \sum_{\{k:t_{k+1}-t_k>
  2^{-j}\}} 1_{t > t_{k+1}} u_k, & \\ &\displaystyle v_3 =  \sum_{\{k:t_{k+1}-t_k\leq
  2^{-j}\}} 1_{t > t_{k+1}} u_k.&
\end{eqnarray*}
The terms in $v_1$ have disjoint supports, and the square 
summability with respect to $k$ is inherited from $f_k$. 
Hence it suffices to consider a single $f_k$, for which the bound 
follows from Lemma~\ref{longt}.

The terms in $v_2$ do not have disjoint supports. However, they decay
in time on the $2^{-j}$ scale while their starting points $t_k$ are at
least $2^{-j}$ separated because the intervals $[s_k,t_k]$ are
disjoint. Hence they are almost orthogonal, and again it suffices to
consider a single $f_k$. But then we can use again Lemma~\ref{longt}.
Note that in this case there is no significant cancellation between the inputs at
times $s_k$ and $t_k$.

The terms in $v_3$ also decay exponentially on the $2^{-j}$ scale.
However, they correspond to intervals $[t_k,t_{k+1}]$ of size less than
$2^{-j}$ which can be closer then $2^{-j}$, so we loose the
orthogonality with respect to $k$.  We partition the unit interval in subintervals of
length $2^{-j}$, and group the intervals $[t_k,t_{k+1}]$ together
in bunches contained in single $2^{-j}$ subintervals.  The outputs of
different bunches are almost orthogonal, so we only need to worry
about a single bunch.

Within a single bunch the orthogonality is lost.  However, the
intervals are disjoint so within each $2^{-j}$ subinterval we retain
control of the sum of the lengths of $[t_k,t_{k+1}]$. Another redeeming
feature is that now there is some cancellation between the input at
times $t_k$ and $t_{k+1}$.

Then it suffices to show that
\begin{lemma}
Let $u_k$ be as in \eqref{uk}. Then
\[
\| 1_{t>t_{k+1}} u_k  \|_{V^2}  
 \lesssim |t_{k+1}-t_k| 2^j   \|f_k\|. 
\]
\label{shortt}\end{lemma}
This lemma is only interesting if $|t_{k+1}-t_k| \leq 2^{-j}$,
otherwise it follows from Lemma~\ref{longt}.
To obtain  \eqref{k} it remains to prove Lemma \ref{longt}
and Lemma \ref{shortt}.

\begin{proof}[Proof of Lemma~\ref{longt}:]
Recall that
\[
v(t) = 1_{t > t_0} \kappa_j^-(B(t))  \kappa_j^- (B(t_0))
e^{(t-t_0)B(s)} f_0
\]
and note the trivial  bound,
\begin{equation}
\|e^{2^j (t- t_0)} v(t)\| \lesssim  \|f_0\|
\label{linf0}\end{equation}
To obtain the conclusion of the lemma we prove a stronger result,
which asserts that $v$ is Lipschitz on the $2^{-j}$ scale and decays
exponentially on the same scale.  More precisely, we claim that
\[
\|v(\tau_1)-v(\tau_2)\| \lesssim |\tau_1 - \tau_2| 2^j e^{2^j (t_0-\tau_1)}
\|f_0\| \qquad t_0 \leq \tau_1 \leq \tau_2 
\]
Indeed,
\begin{eqnarray*}
v(\tau_1)-v(\tau_2) &=& [\kappa^-_j(B(\tau_1)) - \kappa_j^-(B(\tau_2))]  \kappa_j^-
 (B(t_0))e^{(\tau_1-t_0)B(t_0)} f_0
\\
& &+ \kappa^-_j(B(\tau_2))  \kappa_j^- (B(t_0)) [e^{(\tau_1-t_0)B(s)} -
e^{(\tau_2-t_0)B(s)}] f_0
\end{eqnarray*}
For the first term we use (\ref{chidiffh}), while the bound for the
second term is trivial since for $\xi \in [-2^{j+2},-2^j]$ we have
\[
 |\kappa_j^-(\xi) [e^{(\tau_1-t_0)\xi} -e^{(\tau_2-t_0)\xi}]| \lesssim
 2^{j}|\tau_1-\tau_2| e^{2^j(t_0-\tau_1)}
\]
\end{proof}

\begin{proof}[Proof of Lemma~\ref{shortt}:]
As before,  begin with a pointwise estimate for $t \ge t_{k+1}$, 
\begin{equation}
\|u_k(t)\| \lesssim |t_{k+1}-t_k| 2^{j}  e^{2^j (t_{k+1} -t)}\|f_k\|
\label{linf}\end{equation}
To prove it we write
\[
u_k  = w_1+w_2
\]
where
\begin{eqnarray*}
w_1(t) &=& \kappa^-_j(B(t))\kappa_j^-(B(t_{k+1}))
(e^{(t-t_{k+1})B(t_{k+1})} - e^{(t-t_k)B(t_{k+1})}) f_k 
\\[2mm]
w_2(t)&=& \kappa_j^-(B(t)) [\phi_j(B(t_{k+1}) - \phi_j(B(t_k))] f_k
\end{eqnarray*}
where
\[
\phi_j(b) = \kappa_j^-(b) e^{(t-t_k)b}.
\]
The bound for  $w_1$ follows from the inequality
\begin{equation}
| e^{(t-t_{k+1}) \xi} - e^{(t-t_k) \xi}| \lesssim |t_{k+1}-t_k|
2^{j}  e^{2^j (t_{k+1} -t)}, \qquad \xi \in [-2^{j+2},-2^{j}]
\label{ety}\end{equation}
The estimate for $w_2(t)$  follows from  
(\ref{chidiffh}) applied to  $\phi_j$.

Next we seek a similar Lipschitz bound for $ t_{k+1} < \tau_1 < \tau_2 $,
namely 
\begin{equation} 
\|u_k(\tau_1) - u_k(\tau_2)\| \lesssim  |\tau_1-\tau_2|  |t_{k+1}-t_k|
 2^{2j} 
e^{2^j (t_k -\tau_1)}\|f_k\|
\label{lip} \end{equation}
We split $u_k$ as above, $u_k = w_1+w_2$. For $w_1$ this bound is
obtained as in Lemma~\ref{longt}, using \eqref{chidiffh} and symbol
bounds. It remains to prove it for $w_2$. We denote
\[
\phi_j(\xi)=\kappa_j^-(\xi) e^{(\tau_1-t_{k})\xi}, \qquad \psi_j(\xi)
= \kappa_j^-(\xi) e^{(\tau_1-t_{k})\xi} (1- e^{(\tau_1-\tau_2)\xi})
\]
and represent
\[
w_2(\tau_1)-w_2(\tau_2) = w_3 + w_4 
\]
where
\begin{eqnarray*}
w_3 &=&  (\kappa^-_j(B(\tau_1))-\kappa_j^-(B(\tau_2))) [\phi_j(B(t_{k+1}))
-  \phi_j (B(t_k))] f_k
\\[2mm]
w_4 &=& \kappa_j^-(B(\tau_2)) (\psi_j(B(t_{k+1})) - \psi_j(B(t_k))) f_k
\end{eqnarray*}
For $w_3$ we use (\ref{chidiffh}) twice together with the fact that
$\phi_j$ is a bump function on the $2^j$ scale, of size $e^{2^j
  (t_{k+1} -\tau_1)}$.  Finally, the bound for $w_4$ needs
(\ref{chidiffh}) for $\psi_j$, which is a bump function on the $2^j$
scale and of size $|\tau_1-\tau_2| 2^j e^{2^j(t_{k+1}-\tau_1)}$.
\end{proof} 

We now continue with the proof of Proposition~\ref{precise},
\eqref{pk}. Recall that we have reduced the problem to the case
when $A=0$, and denote
\[
L = -i [I- (D_t +i B) H]
\]
It suffices to look at the forward part of $L$, 
\begin{eqnarray*}
L(t,s)& =& 1_{t> s} \sum_j  \kappa_j^-(B(t))(B(t)-B(s)) \kappa_j^-(B(s))
e^{(t-s)B(s)} \\ & &+ \{ \d_t \kappa_j^-(B(t)) \} \kappa_j^-(B(s))
e^{(t-s)B(s)}. 
 \end{eqnarray*}
We need  to prove that
\begin{equation}
\|L f\|_{L^\infty} \lesssim \|f\|_{(V^2)^*}.
\label{khf}\end{equation}
It suffices to do this in the special case when $f$ is an atom.
We denote
\[
L^1_j(t,s) =  1_{\{t>s\}}\kappa_j^-(B(t))  (B(t)-B(s)) \kappa_j^- (B(s))
e^{(t-s)B(s)} 
\]
\[
L^2_j(t,s)=  1_{\{t>s\}} \{\d_t  \kappa_j^-(B(t)) \} \kappa_j^-(B(s))
e^{(t-s)B(s)} 
\]
The difference between these two components is that $L_j^1$ keeps the
size of the frequency, but $L_j^2$ does not, so we need to gain some
decay off the diagonal.  Arguing exactly as in the case of (\ref{k}),
the problem reduces to the two counterparts of the estimates \eqref{linf0} and
\eqref{linf} in Lemma \ref{longt}, respectively Lemma \ref{shortt}.
These are stated in the next two lemmas. The first lemma implies
 \eqref{khf} for $L^1$ atoms,  and also for 2-variation type atoms 
with $t_{k+1} - t_k \geq 2^{-j}$.  
\begin{lemma}
For  $t > s$ we have
\begin{equation}
\|L_j^1(t,s) g\| \lesssim e^{ 2^j (s-t)} \|g\|
\label{l1ts}\end{equation}
\begin{equation}
\|\kappa_i(B(t)) L_j^2(t,s) g\| \lesssim 2^{- |i-j|} 
e^{2^j (s-t)} \|g\|.
\label{l2ts}\end{equation}
\label{longt1}\end{lemma}
The second lemma allows us to prove  \eqref{khf} for 2-variation type atoms 
with $t_{k+1} - t_k \leq 2^{-j}$.  

\begin{lemma}
Suppose $|t_{k+1}-t_k| \leq 2^{-j}$ and $t > t_{k+1}$. Then 
\begin{equation}
\|(L_j^1(t,t_{k+1})-L_j^1(t,t_k)) g \| 
  \lesssim |t_{k+1}-t_k| 2^j   e^{2^j (t_{k+1}-t)}  \|g\|
\label{l1tts} \end{equation}
and
\begin{equation}
\| \kappa_i(B(t)) (L_j^2(t,t_{k+1})-L_j^2(t,t_k)) g \|
 \lesssim |t_{k+1}-t_k| 2^j 2^{- |i-j|} e^{2^j (t_{k+1}-t)}  \|g\|.
\label{l2tts}
\end{equation}
\label{shortt1}\end{lemma}

\begin{proof}[Proof of Lemma~\ref{longt1}:]
The bound (\ref{l1ts}) follows from (\ref{bt1}),
\begin{eqnarray*} 
\|L_j^1(t,s)g\| &\lesssim & |t-s|\left( \| B(s) \kappa_j^- (B(s))
e^{(t-s)B(s)} g\| \right. \\ & &+\ \left.\|\kappa_j^- (B(s)) e^{(t-s)B(s)}g\| \right)
\\ & \lesssim& 2^j |t-s|e^{2^j (s-t)}\|g\|.
\end{eqnarray*} 
For (\ref{l2ts}) we note that 
\[
\|\kappa_j^-(B(s)) e^{(t-s)B(s)}\| \lesssim  e^{2^j  (s-t)}.
\]
Then it remains to show that
\begin{equation}
\|\kappa_i(B(t))\d_t \kappa^-_j(B(t)) \| \lesssim 2^{-|i-j|}.
\label{cdtc}\end{equation}
We consider two cases. If $i > j$ then we use
(\ref{chidiffhh}):
\begin{eqnarray*}
\|\kappa_i(B(t)) \d_t \kappa_j^-(B(t)) g\| \lesssim 
2^{-i} \| B(t) (\d_t \kappa_j^-(B(t)))   g\|
\lesssim
2^{j-i}\|g\|.
\end{eqnarray*}
If $i \leq j$ then we use duality and (\ref{chidiff}):
\begin{eqnarray*}
\|\d_t \kappa^-_j(B(t)) \kappa_i(B(t)) g\| \lesssim
2^{-j} \|B(t) \kappa_i(B(t)) g\|
\lesssim 2^{i-j}\| g\|.
\end{eqnarray*}
\end{proof}
\begin{proof}[Proof of Lemma~\ref{shortt1}:]
For (\ref{l1tts}) we write
\[
(L_j^1(t,t_{k+1})-L_j^1(t,t_k)) g = w_1 + w_2 + w_3
\]
where
\begin{eqnarray*}
w_1 &=&  \kappa^-_j(B(t))  (B(t)-B(t_{k+1})) (\phi_j (B(t_{k+1}))
- \phi_j (B(t_k)) )g
\\[2mm]
w_2 &=&  \kappa^-_j(B(t))  (B(t)-B(t_{k+1})) \phi_j(B(t_k))
(1-e^{(t_{k+1}-t_k)B(t_k)})g
\\[2mm] 
w_3 &=&  \kappa^-_j(B(t)) (B(t_{k+1}) - B(t_k))  \kappa_j^- (B(t_k))
e^{(t-t_k)B(t_k)} g
\end{eqnarray*}
with 
\[
\phi_j(\xi)=\kappa_j^- (\xi) e^{(t-t_{k+1})\xi}.
\]
For $w_1$ we use (\ref{bt1}) to get
\begin{equation}
\| \kappa_j^-(B(t))  (B(t)-B(t_{k+1}))\| \lesssim 2^j |t-t_{k+1}| 
\label{pr}\end{equation}
and (\ref{chidiffh}) to obtain
\[
\|\phi_j (B(t_{k+1}))- \phi_j (B(t_k))\| 
\lesssim 2^j |t_{k+1}-t_k| e^{2^j (t_{k+1} -t)}.
\]
For $w_2$ and $w_3$ we combine (\ref{pr}) with straightforward
symbol bounds.

It remains to prove (\ref{l2tts}). Set
\[
\kappa_i(B(t)) (L_j^2(t,t_{k+1})-L_j^2(t,t_k)) g = w_4+w_5
\]
where
\begin{eqnarray*}
w_4&=& \kappa_i(B(t)) (\d_t \kappa^-_j(B(t)))(\phi_j(B(t_{k+1})) -\phi_j(B(t_k))) g
\\[2mm]
w_5&=&\kappa_i(B(t)) (\d_t \kappa^-_j(B(t)))\phi_j(B(t_k))(1 - e^{(t_{k+1}-t_k)B(t_k)}) g
\end{eqnarray*}
To estimate $w_4$ we use (\ref{cdtc}) for the first two factors
combined with (\ref{chidiffh}) for the last.  For $w_5$ we combine
(\ref{cdtc}) with a symbol bound for the rest.
\end{proof}

Now the proof of estimate \eqref{pk} and hence the proof of 
Proposition~\ref{precise} are  complete. 
\end{proof}

\section{The dispersive estimates}
\label{six}
In this section we prove Theorems \ref{tpq1},\ref{tpq},\ref{tpdq}
using the parametrices constructed in the previous section.  By Lemma
\ref{co} it suffices to prove the estimates \eqref{estrsk} and
\eqref{estrspk} for the range of exponents given in the Theorems.

\subsection{The parametrices.} 
In order to prove Theorems \ref{tpq1},\ref{tpq},\ref{tpdq} we need to
use the parametrix in Proposition~\ref{pari} or the one in
Proposition~\ref{precise}. In either case we have to verify that the
estimate (\ref{comut}) holds.   In what
follows we assume that the operator $P$ is in cannonical form,
\[
P = D_t +a^w(t,x,D_x) +i b^w(t,x,D_x)
\]

\begin{lemma}\label{division}  Suppose that $p$ is in cannonical form
  and that either {\bf (A1)} and {\bf (A2)'} hold or {\bf (A1)'}
  holds. Then the following fixed time estimate is valid:
\[ 
\|[D_t+a^w,b^w] u\|_{L^2} \lesssim \|b^w u\|_{L^2} + \|u\|_{L^2}  
\]
\label{comt}\end{lemma} 

\begin{proof} 
  We first show that {\bf (A1)} and {\bf (A2)'} imply {\bf (A1)'}.
  Suppose that $ p = \tau + a(t,x,\xi) $ and $q= b(t,x,\xi)$. The
  principal normality condition {\bf (A1)} takes the form
\[ 
 |\{ \tau + a , b\}| = |\{ a,b\} + b_t| \lesssim |\tau + a| + |b| + 1 
\]
Setting $\tau = - a$ this  reduces to 
\begin{equation} 
 |\{ \tau + a , b\}| = |\{ a,b\} + b_t| \lesssim  |b| + 1 
\end{equation} 
If $\tau + a, b \in \l \Sl$ then $\{ \tau + a , b\} \in \l
S_\lambda^1$.  Since $b$ is of principal type the stronger condition
{\bf (A1)'} holds by a simple division argument.

Suppose now that {\bf (A1)'} holds. Then
\begin{eqnarray*}
\{ a,b\} + b_t &=& r_1(t,x,\tau,\xi)+ r_2(t,x,\tau,\xi)
(\tau+a(t,x,\xi)) \\ &+&
r_3(t,x,\tau,\xi) b(t,x,\xi) + r_4(t,x,\tau,\xi)
\end{eqnarray*}
where
\[
  |r_1(t,x,\xi,\tau)| \lesssim |\tau+a(t,x,\xi)|+ |b(t,x,\xi)| + 1 
\]
The left hand side is independent of $\tau$. Setting
$\tau=-a(t,x,\xi)$ on the right we obtain a similar relation of the
form
\[
\{ a,b\} + b_t = r_1(t,x,\xi) +
r_3(t,x,\xi) b(t,x,\xi) + r_4(t,x,\xi)
\]
where
\[
  |r_1(t,x,\xi)| \lesssim |b(t,x,\xi)| + 1, \qquad r_1 \in \l \Sl 
\]
Hence, after rescaling, by Theorem 3 of \cite{MR1944027}, we obtain
\[ 
\|r_1^w(t,x,D)\|_{L^2} \lesssim \| b^w(t,x,D) u\|_{L^2} + \|u\|_{L^2}. 
\]
Moreover the operator
\[ 
(r_3 b)^w(t,x,D) - r_3^w(t,x,D) b^w(t,x,D) 
\]
is bounded in $L^2$. This implies the conclusion of  Lemma~\ref{comt}.
\end{proof}

\subsection{The general case: Proof of Theorem \ref{tpq1}}

Recall that in this case the cannonical form of $p$ is
\[
p = \tau +a(t,x,\xi) + \alpha b(t,x,\xi) + i b(t,x,\xi) \qquad \alpha
\in \R
\]
where the symbol $ \tau +a(t,x,\xi)$ satisfies the curvature condition
{\bf (A2)}.

We  prove that we can find a parametrix $K$ which satisfies
(\ref{estrsk}) and (\ref{estrspk}).  We choose a second cutoff multiplier 
 $\tilde \chi$, identically $1$ in the support of $\chi$.  Then we define
\[
K(t,s) = \tilde \chi^w(x,D)  H(t,s).  
\]
where $H$ is the parametrix in Proposition~\ref{precise}, modified 
as described in Remark~\ref{calpha} if $\alpha \neq 0$.

Given the estimates (\ref{k}) and (\ref{pk}), in order to prove
(\ref{estrsk}) and (\ref{estrspk}) it suffices to show that
for $(r,s)$ satisfying (\ref{rs}) we have
\begin{equation}
\chi^w: \l^{-\rho} L^{q'}L^{r'} \to (V^2_A)^*,
\qquad \tilde \chi^w:  V^2_A \to \l^{\rho} L^{q}L^{r} ,
\label{embed}\end{equation}
where $A = a^w$.  These are dual estimates, and Lemma \ref{v2a} below
asserts that they are a consequence of the dispersive estimates for
$D_t + a^w$.

\begin{prop}
  Let $q > 2$. Suppose that the estimate \eqref{l2l2e} in
  Corollary~\ref{l2l2} holds for the operator $D_t +a^w$. Then for any
  symbol $\chi \in \Sl$ with support in $B_\l$ we have the following
  microlocal embeddings:
\begin{equation}
\chi : V^2_A \to \l^{\rho(r,s)} L^q L^r, \qquad
\chi:\l^{-\rho(r,s)} L^{q'} L^{r'}
\to (V^2_A)^{*} 
\end{equation}
\label{v2a}\end{prop}

\begin{proof}
It suffices to prove the first embedding, the second follows by
duality. We do not use the full strength of (\ref{estrs}), instead
we only use $L^2$ norms on the right hand side. Using the cannonical
form of $p$ we write it as
\[
\|\chi^w(t,x,D_x) u\|_{\l^{\rho(q,r)} L^{q}L^{r}} \lesssim 
\| (D_t+a^w(t,x,D_x)) u\|_{L^2} +\|u\|_{L^2}
\]
We further specialize this to solutions to the homogeneous equation,
\[
D_t+a^w(t,x,D_x)) u=0, \qquad u(0) = u_0
\]
for which the $L^2$ norm of the solutions
is preserved in time. Then
\begin{equation}
\|\chi^w(t,x,D_x) u\|_{\l^{\rho(q,r)} L^{q}L^{r}} \lesssim
\|u_0\|_{ L^2}
\label{hmg}\end{equation}

We define the atomic space $U_A^q \subset L^\infty$ whose atoms have
the form
\[
u = \sum_{(t_k) \in \TT} 1_{[t_k,t_{k+1})} S(t,0) u_k, \qquad
\sum_k \|u_k\|_{L^2}^q = 1 
\]
Thus the atoms are step functions where each step is a solution 
to the homogeneous equation. 

\begin{lemma}
Suppose that \eqref{hmg} holds. Then 
\[
\chi^w: U_A^q \to \l^{\rho(q,r)} L^{q}L^{r}
\]
\end{lemma}
The proof of the Lemma is straightforward. It suffices to prove it for
each $U_A^q$ atom. But then we  apply \eqref{hmg} to each
step of the atom and then sum up the $q$'th power of the results.

To conclude the proof of Proposition~\ref{v2a} we still need 
a second result, namely

\begin{lemma}
Suppose that $q > 2$. Then $V^2_A \subset U^q_A$.
\end{lemma}

\begin{proof}
We can conjugate by the evolution operator $S(t,s)$ associated  to $A$
and reduce the problem to the case when $A=0$. Hence we replace
$V^2_A$ by $V^2$ and $U^q_A$ by $U^q$. Recall also that according
to our convention, all $V^2$ functions are right continuous. The same
holds for all $U^q$ functions because each atom is right continuous.

Let $u \in V^2$ with norm $1$. For each nonnegative integer $j$ we
inductively construct functions $u_{j}$ and $v_j$ and a finite
disjoint partition $\I_j$ of the time interval $[0,1]$ with the
following properties:

\begin{enumerate}
\item  The functions $u_j$ are right continuous, $\|u_j\|_{L^\infty} \leq 2^{-j}$ and
$u - u_j$ is constant on any interval $I \in \I_j$. \label{item1}
\item The functions $v_j$ are right continuous step functions associated to the
partition $\I_j$. \label{item2}
\item We have $u_{j+1} = u_j - v_{j+1}$.\label{item3}
\item For each $j$, $I_{j+1}$ is a subpartition of $\I_j$. \label{item4}
\end{enumerate} 

This partition is constructed as follows. We initialize $u_0=u$,
$v_0=0$, $\I_0 = \{[0,1]\}$. It remains to do the inductive step.
Suppose we have $u_j$ and $\I_j$.  We partition each interval $I \in
\I_j$ according to the following criteria. Begin with the left
endpoint $t_I^0$. Then choose the next point $t_I^1$ minimal with the
property that $\|u_j(t_I^0) - u_j(t_I^1)\| \geq 2^{-j-1}$, and continue
until no such point can be found (i.e. we have reached the right end
of $I$). This process ends after finitely many steps, as (i) shows that $u_j
\in V^2(I)$.

The finer partition of $[0,1]$ obtained in this way is denoted by 
$\I_{j+1}$. The function $v_{j+1}$ is defined by
\[
v_{j+1}(t) = u_j(t_I^k), \qquad t \in [t_I^k,t_I^{k-1}]
\]
Then we set 
\[
u_{j+1} = u_j - v_{j+1}
\]

It is clear that the properties \eqref{item1}-\eqref{item4} are satisfied by
construction. By \eqref{item1} and \eqref{item3} we obtain the representation
\[
u = \sum_{j=1}^\infty v_j
\]
which converges in $L^\infty$ since 
\[
\|v_j\|_{L^\infty} \lesssim \|u_{j-1}\|_{L^\infty} +
\|u_j\|_{L^\infty} \lesssim 2^{1-j} + 2^{-j}
\]

Next we measure $v_j$ in $U^q$ as multiples of atoms.
We obtain
\[
\|v_j\|_{U^q} \lesssim 2^{-j} n_j^{\frac{1}q}
\]
where $n_j$ is the number of intervals in $I_j$.  To estimate $n_j$ we
compute the 2-variation of $u$ with respect to the $\I_j$ partition.
This is where we use the inductive choice of the partition.
Precisely, take $I \in \I_{j-1}$. By \eqref{item1} we know that for
all intervals $J \in \I_j$ with $J \subset I$, except possibly for the
last one, the variation of $u$ between its endpoints equals the
variation of $u_j$ between its endpoints, which is at least $2^{-j}$.
Hence we obtain
\[
1 = \|u\|_{V^2}^2 \geq (n_j - n_{j-1}) 2^{-2j}
\]
Therefore $n_j - n_{j-1} \leq 2^{2j}$, which after summation leads
to $n_j \lesssim 2^{2j}$. Going back to $v_j$ this yields
\[
\|v_j\|_{ U^q} \lesssim 2^{(\frac{2}{q}-1)j},
\]
which in turn implies that
\[
\|u\|_{U^q} \lesssim 1
\]
\end{proof} 

This concludes the proof of Proposition~\ref{v2a}.
\end{proof}

\subsection{ The involutive case: Proof of Theorem \ref{tpq}}

We begin with a discussion of the geometric conditions.  The condition
{\bf (A3)'} guarantees that for each $(x,\xi) \in \charpq$ there exist
real $\alpha,\beta$ such that the Hessian $\d_\xi^2 (\alpha
p_{re}(x,\xi) +\beta p_{im}(x,\xi))$ restricted to $T \charpq_x$ has
rank at least $n-2-k$. For operators in cannonical form this says that
if $b(t,x,\xi)=0$ then the Hessian $\d_\xi^2 (\alpha a(t,x,\xi) +\beta
b(t,x,\xi))$ restricted to the orthogonal complement of
$b_\xi(t,x,\xi)$ has rank at least $n-2-k$.  The stronger condition
{\bf (A5)'} says that for operators in cannonical form the same holds
with $\alpha = 0$ and $\beta=1$. This is the same as saying that the
characteristic set of $b$ has at least $n-2-k$ nonvanishing
curvatures. This will allow us to use Theorem~\ref{tp} for $b^w$.

Here the dispersive estimates will follow using the simpler parametrix
defined in (\ref{h}), the $L^2$ estimates of Proposition~\ref{pari}
and the dispersive estimates for operators with real symbols in
Theorem~\ref{tp}.

We consider a multiplier $\tilde \chi$ supported in $B_\l$ and whose
symbol equals $1$ near the support of $\chi$.  Then we define the
localized parametrix $K$ by
\[
K = \tilde \chi^w H.
\]
where $H$ is defined in (\ref{h})
with $A=a^w$ and $B=b^w$. We begin with fixed time estimates.

\begin{prop}
Assume that {\bf (A3-5)'} hold. Then the  parametrix $K$ defined above
satisfies the bounds
\begin{eqnarray}
\|K(t,s)\|_{L^2 \to L^2} &\lesssim& 1 \nonumber\\
\|K(t,s) \chi^w\|_{L^{\frac{2(n-k)}{n-k+2}} \to L^2} &\lesssim& |t-s|^{-\frac12}
\l^{\frac{n+k-2}{2(n-k)}}, \nonumber \\  \|  K(t,s)\|_{L^2 \to
  L^{\frac{2(n-k)}{n-k-2}}} &\lesssim& |t-s|^{-\frac12} \l^{\frac{n+k-2}{2(n-k)}}
\label{ksp1}
\\ 
\| K(t,s) \chi^w\|_{L^{\frac{2(n-k)}{n-k+2}} \to L^{\frac{2(n-k)}{n-k-2}}} 
&\lesssim & |t-s|^{-1}\l^{\frac{n+k-2}{n-k}}
\nonumber \end{eqnarray}
\begin{equation}
\| [(D_t+A+iB) K] (t,s) \chi^w  \|_{L^{\frac{2(n-k)}{n-k+2}} \to L^2}
 \lesssim |t-s|^{-1/2} \l^{\frac{n+k-2}{2(n-k)}}
\label{ksp3}\end{equation}
\label{pos}\end{prop}
\begin{proof}
  For each $t$ and each $1 \leq \mu \leq \l^\frac12$ we define the
  Hilbert space
\[
X_\mu(t) = \{ u \in L^2;\ Bu \in L^2\}, \quad \|u\|^2_{X_\mu(t)} =
\mu \|u\|_{L^2}^2 +  \mu^{-1} \|Bu\|_{L^2}^2 
\]
Its dual is given by
\[
X_\mu(t)^{*} = \{ u = f_1+ B(t) u_2; \ f_1, f_2  \in L^2\},
\]
\[
\|f\|_{X_\mu(t)^{*}}^2  = \inf_{f=f_1+B(t) f_2}  \mu^{-1}
\|f_1\|_{L^2}^2 + \mu \|f_2\|_{L^2}^2
\]
Set $\mu = |t-s|^{-\frac12}$. The $L^2 \to L^2$ estimates for $H$ in
Proposition~\ref{pari} lead to
\[
\|H(t,s)\|_{L^2 \to X_\mu(t)} \lesssim |t-s|^{-\frac12}, \quad
\|H(t,s)\|_{X_\mu(t)^* \to L^2} \lesssim |t-s|^{-\frac12},
\]
\[
\|H(t,s)\|_{X_\mu(s)^* \to  X_\mu(t)}
\lesssim |t-s|^{-1}
\] 
A short computation using Lemma~\ref{l0} (see also Lemma~\ref{co})
also shows that
\[
\| [(D_t+A+iB) \tilde \chi^w H](t,s)\|_{X_\mu(t)^* \to L^2}  \lesssim |t-s|^{-\frac12},
\]
Then the conclusion of the proposition follows from 
the next lemma:

\begin{lemma}
Assume that {\bf (A3-5)'} hold. Then 
\begin{equation} 
\chi^w: X_\mu \to \l^{\frac{n+k-2}{2(n-k)}}  L^{\frac{2(n-k)}{n-k-2}}
\end{equation}
\end{lemma}

\begin{proof}
 We  need to prove the estimate
\[
\l^{-\frac{n-1}{n-k}}  \|\chi^w u\|_{L^\frac{2(n-k)}{n-k-2}} 
\lesssim \mu
\|u \|_{L^2} + \mu^{-1} \|B(t) u\|_{L^2} \quad 1 \leq \mu \leq \l^{\frac12}
\]
We can localize spatially on the $\mu^{-2}$ scale.  Then we rescale
back to scale $1$. After doing this we have reduced the problem to a
similar problem but with $\mu = 1$ and $B(x,\xi):= \mu^{-2} B(\mu^{-2}
x, \mu^2 \xi)$.  Then the conclusion follows from Theorem~\ref{tp}
with $\l:= \mu^{-2} \l$.
\end{proof}
This completes the proof of Proposition~\ref{pos} and hence the proof
of Theorem~\ref{tpq}.
\end{proof}

\subsection{The degenerate involutive case: Proof of Theorem \ref{tpdq} }
\label{degenerate}

As before, we consider a multiplier $\tilde \chi$ supported in $B_\l$ and whose
symbol equals $1$ near the support of $\chi$.   Then we define the
localized parametrix $K$ by
\[
K = \tilde \chi^w H
\]
where $H$ is the better parametrix introduced in (\ref{param1}) with
$A=a^w$ and $B=b^w$.  Because of Lemma~\ref{comt} and
Theorem~\ref{exact}, this is a good $L^2$ parametrix. We need to show
that it satisfies (\ref{estrsk}) and (\ref{estrspk}).

Our strategy is as follows.  Due to the presence of the small
parameter $\dd$, one expects that most of the kernel of the parametrix
$K$ is concentrated in phase space in a small angular neighbourhood 
of the Hamilton flow for $D_t +a^w$. Our assumption {\bf (A6)'} shows
that this is not a degenerate decay direction, therefore this part
of the parametrix should satisfy good pointwise estimates.
The rest of the parametrix, on the other hand, may contain bad decay
directions. However, it has the reeedeming quality that it is small,
i.e. satisfies the same bounds as the parametrix we have constructed
in the nondegenerate case.

As in the nondegenerate case, we begin with a discussion of the
geometric assumptions.  From {\bf (A3)'} we know that, given
$(t,x,\xi)$ with $b(t,x,\xi)=0$,  there are real $\aa,\bb$ so that
the Hessian $\d_\xi^2 (\alpha a(t,x,\xi) +\beta b(t,x,\xi))$
restricted to the orthogonal complement of $b_\xi(t,x,\xi)$ has rank
at least $n-2-k$.  The condition {\bf (A5)'} says that the
same holds with $\alpha = 0$ and $\beta=1$. {\bf (A6)'}, on the 
other hand, says that we can also choose $\aa = 1$ and $\bb = 0$.
In other words, $\d_\xi^2 a(t,x,\xi)$ must have
rank $n-2-k$ on the orthogonal complement of $b_\xi$.

\subsubsection{ A pointwise fixed time bound.}

Here we set up the first building block of our estimates for $K$,
namely the pointwise estimates in directions which are close
to the Hamilton flow of $D_t + a^w$. In the following proposition
the notation $B^{1,1}_1$ stands for a Besov space.

\begin{prop}
Let $\e > 0$, small. Let $\kappa^l, \kappa^r$ be symbols whose Fourier transforms 
satisfy 
\[
\|\hat \kappa^l\|_{B^{1,1}_1}
\leq 1, \qquad \|\hat \kappa^r\|_{L^1} \leq 1, \qquad \text{supp }
\hat \kappa^l, \hat \kappa^r \subset [-\e,\e]
\]
Then 
\begin{equation}
\|\tilde\chi^w \kappa^l((t-s)b^w(t)) S(t,s) \kappa^r((t-s)b^w(s))
\chi^w \|_{L^1 \to L^\infty} \lesssim     
\l^{\frac{n+k-2}{2}}
 |t-s|^{-\frac{n-k}{2}}   
\label{simest}\end{equation} 
\label{sim}\end{prop}
\begin{proof}
  We change both the temporal and spatial scale, exactly as in the
  proof of Proposition~\ref{fixedtime}. In the new rescaled setting
  the spatial and frequency scales are both equal to
\[
  \mu = \sqrt{|t-s| \lambda},
\]
$s$ becomes $0$ and $t$ becomes $1$.
The  estimate \eqref{simest} changes to
\begin{equation}
\|\tilde \chi^w \kappa^l(b^w(1)) S(1,0) \kappa^r(b^w(0))\chi^w \|_{L^1 \to L^\infty}
\lesssim     
\mu^{k-1}
\label{resceq}\end{equation}
The rescaled symbols $a,b$ satisfy \eqref{a11}, and in addition
\begin{equation}
 |b_x| + |b_\xi| \lesssim \mu
\label{bxxi}\end{equation}
The symbols $\tilde \chi$ and $\chi$ are compactly supported inside
a ball
\[
B_\mu = \{ |x| \leq \mu, \ |\xi| \leq \mu\}
\]
and are smooth on the scale of $B_\mu$.

We rewrite the above operator in the form
\[
 \int   \hat{\kappa}^l(\theta) \tilde \chi^w 
e^{i \theta b^w(1)} S(1,0)  e^{ih b^w(0)}   \chi^w \hat{\kappa}^r(h)  \, d
\theta\, dh
\]
Since $ \hat{\kappa}^r$ is integrable 
we can neglect the $h$ integration and seek a bound for
\[
H_h(t,s) = \int  \hat{\kappa}^l(\theta)  \tilde\chi^w
e^{i \theta b^w(1)} S(1,0)  e^{ih b^w(0)} \chi^w  d \theta
\]
To obtain pointwise bounds we need some representation for the kernel
$W_\theta (y,\tilde y)$ of $ \tilde\chi^w e^{i \theta b^w(1)} S(t,s)
e^{i h b^w(0)} \chi^w$.  In the phase space this means we move from
$(x,\xi)$ to $(x_h,\xi_h)$ along the $D_h- b^w(0)$ flow, to
$(x_{h,1},\xi_{h,1})$ along the $D_t +a^w$ flow and then further to
$(x_{h,1,\theta}, \xi_{h,1,\theta})$ the Hamilton flow for $D_\theta -
b^w(1)$. This motivates the following

\begin{lemma}
  Let $a,b$ be real symbols satisfying \eqref{a11} and \eqref{bxxi}.
  The kernel of $\tilde \chi^w e^{i \theta b^w(1)} S(1,0) e^{i h
    b^w(1)} \chi^w$ has the form
\[
W_\theta (y,\tilde y) =\int_{B_\mu} G(\theta,x,\xi,
y)e^{i\Psi}
e^{-\frac12(\tilde y-x)^2}  e^{i\xi(x-\tilde y)}
\,dx\, d \xi +O(\mu^{-\infty})
\]
with
\begin{eqnarray*}
\Psi &=& - \xi_{h,1,\theta}(x_{h,1,\theta}-y) + \xi(x-\tilde y)
+ \psi_{-B(0)}(h,x,\xi) \\&& +\
\psi_A(1,x_{h},\xi_{h})+
 \psi_{-B(1)}(\theta, x_{h,1},\xi_{h,1}) 
\end{eqnarray*}
where $G$ is smooth, bounded, compactly supported in $(x,\xi) \in
B_\mu$, rapidly decreasing away from $y = x_{h,1,\theta}$ and
\begin{equation}
\left| (y - x_{h,1,\theta})^\bb  \frac{\d^\aa}{\d\theta^\aa} G(\theta,x,\xi,
y)\right| \lesssim c_{\aa,\bb} \mu^{|\aa|}
\label{gtheta}\end{equation}
\end{lemma}

\begin{proof}
We interpret $e^{i \theta b^w(1)} S(1,0) e^{i h b^w(1)}$ as a single
evolution operator where the generator is $\pm b^w(0)$ up to time $0$,
$a^w(t)$ for $t \in [0,1]$ respectively $\pm b^w(1)$ beyond time $1$.
Then the representation is the one given by Proposition~\ref{exact}. In
addition, in order to obtain \eqref{gtheta} we also need the
supplimentary result in Proposition~\ref{tdif}, which we can apply
due to \eqref{bxxi}. The role of the operators $\tilde \chi$ and
$\chi$ is simply to restrict the nontrivial part of $G$ to a compact 
subset of $B_\mu$. 
\end{proof}

We return to the proof of Proposition~\ref{sim}. We have
\[
H_h (y,\tilde y) =  \int \hat{\kappa}^l(\theta) W_\theta
(y,\tilde y) d\theta
\] 
and we want to prove that
\begin{equation}
|H_h (y,\tilde y) | \lesssim \mu^{k-1}
\label{hhyy}\end{equation}
Here $\hat{\kappa}^l(\theta) \in B^{1,1}_{1}$, with support in
$[-\e,\e]$. But $B^{1,1}_1$ can be thought of as an atomic space where
an atom $\omega_j$ is a bounded function supported in an interval $I$ of
size $2^{-j}$ and which is smooth on the same scale. The size of $j$ 
is limited by the fact that $I$ must be contained in $[-\e,\e]$.
Without any restriction in generality we can assume that
$\hat{\kappa}^l(\theta)$ is such an atom,
\[
\hat{\kappa}^l(\theta) = \omega_j(\theta), \qquad \text{supp }
\omega_j \subset [\theta_0-2^{-j}, \theta_0+2^{-j}] \subset [-\e,\e]
\]
We need to consider three cases depending on the size of $j$.

\bigskip

{\bf I. The case $2^{j} \leq \mu$.}
Then we first need to integrate by parts with respect to $\theta$.
The derivative of the phase with respect to $\theta$ is
\begin{eqnarray*}
&\displaystyle  \frac{d}{d\theta}[ \xi_{h,1,\theta}(x_{h,1,\theta}-y) +
\psi_{-B(1)}(\theta,x_{h,1},\xi_{h,1})] =&\\ & b(x_{h,1},\xi_{h,1}) +
O(\mu)(1+|(x_{h,1,\theta}-y)|)&
\end{eqnarray*}
while for higher order derivatives
we get
\[
\left|\frac{d^\aa}{d\theta^\aa} [ \xi_{h,1,\theta}(x_{h,1,\theta}-y) +
\psi_{-B(1)}(x_{h,1},\xi_{h,1})] \right| \leq c_\aa \mu^{|\aa|}
(1+|(x_{h,1,\theta}-y)|)
\]
Also we can use \eqref{gtheta} for $G$.
Hence in the region where $|b| \gg \mu$ we can integrate by parts
and get a rapidly decaying contribution. At this point we simply take
absolute values and write
\begin{eqnarray*}
|H_h(y,\tilde y) | \! \lesssim \! \int_{B_\mu} \!\!\!\!\!\!&&
 ( 1+ \mu^{-1}|b(x_{h,1},\xi_{h,1})|)^{-N}
(1+|y-x_{h,1,\theta}|)^{-N} \\ && (1+|\tilde y-x|)^{-N} 
dx d\xi d\theta
\end{eqnarray*}
Since $x_{h,1,\theta}$ is a Lipschitz function of $x$, it follows that
the integration with respect to $x$ is trivial, so we get
\[
|H_h(y,\tilde y) | \lesssim \int_{B_\mu}
 ( 1+ \mu^{-1}|b(\tilde y_{h,1},\xi_{h,1})|)^{-N} 
(1+| y-\tilde y_{h,1,\theta}|)^{-N} 
d\xi d\theta
\]
Since $|\d_\xi b| \approx \mu$ in $\{b = 0\}$, the first factor
in the integrand essentially restricts $\xi$ to a neighbourhood 
of size $1$ of $\{b = 0\}$. Then we can evaluate the above integral by
a similar integral on $\{b = 0\}$,
\[
|H_h(y,\tilde y) | \lesssim \int_{\{b (\tilde y,\xi)=0\} \cap B_\mu}
(1+| y-\tilde y_{h,1,\theta}|)^{-N}  d\H^{n-2}(\xi)\, d\theta
\]
Also on the set $b = 0$ we can use the principal normality condition
to absorb $h$ into $\theta$. More precisely, we can write 
\[
\tilde y_{h,1,\theta} = \tilde y_{0,1,\theta+\tilde \theta(h,y,\xi)}
\qquad \tilde \theta \approx h
\]
Thus without any restriction in generality we assume that $h=0$.
Dropping the subscript $h$, it remains to prove that 
\[
\int_{\{b (\tilde y,\xi)=0\} \cap B_\mu} (1+| y-\tilde y_{1,\theta}|)^{-N}  d\H^{n-2}(\xi) \,
d\theta \lesssim \mu^{k-1}
\]
We also rescale $\theta$ by a $\mu$ factor, so that it varies on the
$\mu$ scale (same as for $\xi$). Then the desired estimate becomes
\[
\int_{\{b (\tilde y,\xi)=0\} \cap B_\mu}
(1+| y-\tilde y_{1,\mu^{-1} \theta}|)^{-N} d\H^{n-2}(\xi)\, 
d\theta \lesssim \mu^{k}
\]
Hence we need to investigate the map 
\[
\xi, \theta \to \tilde y_{1,\mu^{-1} \theta}   \xi \qquad
\text{in char}_{\tilde y} \  B
\]
As in the proof of Proposition~\ref{fixedtime} one shows that this is a small
Lip\-schitz perturbation of
\[
\xi,\theta  \to a_\xi + \mu^{-1} \theta b_\xi
\]
The differential of this map is given by the matrix
\[
\left( \begin{array}{cc}a_{\xi\xi} + \mu^{-1} \theta b_{\xi \xi}, &
    \mu^{-1}  b_\xi
  \end{array} \right)
\]
acting from  $b_\xi^\perp \times \R$ to $\R^n$.  Since
$\theta$ is small it suffices to study this at $\theta = 0$.  Then we
need this matrix to have rank $n-k-1$.  Equivalently, $a_{\xi\xi}$
must have rank $n-k-2$ as a quadratic form acting on the orthogonal
complement of $b_\xi$. But this follows from our assumption {\bf
  (A6)'}.

{\bf II. The case $\mu \leq 2^j \leq \mu^2$.}
For this range of $j$ we can still integrate by parts with respect to
$\theta$ but only in the region where $|b| \gg 2^j$. Neglecting all other
oscillations as well as the $x$ integration this leads to
\[
|H_h(y,\tilde y) | \! \lesssim \! \int_{I_j}
\int_{B_\mu} ( 1+ 2^{-j}|b(\tilde y_{h,1},\xi_{h,1})|)^{-N}
(1+|y-\tilde y_{h,1,\theta}|)^{-N} 
d\xi d\theta
\]
Since $|\theta - \theta_0| \leq \mu^{-1}$ in $I_j$
it follows that $| y_{h,1,\theta}- y_{h,1,\theta_0}| \lesssim 1$.
Then the $\theta$ integration is also trivial and we obtain
\begin{equation}
|H_h(y,\tilde y) | \! \lesssim \! 2^{-j} 
\int_{B_\mu} ( 1+ 2^{-j} |b(\tilde y,\xi)|)^{-N}
(1+|y-\tilde y_{h,1,\theta_0}|)^{-N} 
d\xi
\label{hhyys}\end{equation}
We split this integral into two regions, one where $b$ is small,
$|b| \ll \mu^2$, and one where $b$ is large, $|b| \gtrsim \mu^2$.  

In the first region the level sets of $b$ are nondegenerate and close
to zero level sets. Then we can use the coarea formula to reduce
\eqref{hhyy} in this region to bounds for integrals over level sets of
$b$:
\begin{equation}
\int_{\{b(\tilde y,\xi) = b_0\} \cap B_\mu} 
(1+|y-\tilde y_{h,1,\theta_0}|)^{-N} 
 d\H^{n-2}(\xi)   \lesssim \mu^k, \qquad |b_0| \ll \mu^2
\label{lowb}\end{equation}
Hence we need to consider the map
\[
\{b(\tilde y,\xi) = b_0\} \ni \xi \to \tilde y_{h,1,\theta_0}
\]
Since $h,\theta_0 \ll 1$, this is a small $\dot C^1$ 
perturbation of the map $\xi \to a_\xi$ on $\{b(\tilde y,\xi) =
b_0\}$. The differential of this map is
\[
\eta \to a_{\xi \xi} \eta, \qquad \eta \perp b_\xi
\]
Since $b_0 \ll \mu^2$, this is a small $\dot C^1$ perturbation of a similar map 
for $b_0=0$, and by {\bf (A6)'} it has rank at least $n-k-2$.
The domain of integration is an $n-2$ dimensional cube of size $\mu$, so \eqref{lowb}
follows.

It remains to consider the large values of $|b|$, i.e. where $b
\gtrsim \mu^2$. For this the right hand side in \eqref{hhyys} is
largest when $2^j = \mu^2$, in which case we need to prove that
\begin{equation}
\int_{B_\mu}  (1+|y-\tilde y_{h,1,\theta_0}|)^{-N}  d\xi   \lesssim \mu^{k+1}
\label{highb} \end{equation}
Now  we need to consider the map
\[
B_\mu  \ni \xi \to \tilde y_{h,1,\theta_0}.
\]
Since $h,\theta_0 \ll 1$, this is a small small $\dot C^1$ 
perturbation of the map $\xi \to a_\xi$. The domain of integration is
$n-1$ dimensional, therefore we need $a_{\xi \xi}$ to have rank at
least $n-k-2$.
 
{\bf III. The case $\mu^2 \leq 2^j$.} This is the easiest case, as
there is no integration by parts. We freeze $\theta$ as in the
previous case and the estimate quickly reduces to \eqref{highb}. 
 \end{proof}

\subsubsection{ A dyadic fixed time $L^{r'} \to L^r$ bound.}

Here we show how to combine the above pointwise bounds
with the $L^{r'} \to L^r$ previously obtained in the nondegenerate
case.

\begin{prop}
Suppose that $\phi_j$, $\psi_j$ are smooth bump functions on the $2^j$
scale. Let $2 \leq r \leq \frac{2(n-k)}{n-k-2}$  and $q$ subject to \eqref{rs1}.
Then we have
\[
\|\tilde \chi^w \phi_j(b^w(t)) S(t,s) \psi_j(b^w(s)) \chi^w \|_{L^{r'} \to L^{r} }
\!\lesssim\! \l^{2\rho(q,r)}
2^{ \frac{2j}q} \! (1+2^j|t-s|)^{-\frac{2(n-k-2)} 
{q(n-k)}}
\]
\label{ftrr}\end{prop}

\begin{proof}
  Let $\e > 0$, small.  Let $\kappa$ be a smooth function supported in
  $[-\e,\e]$ which equals $1$ near the origin.  We define modified
  functions $\tilde \phi_j$, $\tilde \psi_j$ by
\[
\hat{\tilde{\phi_j}}(\theta) =
\left\{ \begin{array}{lc} 0 & |t-s| < 2^{-j} \cr \cr
\kappa(|t-s|^{-1} \theta) \hat{\phi_j}(\theta) & \text{otherwise}
\end{array} \right.
\]
and similarly for $\tilde \psi_j$. The difference satisfies
\[
|\phi_j(b) - \tilde \phi_j(b)| \lesssim (1+ 2^{-j}|b|)^{-N} (1+2^j
|t-s|)^{-N}
\]
 
We want to substitute $\phi_j$, $\psi_j$ by $\tilde{\phi_j}$, 
respectively $ \tilde \psi_j$. To do this we consider the
error term
\[
E(t,s)= \phi_j(b^w(t)) S(t,s) \psi_j(b^w(s)) - 
\tilde \phi_j(b^w(t)) S(t,s) \tilde \psi_j(b^w(s))
\]
The above difference bound easily leads  the estimates
stronger than  \eqref{ks1}, \eqref{ks2}, namely 
\[
\|E(t,s)\|_{L^2 \to L^2} \lesssim  (1+2^j |t-s|)^{-N},
\]
\[
\|b^w(t) E(t,s)\|_{L^2 \to L^2} \lesssim  2^j (1+2^j |t-s|)^{-N},
\]
\[
\|E(t,s)b^w(s)\|_{L^2 \to L^2} \lesssim  2^j (1+2^j |t-s|)^{-N}
\]
\[
\|b^w(t) E(t,s)b^w(s)\|_{L^2 \to L^2} \lesssim  2^{2j} (1+2^j |t-s|)^{-N},
\]
Using these relations the $L^{r'} \to L^{r}$ bound for $\tilde \chi^w
E(t,s) \chi^w$ follows from the curvature condition on the
characteristic set of the symbol $b$ simply by repeating the arguments
in the nondegenerate case in the proof of Proposition~\ref{pos}.

It remains to prove the estimate for the operator 
\[
\tilde \chi^w  \tilde \phi_j(b^w(t)) S(t,s) \tilde \psi_j(b^w(s)) \chi^w
\]
which is zero if $|t-s| < 2^{-j}$. For this we interpolate
the trivial $L^2 \to L^2$ bound with an $L^1 \to L^\infty$
bound derived from Proposition~\ref{sim}. If we set
\[
\tilde \phi_j(b) = \kappa^l_j ((t-s)b), \qquad \tilde \psi_j(b) = \kappa^r_j ((t-s)b),
\]
then by definition both $\kappa^l_j $ and $\kappa^r_j$ are bump functions
on the $ 2^j|t-s|$ scale and are supported in $[-\e,\e]$.
Then
\[
\| \hat \kappa^l_j\|_{L^1} \lesssim 1, \qquad \| \hat
\kappa^r_j\|_{B^{1,1}_1} \lesssim 2^j|t-s|,
\]
therefore Proposition~\ref{sim} yields
\[
\|\tilde \chi^w \tilde \phi_j(b^w(t)) S(t,s) \tilde \psi_j(b^w(s))\chi^w \|_{L^1 \to L^\infty}
\lesssim  2^j|t-s| \l^{\frac{n+k-2}{2}} |t-s|^{-\frac{n-k}{2}}   
\]
which is exactly what we need since $|t-s| > 2^{-j}$ and
the following relations hold:
\[
2\rho(q,r) = \frac{n+k-2}{2}\left(\frac{1}{r'}-\frac{1}r\right), \qquad
(n-k)\left(\frac{1}{r'}-\frac{1}r\right) = \frac{4}q.
\]
\end{proof}

\subsubsection{ A dyadic $L^2 \to L^q L^r$ bound.}
Here we use a $TT^*$ argument to derive an $L^2 \to L^q L^r$ bound
from the above fixed time bound.
\begin{prop} 
Let $\phi_j(t,b)$ be smooth bump functions on the $2^j$ scale. Let $l \leq j$.
Then
\begin{equation}
\| 1_{|t-s| < 2^{l-j}} \chi^w \phi_{j}(t,b^w(t)) S(t ,s)\|_{L^2 \to L^q L^r} \lesssim
 2^{\frac{2l}{q(n-k)}} \l^{\rho(q,r)} 
\label{eqkht}\end{equation}
\label{kht}\end{prop}

\begin{proof}
  Using a $TT^*$ argument this reduces to a bound
\[
\| \chi^w \phi_j(t,b^w(t)) S(t ,s) \phi_{j}(s,b^w(s)) \chi^w\|_{L^{q'}L^{r'} \to
L^q L^r}  \lesssim
 2^{\frac{4l}{r(n-k)}} \l^{2\rho(r,s)}
\]
where $t,s$ are restricted to a $2^{l-j}$ interval.  For this we use
the fixed time bound in Proposition~\ref{ftrr}, which yield
\begin{eqnarray*}
&&\| \chi^w \phi_j(t,b^w(t)) S(t ,s) \phi_{j}(s,b^w(s)) \chi^w\|_{L^{q'}L^{r'} \to L^q L^r}
 \\ &&\hspace{3cm}\lesssim  \l^{2 \rho(q,r)}
2^{\frac{2j}q} \left( \int_0^{2^{l-j}}   (1+2^j|t|)^{-\frac{n-k-2}{n-k}} dt\right)^{\frac{2}{q}} 
\\
&&\hspace{3cm}\approx   \l^{2 \rho(q,r)}
2^{\frac{2j}q} 2^{\frac{2(l-j)}q} 
2^{- l \frac{2(n-k-2)}{q(n-k)}} =  2^{\frac{4l}{q(n-k)}} \l^{2\rho(q,r)}.
\end{eqnarray*}
We note that in effect one could obtain better bounds for $l < 0$,
but they are not needed here.
\end{proof}

\subsubsection{The parametrix bound (\ref{estrsk}).}

We recall that the global parametrix $H$ is given by
\begin{eqnarray*}
H(t,s) &=&  \sum_j 1_{t < s} \kappa^-_j(\dd b^w(t)) S(t,s) 
 \kappa^-_j(\dd b^w(s)) e^{\dd (t-s)b^w(s)} \\ &&-\ 1_{t > s} 
\kappa^+_j(\dd b^w(t)) S(t,s) 
 \kappa^+_j(\dd b^w(s)) e^{\dd (t-s)b^w(s)}
\end{eqnarray*}
It suffices to consider the first term. We set 
\[
H_j(t,s) = \kappa^-_j(\dd b^w(t)) S(t,s) 
 \kappa^-_j(\dd b^w(s)) e^{\dd (t-s)b^w(s)}
\]
and apply Proposition~\ref{ftrr} for each $j$.
The symbols $ \kappa^-_j(\dd \cdot )$ are bump functions on the
$\dd^{-1} 2^j$ scale, while the exponential factor
provides exponential decay on the $2^{-j}$ time scale.
Hence we obtain
\begin{eqnarray*}
&&\!\!\!\!\|\tilde\chi^w H_j(t,s) \chi^w\|_{L^{r'} \to L^r} \\ && \hspace{1cm}\lesssim 
 \l^{2 \rho(q,r)}
(\dd^{-1}2^j)^{\frac{2}{q} }  (1+\dd^{-1}
2^j|t-s|)^{-\frac{2(n-k-2)}{q(n-k)}} (1+ 2^j|t-s|)^{-N}
\\ && \hspace{1cm}\lesssim \l^{2 \rho(q,r)} \dd^{-\frac{4}{q(n-k)}} 2^{\frac{2j}q} 
(2^j |t-s|)^{-\frac{2(n-k-2)}{q(n-k)}} (1+ 2^j|t-s|)^{-N}
\end{eqnarray*}

Next we sum this up with respect to $j$. For fixed $t-s$ 
the largest contribution comes from the indices $j$ which
satisfies $2^j |t-s| \approx 1$. This gives

\begin{prop}
 For $2 \leq r \leq \frac{2(n-k)}{n-k-2}$ the operator $H$ satisfies the estimate
\begin{equation}
\|\tilde \chi^w H(t,s)\chi^w\|_{L^{r'} \to L^r} \lesssim  \dd^{-\frac{4}{q(n-k)}} \l^{2\rho(q,r)}
 |t-s|^{-\frac{2}q}   
\end{equation}
\end{prop}
If we have this then we use the Hardy-Littlewood-Sobolev
inequality to obtain the $L^{q'}L^{r'} \to L^q L^q$ bounds for $H$.

Next we prove the $L^2 \to L^r L^s$ bounds for $H$. The $ L^{r'}
L^{s'} \to L^2$ bounds are essentially dual, and their proof is similar.
We decompose the operator $H$ as
\[
H = \sum_{l}  H^l
\]
where
\[
  H^l(t,s) = \sum_j 1_{2^{l-1-j} \leq |t-s| \leq 2^{l-j}} H_j(t,s)
\]
We apply Proposition~\ref{kht} to the terms on the right.
The symbols $\kappa_j^-(\dd \cdot)$ are smooth bumps on the
$\dd^{-1} 2^j$ scale, and the exponential provides rapid decay
on the $2^{-j}$ time scale. Then we obtain
\[
\|1_{2^{l-1-j} \leq |t-s| \leq 2^{l-j}} \chi^w H_j(t,s)\|_{L^2 \to L^q L^r}
\lesssim  \dd^{-\frac{2}{q(n-k)}}    2^{\frac{2l}{q(n-k)}} \l^{\rho(q,r)} (1+2^l)^{-N}
\]
We sum this with respect to $j$ using the fact that
we have orthogonality which gives square summability in $j$ on the
$L^2$ side, while on the $ L^q L^r$ side we only need $l^q$ 
summability in $j$ because the functions we add live in disjoint time
intervals. Thus we obtain
\[
\| \chi^w H^l \|_{L^2 \to L^q L^r}
\lesssim  \dd^{-\frac{2}{q(n-k)}}   2^{\frac{2l}{q(n-k)}} \l^{\rho(q,r)} (1+2^l)^{-N}
\]
Adding these estimates with respect to $l$ we finally obtain
\[
\| \chi^w H\|_{L^2 \to L^q L^r}
\lesssim  \dd^{-\frac{2}{q(n-k)}}   \l^{\rho(q,r)}.
\]

\subsubsection{The error estimates \eqref{estrspk}.}

Here we prove that for each $t \in [0,1]$ the parametrix $H$ satisfies
the error estimates
\begin{equation}
\| (D_t +a^w +i \dd b^w)H(t ,s) \tilde \chi\|_{L^{q'}L^{r'} \to L^2} \lesssim
 \dd^{ - \frac{2}{q(n-k)}}
\l^{\rho(q,r)}
\label{err}\end{equation}

\begin{proof}
It suffices to look at the forward part of the error. Then for $t > s$
we denote
\[
E(t,s)= -i (D_t -a^w +i \dd b^w)H(t,s) =  \sum_j E_j(t,s) 
\]
where
\[
E_j(t,s)= \dd \kappa_j^-(b^w(t))  (b^w(t) S(t,s)- S(t,s) b^w(s))
\kappa^-_j(b^w(s)) e^{\dd (t-s)b^w(s)}
\]
Let $\tilde \kappa_j(\eta)$ be a symbol which equals $2^{-j} \eta$ 
in the support of $\kappa_j^-$ and has slightly larger support.
Then
\begin{eqnarray*}
E_j(t,s)& = &   2^j \kappa_j^-(\dd b^w(t))  
(\tilde \kappa_j(\dd b^w(t)) S(t,s)-
S(t,s) \tilde \kappa_j(\dd b^w(s))) \\ && 
 \kappa^-_j(\dd b^w(s)) e^{\dd (t-s)b^w(s)} 
\end{eqnarray*}
which we rewrite as
\begin{eqnarray*} 
E_j(t,s)& = &  \int_s^t  2^j \kappa_j^-(\dd b^w(t))S(t,h) [D_t+a^w,\tilde
\kappa_j(\dd b^w)](h) S(h,s) 
\\ & &  \kappa^-_j(\dd b^w(s)) e^{(t-s)\dd b^w(s)}dh
\end{eqnarray*}
We further split the $E_j$'s into dyadic pieces based on the distance $t-s$,
\[
E_j^l(t,s) = 1_{\{2^{l-j-1} \leq t-s \leq 2^{l-j}\}} E_j(t,s)
\]
Then we add them back interchanging the order of summation,
\[
E^l = \sum_{j > l} E_j^l, \qquad E = \sum_l E^l
\]

To obtain an  $L^{q'}L^{r'} \to L^2$ bound for $1_{t > s} E_j^l(t,s) \chi^w$
we first use the fact that the operators $ \kappa_j^-(\dd b^w(t))$,
$S(t,h)$ are $L^2$ bounded. By \eqref{chidiffh}  the operator
$[D_t+a^w, \tilde \kappa_j(\dd b^w)(h)]$ is also bounded.
In addition, the $h$ integration occurs on an interval of size
$2^{l-j}$. Then we have
\begin{eqnarray*}
&& \| E_j^l(t,\cdot)\chi^w\|_{L^{q'}L^{r'} \to L^2}
\\ && \hspace{1cm}
 \lesssim 2^l \sup_{h<t} \|  1_{s \in I^j(t,h)}  S(h,s) 
\kappa^-_j(\dd b^w(s)) e^{(t-s)\dd b^w(s)} \chi^w\|_{L^{q'}L^{r'} \to L^2}
\end{eqnarray*}
where
\[
I^j(t,h) = [t-2^{l-j-1}, \min\{ t-2^{l-j},h\}]
\]
For the term on the right we use the dual of the estimate \eqref{kht}.
The symbol $\kappa^-_j(\dd \eta) e^{(t-s)\dd \eta}$ is a smooth 
bump function on the $\dd^{-1} 2^j$ scale, and whose size
can be bounded by $(1+2^j|t-s|)^{-N}$. Hence we obtain
\begin{eqnarray*}
& \|  1_{s \in I^j(t,h)}  S(h,s) 
\kappa^-_j(\dd b^w(s)) e^{(t-s)\dd b^w(s)} \chi^w\|_{L^{q'}L^{r'} \to L^2} &
\\ &\lesssim \dd^{-\frac{2}{q(n-k)}} \l^{\rho(q,r)} (1+ 2^l)^{-N}&
\end{eqnarray*}
which implies that
\[
\|E_j^l \chi^w\|_{L^{q'}L^{r'} \to L^2} \lesssim   \dd^{-\frac{2}{q(n-k)}}
\l^{\rho(q,r)}  2^l (1+ 2^l)^{-N}
\]
Summing up with respect to $j$ we obtain
\[
\|E^l \chi^w\|_{L^{q'}L^{r'} \to L^2} \lesssim   \dd^{-\frac{2}{q(n-k)}}
\l^{\rho(q,r)}  2^l (1+ 2^l)^{-N}
\]
There is no loss in the summation. On one hand
the inputs come from disjoint time intervals, so we gain
an $\ell^{q'}$ summation with respect to $j$. On the other
hand the outputs are almost $L^2$ orthogonal,
so we only need an $\ell^2$ summation in $j$.

The last step is to perform the summation with respect to $l$,
which is trivial.
\end{proof}

\section{Applications to local solvability}

In local solvability problems one considers a partial differential
operator or a pseudodifferential operator $P(x,D)$ and seeks 
to find a local solution $u$ to the equation
\[
P(x,D) u =f
\]
for $f$ with sufficiently small support. Local solvability is known to
hold for principally normal pseudodifferential operators with
$S_{1,0}$ type symbols, see H\"ormander~\cite{MR24:A434} and also for low
regularity symbols of type $C^2 S_{1,0}$, see Tataru~\cite{MR1944027}. We first
state the $L^2$ solvability result proved in \cite{MR1944027}.

\begin{theorem} Let $p \in C^2 S^1_{1,0}$ be a principally normal 
pseudodifferential operator. Then $P(x,D)$ is locally solvable with loss
of one derivative, in the sense that for sufficiently small $\e > 0$,
any ball $B_\e$ of radius $\e$ and any $f \in L^2$ supported in $B_\e$ 
there is $u \in L^2(\R^n)$ so that $P(x,D) u = f$ in $B_\e$.
\label{l2solve}\end{theorem} 

By duality this theorem reduces to proving an $L^2$ bound from below
for the adjoint operator, namely
\begin{equation}
\| v\|_{L^2} \lesssim \e \|P(x,D)^* v \|_{L^2}, \qquad \text{ supp } v
\subset B_\e
\end{equation}
This estimate is stable with respect to $L^2$ bounded perturbations of
$P(x,D)$.  An immediate consequence of it is the following

\begin{cor}
Under the same assumptions as in Theorem~\ref{l2solve}, local
solvability holds for $P(x,D)+V$ for any potential $V \in L^\infty$.
\label{vlinf}\end{cor}

In this section we consider principally normal pseudodifferential
operators $P$ and use geometric information about their
characteristic sets in order to derive local solvability for operators
of the form $P(x,D)+V$ where $V$ is an unbounded potential.
Our main result is

\begin{theorem} a) Let $q$, $\rho(q)$  be as in \eqref{r}. Let $k \geq
  0$ and $p \in C^2S^{1+2\rho(q)}_{1,0}$ be a symbol whose restriction
  to frequency $\l$ satisfies (A1)',(A2),(A3).  Let $V \in L^{s}$
where 
\[
\frac{1}s = \frac{1}{q'}-\frac{1}q
\]
Then $P(x,D) +V$ is locally solvable with loss of one derivative in
the sense that for sufficiently small $\e > 0$, any ball $B_\e$ of
radius $\e$ and any $f \in H^{-\rho(q)}+L^{q'}$ supported in $B_\e$
there is $u \in H^{\rho(q)} \cap L^q$ so that $p^w u = f$ in $B_\e$.

b) The same result holds if $q$, $\rho(q)$ are as in \eqref{r1}
and the restriction of $p$ to frequency $\l$ satisfies (A1), (A2)',(A3)'.
\label{lpsolve}\end{theorem}

Results of the same type but with different assumptions have been
obtained also by Dos Santos~\cite{DS}.

\begin{remark}
Implicit in this theorem is the assumption that $P$ has order
$1+2\rho(q)$. However, this is the most interesting case. If the order
of $P$ is different then one should consider the two possibilities:

a) If the order is larger then $p \in C^2S^{1+2\rho(\tilde q)}_{1,0}$
for some $\tilde q > q$. In this situation the above result holds
with $q$ replaced by $\tilde q$. To prove it 
one  simply needs  to relax the dispersive estimates
using Sobolev embeddings.

b) If the order of $P$ is smaller then, under the same assumptions as
in Theorem~\ref{lpsolve}, one can prove that the result obtained by
formally interpolating the conclusions of Corollary~\ref{vlinf} and
Theorem~\ref{lpsolve} is true.
\end{remark}

\begin{proof}[Proof of Theorem~\ref{lpsolve}]
  We prove part (a), which uses Theorem~\ref{tpq1}. The proof of part
  (b) is similar, with the only difference that it uses
  Theorem~\ref{tpq} instead.

By duality the theorem reduces to proving a bound from below for the
adjoint operator. Our main estimate is
\begin{equation}
\| v\|_{\e^{\frac18} H^{\rho(q)} \cap L^q} \lesssim \| P(x,D)^*
v\|_{\e^{-\frac18} H^{-\rho(q)}+L^{q'}}, \qquad \text{ supp } v
\subset B_\e
\end{equation}
Multiplication by $V$ maps $L^q$ into $L^{q'}$. If $\e$ is
sufficiently small then $V$ must be small in $B_\e$. Hence we can
freely replace $P(x,D)^*$ by $P(x,D)^*+V$ and conclude the proof of
the theorem.

Another useful observation is that we can modify $P(x,D)^*$ by any
perturbation which is bounded from $H^{\rho(q)}$ into $H^{-\rho(q)}$.
We make use of this in order to truncate the symbol $p(x,\xi)$ at
frequency $\leq \sqrt \l$ with respect to the $x$ variable. After this
reduction, the frequency $\l$ part of $p$ belongs to $\l^{1+2\rho(q)}
\Sl$.

We fix a ball $B_r$ of fixed sufficiently small radius $r$, and we
assume that $B_{\e} \subset \frac12 B_r$.  We consider a locally
finite covering of the frequency space with balls
\[ 
\R^n = \bigcup_{j=0}^\infty B_j
\]
where $B_0=B(0,1)$ while for each $j > 0$ there exists some $\l > 1$
so that $B_j \subset \{|\xi| \approx \l\}$, the radius of $B_j$ is
comparable with $\l$ and Theorem~\ref{tpq1}(b) can be applied in $B_r
\times B_j$ with respect to a suitable coordinate system (which may
depend on $j$). This is always possible if $r$ is sufficiently small,
as discussed in Section~\ref{norm}.  Correspondingly we choose a
smooth partition of unity in the frequency space
\[ 
1 = \chi_0(\xi) + \sum_{j} \chi_j(\xi) 
\]

Let $\chi(x,\xi)$ be a smooth symbol supported in $B_r \times B_j$ and
which equals $1$ in $\frac12 B_r \times \text{ supp } \chi_j$.  We
also choose symbols $p_j \in \l^{1+2\rho(q)} \Sl$ which which agree
  with $p$ in $B_r \times B_j$. From Theorem~\ref{tpq1}(b) we
  obtain
\[
\| \chi^w w\|_{\l^{-\rho(q)} L^\infty L^2 \cap L^q} \lesssim
\|\bar p^w_j(x,D) w\|_{ \l^{\rho(q)}L^2 + \chi^w(\l^{\rho(q)} L^1 L^2 +
  L^{q'})} +\|w\|_{\l^{-\rho(q)} L^2}
\]
We apply the above inequality to $w = \chi_j(D) v$ with $v$ supported
in $\frac12 B_r$. After including some rapidly decreasing tails in the
last right hand side term we obtain
\begin{eqnarray*}
\| \chi_j (D) v\|_{\l^{-\rho(q)} L^\infty L^2 \cap L^q} & \lesssim& 
\|\bar p^w_j(x,D) \chi_j(D) v\|_{  \l^{\rho(q)} L^2 + \l^{\rho(q)} L^1 L^2 + L^{q'})}
 \\ && +\ \|\chi_j(D) v\|_{\l^{-\rho(q)} L^2} + \l^{-N} \|v\|_{L^2}
\end{eqnarray*}
We take a new multiplier $\tilde \chi_j$ with slightly larger support,
and which equals $1$ in the support of $\chi_j$. Replacing $\chi_j(D)$
by $ \chi_j(D)\tilde \chi_j(D)$ in the right hand side, after some
commutations we get
\begin{eqnarray*}
\| \chi_j (D) v\|_{\l^{-\rho(q)} L^\infty L^2 \cap L^q}\! &\lesssim& 
\|\chi_j(D) \bar p^w_j(x,D) \tilde \chi_j(D) v\|_{  
\l^{\rho(q)} L^2 + \l^{\rho(q)} L^1 L^2 + L^{q'}} \\ && +\
\|\tilde \chi_j(D) v\|_{\l^{-\rho(q)} L^2} + \l^{-N} \|v\|_{L^2}
\end{eqnarray*}
We can also replace $\bar p^w_j(x,D)$ first by $P_j(x,D)^*$ and then by 
$P(x,D)^*$ to obtain
\begin{eqnarray*}
\| \chi_j (D) v \|_{\l^{-\rho(q)} L^\infty L^2 \cap L^q}\! &\lesssim& 
\|\chi_j(D) P(x,D)^* \tilde \chi_j(D) v\|_{  \l^{\rho(q)} L^2 + \l^{\rho(q)}
  L^1 L^2 + L^{q'}}\\ && +\
\|\tilde \chi_j(D) v\|_{\l^{-\rho(q)} L^2} + \l^{-N} \|v\|_{L^2}
\end{eqnarray*}
Finally we drop $\tilde \chi_j$ in the last term at the expense of a rapidly 
decreasing contribution,
\begin{eqnarray*}
\| \chi_j (D) v \|_{\l^{-\rho(q)} L^\infty L^2 \cap L^q} &\lesssim&  
\|\chi_j (D) P(x,D)^* v\|_{  \l^{\rho(q)} L^2 + \l^{\rho(q)} L^1 L^2 + L^{q'}}
\\ && +\ \|\tilde \chi_j(D) v\|_{\l^{-\rho(q)} L^2} + \l^{-N} \|v\|_{L^2}
\end{eqnarray*}
At this point we use the assumption on the support of $v$.  The
kernels of $\chi (D) v$ and of $\chi(D) P(D,x)$ decay rapidly on the
$\l$ scale. Hence if $\l^{-1} < \e^{\frac12} $ then all the tails
beyond the $\e^{-\frac14}$ scale are negligible.  Then we use Holder's
inequality to turn the $L^\infty L^2$ and the $L^1 L^2$ norms into
$L^2$ norms and obtain
\begin{eqnarray*}
\| \chi_j (D) v \|_{\e^{\frac18} \l^{-\rho(q)} L^2 \cap L^q} &\lesssim& 
\|\chi_j(D) P(x,D)^* v\|_{  \e^{-\frac18} \l^{\rho(q)} L^2 + L^{q'}} 
\\ &&+\ \|\tilde \chi_j(D) v\|_{\l^{-\rho(q)} L^2} 
+\l^{-N} \|v\|_{L^2}
\end{eqnarray*}
Summing up using Littlewood-Paley theory yields
\[
\| \chi_{>\e^{-\frac12}} (D) v \|_{\e^{\frac18} H^{\rho(q)} \cap L^q} 
\lesssim \| P(x,D)^* v\|_{  \e^{-\frac18}
  H^{-\rho(q)} + L^{q'}} + \|v\|_{H^{\rho(q)}} 
\]
If the multiplier $\chi_{>\e^{-\frac12}} (D)$ were not there then for
sufficiently small $\e$ we could absorb the second right hand side
term into the left hand side and conclude the proof.  As it is, we
also need a bound for the low frequencies.  This we get since $v$ has
very small support, therefore most of its energy has to be
concentrated at high frequencies:
\[
\| \chi_{<\e^{-\frac12}} v\|_{H^\rho} \ll \|v\|_{H^\rho}, \qquad 0
\leq \rho < \frac{n}2
\]
Indeed, if $H^\rho \subset L^r$ is a sharp Sobolev embedding then
\[
\| \chi_{<\e^{-\frac12}} v\|_{H^\rho} \lesssim \e^{-\frac{\rho}2}
\|v\|_{L^2} \lesssim \e^{\frac{\rho}{2}} \|v\|_{L^r} \lesssim
\e^{\frac{\rho}{2}} \|v\|_{H^\rho}
\]
The proof is even easier if $\rho=0$.
\end{proof}

\section{Applications to unique continuation}

Consider a partial differential operator $P(x,D)$ of order $m$ in
$\R^n$. Let $\Gamma$ be an oriented hypersurface in $\R^n$, which can be
represented as a nondegenerate level set of a smooth function, $\Gamma =
\{\phi = 0\}$.  The sign of $\phi$ away from $\Gamma$  determines 
the orientation of $\Gamma$.  Denote the two sides of $\Gamma$ by $\Gamma^+
= \{ \phi > 0\}$ and $\Gamma^-=\{\phi < 0\}$. Then we define the unique
continuation property across $\Gamma$ for solutions to $P(x,D)u=0$ as
follows:

\begin{definition} 
  We say that unique continuation property across $\Gamma$ holds for
  the operator $P(x,D)$ if for each $x_0 \in \Gamma$ there exists  a
  neighborhood $V$ of $x_0$ such that the following holds: 
Let $u$ be a solution
  for $P(x,D)u = 0$ in $V$ so that $u = 0$ in $\Gamma^+ \cap V$.  Then $u
  = 0$ near $x_0$.
\end{definition}

In other words, the values of a solution $u$ to $Pu=0$ on one side of
$\Gamma$ (i.e. in $\Gamma^+$) uniquely determine its values on the other
side (i.e. in $\Gamma^-$) near $\Gamma$.  One can also reinterpret this
as an uniqueness result for the Cauchy problem for $P(x,D)$ in $\Gamma^+$
with initial data on $\Gamma$.

Whether the unique continuation property holds depends on the geometry
of the surface $\Gamma$ relative to the operator $P$. One naturally
introduces the pseudoconvexity condition to describe this.
We let $p(x,\xi)$ be the principal symbol of $P$ and introduce the notation 
\[
p_\phi(x,\xi,\tau) = p(x,\xi +i \tau \nabla \phi)
\]
\begin{definition}
We say that the surface $\Gamma$ is strongly pseudoconvex with respect to
$P$ if either

a) $P$ is elliptic and 
\begin{equation}
\{ \Re p_\phi, \Im p_\phi\}
> 0 \quad \mbox{on } T^*_\Gamma \R^n \cap \{ p_\phi=\{p_\phi,\phi\}=0 \} ,\ \tau > 0
\label{is}
\end{equation}

b) $P$ has real principal symbol and both
\begin{equation}
 \{p,\{p,\phi\}\} > 0 \qquad \mbox{ on }
T^*_\Gamma \R^n \cap \{ \{p,\phi\} = p = 0\}
\label{rsf}\end{equation}
and (\ref{is}) hold.
\label{pcs}\end{definition}

Note that the property of pseudo-convexity only depends on $\Gamma$ and its 
orientation and not on $\phi$. 
There is also a version of this which applies to principally normal
operators, but here we choose to keep things simple.
Note also that for anisotropic operators such as the heat or the
Schr\"odinger operator one has to make some obvious adjustments
in the definition of the principal symbol and of the Poisson bracket. 

The pseudoconvexity condition does not preclude surfaces from being
characteristic at least at some points. However, here we assume for
simplicity that this not is the case, namely
\[
p(x, \nabla \phi) \neq 0
\]

Now we can state the main result (see H\"ormander~\cite{MR87d:35002b} and
references therein, and also Isakov~\cite{MR94k:35070} for the anisotropic case):

\begin{theorem}
  Let $P$ be an operator with $C^1$ coefficients which is either
  elliptic or has real principal symbol.  Suppose that the oriented
  surface $\Gamma$ is strongly pseudoconvex with respect to $P$.  Then
  unique continuation across $\Gamma$ holds for $P$.
\label{H}\end{theorem}

A main tool in proving unique continuation results is provided by the
Carleman estimates. To describe them we need to introduce the notion
of pseudoconvex functions. 

\begin{definition}
We say that the function $\phi$ is strongly pseudoconvex with respect to
$P$ if either

a) $P$ is elliptic and 
\begin{equation}
\{ \Re p_\phi, \Im p_\phi\}
> 0 \quad \mbox{on }   \{ p_\phi=0 \} ,\ \tau > 0
\label{if}
\end{equation}

b) $P$ has real principal symbol and both
\begin{equation}
 \{p,\{p,\phi\}\} > 0 \quad \mbox{on }   \{ p=0 \}
\label{rf}\end{equation}
and (\ref{if}) hold.
\end{definition}

The nondegenerate level sets of pseudoconvex functions are
pseudoconvex surfaces. Conversely, any pseudoconvex surface is a
nondegenerate level set of some pseudoconvex function. 

The $L^2$
Carleman estimates below imply the above unique continuation result
via a standard argument.

\begin{theorem}[$L^2$ Carleman estimates]
  Let $P$ be an operator with $C^1$ coefficients which is either
  elliptic or has real principal symbol.  Suppose that $\phi$ is
  strongly pseudoconvex with respect to $P$ in some bounded open 
  $\Omega \subset \R^n$. Given any compact subset $K$ of $\Omega$
  there are $c,\tau_0 > 0$ so that for all
  functions $u$ supported in $K$ we have:

(a) If $P$ is elliptic:
\begin{equation}
\tau^{-1} \| \etf u\|_{H^m_\tau}^2 \leq c \|\etf P(x,D) u\|_{L^2}^2, \quad \tau \geq \tau_0 
\label{ceel}\end{equation}

(b) If $P$ has real principal symbol:
\begin{equation}
\tau \| \etf u\|_{H^{m-1}_\tau}^2 \leq c \|\etf P(x,D) u\|_{L^2}^2, \quad \tau \geq \tau_0 
\label{cerl}\end{equation}
\label{ce}\end{theorem}
Here and below, the spaces $H^k_\tau$ are defined like the usual
Sobolev but giving to $\tau$ the same weight as a derivative.
Precisely,
\[
\|u\|_{H^k_\tau} \approx \|(|D|+\tau)^k u\|_{L^2}
\]
A similar meaning is associated to the notation $H^{k,p}_\tau$.

Our interest lies in replacing the $L^2$ estimates with $L^p$
estimates. This is useful in problems with unbounded potentials. 
These can also arise as linearizations of nonlinear problems.

\begin{theorem}[$L^p$ Carleman estimates, elliptic case]
  Let $P$ be an elliptic operator of order $m$ with $C^1$
  coefficients.  Let $\phi$ be a strongly pseudoconvex with respect to
  $P$ in some compact set $\Omega \subset \R^n$. Assume that the
  characteristic set of $p_\phi$ has $n-2-k$ nonvanishing curvatures,
  and let $r,\rho(r)$ be as in (\ref{r1}).  Then there are $c,\tau_0 >
  0$ so that for all functions $u$ supported in $\Omega$ we have:
\begin{equation}
 \| \etf u\|_{\tau^{\frac14} H^m_\tau \cap
  H^{m-\frac12-\rho(r),r}_{\tau} }^2 
\leq c   \|\etf P(x,D) u\|_{\tau^{-\frac14} L^2+
  H^{\frac12+\rho(r),{r'}}_{\tau}}^2,  \quad \tau \geq \tau_0 
\label{cee}\end{equation}
\label{lpcee}\end{theorem}

\begin{theorem}[$L^p$ Carleman estimates, real case]
  Let $P$ be an operator of order $m$ with real principal symbol and
  $C^2$ coefficients.  Let $\phi$ be a strongly pseudoconvex with
  respect to $P$ in some compact set $\Omega \subset \R^n$.  Assume
  that the characteristic set of $p$ has $n-1-k$ nonvanishing
  curvatures and the characteristic set of $p_\phi$ has $n-2-k$
  nonvanishing curvatures, and let $r,\rho(r)$ be as in (\ref{r}). Then
  there are $c,\tau_0 > 0$ so that for all (smooth) functions $u$ supported in
  $\Omega$ we have:
\begin{equation}
 \| \etf u\|_{\tau^{-\frac14} H^{m-1}_\tau \cap
  H^{m-1-\rho(r),r}_{\tau} }^2 
\leq c  \|\etf P(x,D) u\|_{ \tau^{\frac14} L^2+
  H^{\rho(r),r'}_{\tau}}^2,  \quad \tau \geq \tau_0 
\label{cer}\end{equation}
\label{lpce}\end{theorem}

In most cases one can also obtain mixed norm estimates.  However, we
prefered to have simpler statements for general operators.  Later on
when we consider examples we also state mixed norm estimates for
equations where this is relevant such as the heat equation and the
Schr\"odinger equation. 

Also at this point we contend ourselves with the Carleman estimates.
We state the corresponding unique continuation statements only for the
examples we consider below, and we leave to the reader the task of
deriving the corresponding unique continuation results for other
problems of interest.

Observe that, given any estimate of the form
\[
 \| \etf u\|_{X}^2 \leq c \|\etf P(x,D) u\|_{Y}
\]
the substitution $v = \etf u$ transforms it into
\[
 \| v \|_{X}^2 \leq c \| P_\phi(x,D,\tau) v \|_{Y}
\]
Hence we need to understand the geometry of the operators $P_\phi$. 
Note that only the principal part of $P_\phi$ is important, all the
lower order terms are negligible due to the $L^2$ part of the Carleman
estimates.

The operator $P_\phi$ is an operator with complex symbol, therefore we
would like to apply our results for principally normal operators.
Hence we want $p_\phi$ to be principally normal. However, the
pseudoconvexity condition shows that this is not the case, more
precisely
\[
\{ \Re p_\phi, \Im p_\phi \} > 0 \qquad \text{in} \ \ \ \ \ \  p_\phi = 0
\]

To overcome this difficulty we use a two scale approach which begins
with the observation that the $L^2$ Carleman estimates allow a
localization to the $\tau^{-\frac12}$ scale.  Hence, on one hand we
apply our dispersive estimates on the $\tau^{-\frac12}$ spatial scale.
This scale turns out to be sufficiently small so that the commutator
between $\Re p_\phi$ and $\Im p_\phi$ becomes negligible, i.e. we gain
the principal normality.  On the other hand, in order to combine these
localized results we use the global $L^2$ estimate.

The part of the proof of Theorem~\ref{lpcee} which is obtained by
assembling together spatially localized estimates on the
$\tau^{-\frac12}$ is contained in the following two lemmas.

\begin{lemma}
  Under the assumptions in Theorem~\ref{lpcee} there is a para\-metrix
  $K$ for $P_\phi$ which satisfies the bounds
\begin{equation}
\| K f \|_{\tau^{\frac14} H^m_\tau \cap
  H^{m-\frac12-\rho(r),r}_{\tau} }  \lesssim 
 \|f \|_{\tau^{-\frac14} L^2+
  H^{\frac12+\rho(r),r'}_{\tau}}
\label{lk}\end{equation}
\begin{equation}
\| (I- P_\phi K) f \|_{\tau^{-\frac14} L^2} \lesssim \|f \|_{\tau^{-\frac14} L^2+
  H^{\frac12+\rho(r),r'}_{\tau}}
\label{lh}\end{equation}
\label{las}\end{lemma}

\begin{lemma}
Under the assumptions in Theorem~\ref{lpcee} the operator $P_\phi$ satisfies 
\begin{equation}
\| w \|_{H^{m-\frac12-\rho(r),r}_{\tau} }  \lesssim 
\| w \|_{\tau^{\frac14} H^m_\tau} +  \|P_\phi w \|_{\tau^{-\frac14} L^2} 
\label{lgs}\end{equation}
\label{lad}\end{lemma}
Before proving the lemmas we show how they can be combined
with the $L^2$ Carleman estimates to prove Theorem~\ref{lpcee}.

\begin{proof}[Proof of Theorem~\ref{lpcee}]
We need to prove that
\[
 \| v\|_{\tau^{\frac14} H^m_\tau \cap
  H^{m-\frac12-\rho(r),r}_{\tau} }^2 
\leq c   \| P_\phi(x,D,\tau)  v\|_{\tau^{-\frac14} L^2+
  H^{\frac12+\rho(r),{r'}}_{\tau}}^2,
\]
Let $K$ be as in Lemma~\ref{las}.  We decompose $v$ into
\[  
v= w + KP_\phi v. 
\]
 The function  $KP_\phi v$ satisfies the correct bounds by \eqref{lk}
while
\[
P_\phi w = (I - P_\phi K) P_\phi v.
\]
Using \eqref{lh} we  bound the right hand side in $L^2$: 
\[ 
\Vert P_\phi w \Vert_{\tau^{-1/4}L^2} \le \|P_\phi v \|_{\tau^{-\frac14} L^2+
  H^{\frac12+\rho(r),r'}_{\tau}}
\] 
But the $L^2$ Carleman estimate (\ref{cee}) allows us to 
also obtain an $L^2$ estimate for $w$,
\[
\|w\|_{\tau^{\frac14} H^m_\tau} \lesssim \|P_\phi w\|_{\tau^{-\frac14} L^2}.
\]
Hence we can use Lemma~\ref{lad} to obtain the correct estimate for
$w$ and conclude the proof of Theorem \ref{lpcee}.
\end{proof}

\begin{proof}[Proof of Lemmas~\ref{las},\ref{lad}]
  Without any restriction in generality we assume that
  $\tau^{-\frac12} \ll d(K,\partial \Omega)$.  Then we claim that it
  suffices to prove both lemmas in a ball of radius $\tau^{-\frac12}$.
  For this we consider a locally finite covering of $\Omega$ with
  balls of radius $\tau^{-\frac12}$,
\[
\Omega \subset \bigcup_j  B_j
\]
Correspondingly we consider a smooth partition of unity 
\[
1 = \sum_j \chi_j, \qquad \text{supp } \chi_j \subset B_j
\]
Suppose we know that for each $j$ there is a parametrix $K_j$
so that for $f$ supported in $B_j$ the function $K_j f$ is supported
in $2 B_j$ and the estimates (\ref{lk}), (\ref{lh}) hold. 
Then we can construct a parametrix $K$ for $P_\phi$ by
\[
K = \sum_j K_j \chi_j
\]
To obtain the estimates (\ref{lk}), (\ref{lh}) for $K$ we need to
verify that we can sum up the bounds for $K_j$ in $l^2$,
\[
\sum \|\chi_j f\|_{\tau^{-\frac14} L^2+
  H^{\frac12+\rho(r),r'}_{\tau}}^2 \lesssim \|f\|_{\tau^{-\frac14} L^2+
  H^{\frac12+\rho(r),r'}_{\tau}}^2
\]
respectively
\[
\| \sum K_j \chi_j f\|_{\tau^{\frac14} H^m_\tau \cap
  H^{m-\frac12-\rho(r),r}_{\tau} }^2 \lesssim   \sum 
\|K_j \chi_j f\|_{\tau^{\frac14} H^m_\tau \cap
  H^{m-\frac12-\rho(r),r}_{\tau} }^2
\]

Similarly, if (\ref{lgs}) holds for $w$ supported in $B_j$ then we can
conclude it holds in general using the previous inequality and the
additional estimate
\[
\sum_{j} \| \chi_j w\|_{\tau^{\frac14} H^m_\tau}^2 + 
\|P_\phi \chi_j w\|_{\tau^{-\frac14} L^2}^2 \lesssim 
\| w\|_{\tau^{\frac14} H^m_\tau}^2 + 
\|P_\phi w\|_{\tau^{-\frac14} L^2}^2
\] 
These three inequalities are easy exercises which are left for the
reader. We only observe that while the first two can be localized
further down to the $\tau^{-1}$ spatial scale, this would be
useless as a parametrix satisfying the right bounds cannot be
constructed on a smaller scale. On the other hand for the third bound
the $\tau^{-\frac12}$ scale is optimal because  the commutators
$[P_\phi,\chi_j]$ have to be controlled. 

To prove the lemmas in a ball of radius $\tau^{-\frac12}$ 
we rescale it to the unit ball. The rescaled operator is
\[
\tilde P_\phi(x,D,\mu) = P_\phi(\frac{x}\mu,D,\mu),  \qquad \mu = \tau^\frac12
\]
Then for Lemma~\ref{las} we need a parametrix $\tilde K$ for $\tilde P_\phi$
which satisfies
\begin{equation}
\| \tilde K f \|_{\mu^{\frac12} H^m_\mu \cap
  H^{m-\frac12-\rho(r),r}_{\mu} }  \lesssim 
 \|f \|_{\mu^{-\frac12} L^2+
  H^{\frac12+\rho(r),r}_{\mu}}
\label{lka}\end{equation}
\begin{equation}
\| (I- \tilde P_\phi K) f \|_{\mu^{-\frac12} L^2} \lesssim \|f \|_{\mu^{-\frac12} L^2+
  H^{\frac12+\rho(r),r'}_{\mu}}
\label{lha}\end{equation}
while for Lemma~\ref{lad} we need the estimate
\begin{equation}
\| w \|_{H^{m-\frac12-\rho(r),r}_{\tau} }  \lesssim 
\| w \|_{\mu^{\frac12} H^m_\mu} +  \|P_\phi w \|_{\mu^{-\frac12} L^2} 
\label{lgsa}\end{equation}

Within the unit ball the coefficients of $\tilde P_\phi$ vary by
$O(\mu^{-1})$. Then 
\[
\| (\tilde P_\phi(x,D,\mu) - \tilde P_\phi(0,D,\mu))w\|_{\mu^{-\frac12}
  L^2} \lesssim \| w \|_{\mu^{\frac12} H^m_\mu} 
\]
therefore without any loss we can freeze the coefficients of
$P_\phi$  and replace $\tilde P_\phi(x,D,\mu)$ by $\tilde
P_\phi(0,D,\mu)$.  This is principally normal by default.

The symbol 
\[
\tilde p_\phi(x,\xi,\mu) = p(0,\xi+i\mu \nabla \phi(0))
\]
is elliptic in the region $\tau \ll |\xi|$. Then the only region in
frequency where the problem is nontrivial is $\{ |\xi| \lesssim
\tau\}$. This is where we use the curvature condition.
For low frequencies $|\xi| \lesssim \mu$ we can use the parametrix $K$
given by Theorem~\ref{tpq}, while at higher frequencies $\tilde
P_\phi(0,D,\mu)$ is elliptic. More precisely, we consider a large
enough constant $C$ so that $\tilde P_\phi(0,D,\mu)$ is elliptic in
$\{|\xi| > C \mu\}$. Then we denote by $\chi$ a symbol supported in
$\{|\xi| > 2 C \mu\}$ and which equals $1$ in the region $\{|\xi| < C
\mu\}$.  The curvature condition in the hypothesis of Theorem~\ref{lpcee}
implies that we can use Theorem~\ref{tpq} to produce a parametrix $K_\mu$
in the region $\{|\xi| < 2 C \mu\}$.

Then we define the parametrix $K$ for $P_\phi(0,D,\mu)$ by
\[
K = K_\mu \chi (D) + P_\phi(0,D,\mu)^{-1}(1-\chi(D))
\]
The bounds for $K$ follow easily from the similar bounds in 
Theorem~\ref{tpq}  for $K_\mu$. This concludes the proof of
Lemma~\ref{las}. 

For Lemma \ref{lad} we use the same setup. The bound for
$(1-\chi(D)) w$ follows from ellipticity, while the bound for $\chi(D)
w$ is nothing but (\ref{estr}) in the context of Theorem~\ref{tpq}.
\end{proof}

An argument as the one in the proof of Theorem~\ref{lpcee} shows that
Theorem~\ref{lpce} is a consequence of the following counterparts of
Lemmas~\ref{las},\ref{lad}:

\begin{lemma}
  Under the assumptions in Theorem~\ref{lpce} there is a para\-metrix
  $K$ for $P_\phi$ which satisfies the bounds
\begin{equation}
\| K f \|_{\tau^{-\frac14} H^{m-1}_\tau \cap
  H^{m-1-\rho(r),r}_{\tau} }  \lesssim 
 \|f \|_{\tau^{\frac14} L^2+
  H^{\rho(r),r}_{\tau}}.
\label{xlk}\end{equation}
\begin{equation}
\| (I- P_\phi K) f \|_{\tau^{\frac14} L^2} \lesssim \|f \|_{\tau^{\frac14} L^2+
  H^{\rho(r),r'}_{\tau}}.
\label{xlh}\end{equation}
\label{xlas}\end{lemma}

\begin{lemma}
Under the assumptions in Theorem~\ref{lpce} the operator $P_\phi$ satisfies 
\begin{equation}
\| w \|_{H^{m-1-\rho(r),r}_{\tau} }  \lesssim 
\| w \|_{\tau^{-\frac14} H^{m-1}_\tau} +  \|P_\phi w \|_{\tau^{\frac14} L^2} 
\label{xlgs}\end{equation}
\label{xlad}\end{lemma}

\begin{proof}[Proof of Lemmas~\ref{xlas},\ref{xlad}]
  As in the proof of Lemmas~\ref{las},\ref{lad} we can still localize
  on the $\tau^{-\frac12}$ scale, and then rescale it back to the unit
  ball.  However, now we can no longer freeze the coefficients of the
  rescaled operator
\[
\tilde P_\phi(x,D,\mu) = P_\phi(\frac{x}\mu,D,\mu),
\]
because its characteristic set is nontrivial at high frequencies.
We divide the Fourier space in dyadic regions 
\[
D_\l = \{ |\xi| \leq 2 \mu\}, \qquad \l = \mu
\]
\[
D_\l = \{ \frac\l{2} \leq |\xi| \leq 2\l \}, \qquad \l = 2^j \mu,
\quad j > 0
\]
Correspondingly we consider a partition of unity
in frequency
\[
1 = \sum_{\l=2^j \mu}^{j \geq 0}  \chi_\l(\xi), \qquad \text{supp } \chi_\l \subset D_\l 
\]
We also consider symbols $\tilde \chi_\l$ with slightly larger support
which equal $1$ near the support of $\chi_\l$.  Write the principal
part of $P_\phi$ in the form
\[
P_\phi = \sum_{|\aa| = m} c_\aa(x) (D,\tau)^\aa
\]
where $\tau$ carries the same weight as a derivative.
Then 
\[
\tilde P_\phi = \sum_{|\aa| = m} c_\aa(x/\mu) (D,\mu)^\aa
\]
For $\l = 2^j \mu$ we define the 
regularized  coefficients
\[
c_{\aa,\l} = S_{< \l^{\frac12}} c_\aa
\]
where $ S_{< \l^{\frac12}}$ is a multiplier with smooth symbol which
selects the frequencies $\l^{\frac12}$ and smaller.  These do not
differ much from the original coefficients,
\begin{equation}
\|c_{\aa,\l} - c_\aa\|_{L^\infty} \lesssim \l^{-1} \|c_\aa\|_{C^2}
\label{caa}\end{equation}
We also introduce the modified operators, 
\[
\tilde P_{\phi,\l} = \sum_{|\aa| = m} c_{\aa,\l}(x/\mu) (D,\mu)^\aa
\]
somewhat in the spirit of the paradifferential calculus.

In the region $D_\l$ these symbols are in $\l^{m} \Sl$.
They also satisfy the principal normality condition. Indeed,
since $P$ has real coefficients, it follows that all terms in 
$\Im  \tilde P_{\phi,\l}$ contain at least one power of $\mu$.
Then 
\[
|\{ \Re  \tilde P_{\phi,\l}, \Im P_{\phi,\l}\}| \lesssim \l^{2(m-1)}
\qquad \text{in } D_\l
\]
By (\ref{caa}) the curvature condition in $D_\l$ is easily transfered
from $\tilde P_\phi$ to $P_{\phi,\l}$.  If $\l \approx \mu$ then this
is given in the hypothesis of the Theorem \ref{lpce}. If $\l \gg \mu$ then $\Re
\tilde P_\phi$ is a small perturbation of $\tilde P = P(x/\mu,D)$,
therefore we can use the curvature condition for $P$.

Now we are in a position to use Theorem~\ref{tpq} for $\l \approx \mu$
and Theorem~\ref{tpq1} for $\l \gg \mu$. We obtain parametrices $K_\l$
which satisfy
\[
\l^{m-1} \| K_\l  \chi_\l f\|_{L^2\cap \l^{\rho(r)} L^r} \lesssim
\|f\|_{L^2 + \l^{-\rho(r)} L^{r'}}
\]
\[
 \| (I - \tilde P_{\phi,\l}  K_\l)  \chi_\l  f\|_{L^2} \lesssim
\|f\|_{L^2 + \l^{-\rho(r)} L^{r'}}
\]
Now we can define the parametrix $K$ for $\tilde P_\phi$
as
\[
Kf = \sum_\l  \tilde \chi_\l K_\l  \chi_\l f
\]
It remains to prove that $K$ satisfies the desired bounds.
The bound for $K$ follows from the bound for $K_j$ and
Littlewood-Paley theory. Consider now the error estimates.
We have
\begin{eqnarray*}
I - \tilde P_\phi K &=& \sum_\l I -  \tilde P_{\phi,\l} \tilde \chi_\l K_\l  \chi_\l
\\ &=&  \sum_\l  ( \tilde P_\phi - \tilde P_{\phi,\l} )\tilde \chi_\l
K_\l  \chi_\l +  [\tilde P_{\phi,\l}, \tilde \chi_\l] K_\l  \chi_\l
\\ &+& \sum_\l \tilde \chi_\l ( I- \tilde P_{\phi,\l} K_\l)  \chi_\l
\end{eqnarray*}
For the first term we use (\ref{caa}), for the second an $L^2$ commutator estimate.
Finally the bound for the third is given by the similar estimates for
$K_\l$. For each term we also need to use the Littlewood-Paley theory.
This concludes the proof of Lemma~\ref{xlas}.

For Lemma~\ref{xlad} we use 
Theorem~\ref{tpq} for $\l \approx \mu$ respectively  
Theorem~\ref{tpq1}
for $\l \gg \mu$ (In both case we use the theorems in the form of 
estimate \eqref{estr}). These imply that
\[
\l^{m-1} \|\chi_\l w\|_{ \l^{\rho(r)} L^r} \lesssim 
\l^{m-1} \|\tilde \chi w\|_{L^2} + \|\tilde P_{\phi,\l} \tilde \chi w\|_{L^2}
\]
which after some commuting we can square and sum up using
the Littlewood-Paley theory.
\end{proof}

\subsection{ The Laplace equation}

Consider a second order elliptic operator in $\R^n$,
\[
P = \partial_j g^{jk}(x) \partial_k
\]
Any surface $\Gamma$ is strongly pseudoconvex with respect to $P$.
If $\phi$ is a strongly pseudoconvex function with respect to $P$
with $\nabla \phi \neq 0$ then 
\[
p_\phi(x,\xi,\tau) = g^{jk} \xi_j \xi_k - \tau^2 g^{jk} \partial_j
\phi \partial_k \phi + 2 i \tau g^{jk} \xi_j \partial_k \phi 
\]
The characteristic set of the real part is an ellipsoid centered at
the origin, while the characteristic set of the imaginary part is a
plane through the origin. The characteristic set of the full operator
is the intersection of the ellipsoid with the plane, and has $n-2$
nonvanishing curvatures.  Therefore, we can use Theorem~\ref{lpcee}
with $k=0$:

\begin{theorem}
Let $P$ be a second order elliptic operator with $C^1$ coefficients.
Let $\phi$ be a strongly pseudoconvex function with respect to $P$.
Then for compactly supported $u$ we have
\[
\| \etf u \|_{H^{\frac2{n+2},\frac{2(n+2)}{n}}_\tau\cap
  \tau^{\frac14}H^1_\tau}  \lesssim \|\etf P
u\|_{H^{-\frac2{n+2},\frac{2(n+2)}{n+4}}_\tau+
  \tau^{-\frac14}H^{-1}_\tau}, \quad \tau \geq \tau_0
\]
\end{theorem}
Note that this is not precisely in the form stated in \eqref{cee} but
it can be easily obtained from it by conjugating $P_\phi$ with a first
order elliptic multiplier because $P$ is in divergence form. We prefer
the above formulation because of its symmetry. Similar adjustments are
made in all the other examples we consider.

Applied to unique continuation problems, this yields

\begin{theorem}
Let $P$ be a second order elliptic operator with $C^1$ coefficients.
Let $\Gamma$ be a smooth surface. Then unique continuation
for $P+V$ across $\Gamma$ holds for all potentials $V$ which
have the multiplicative mapping property
\[
V: H^{\frac2{n+2},\frac{2(n+2)}{n}} \to  H^{-\frac2{n+2},\frac{2(n+2)}{n+4}}
\]
\end{theorem}
This includes the case $V \in L^{\frac{n}{2}}$, for which the result
was proved by Wolff~\cite{780.35015}. This can be relaxed to a slightly larger
Morrey space.

For strong unique continuation problems and problems involving
gradient potentials we refer the reader to Wolff~\cite{MR96c:35068},
the authors paper \cite{MR2001m:35075}, and the references therein.

\subsection{The wave equation}

Consider a second order hyperbolic operator in $\R^{d+1}$,
\[
P = \partial_j g^{jk}(x) \partial_k
\]
where the matrix $g^{ij}$ has signature $(d,1)$.  Which
noncharacteristic surfaces $\Gamma$ are strongly pseudoconvex with
respect to the $P$ ? One needs to
distinguish between space like surfaces ($g^{jk} \partial_j \phi
\partial_k \phi < 0$) and time-like surfaces ($g^{jk} \partial_j \phi
\partial_k \phi < 0$). All space-like surfaces are pseudoconvex, since
the Cauchy problem with initial data on a space-like surface is
well-posed. For time-like surfaces, on the other hand, the condition
(\ref{is}) is trivially fulfilled but (\ref{rsf}) may or may not hold.

If $\phi$ is a strongly pseudoconvex function with respect to $P$
with $\nabla \phi \neq 0$ then 
\[
p_\phi(x,\xi,\tau) = g^{jk} \xi_j \xi_k - \tau^2 g^{jk} \partial_j
\phi \partial_k \phi + 2 i \tau g^{jk} \xi_j \partial_k \phi 
\]
The characteristic set of the real part is a hyperboloid, which has
$d$ nonvanishing curvatures.  The characteristic set of the full
operator is the intersection of the hyperboloid with a plane, and has
$d-1$ nonvanishing curvatures.  Hence we can use Theorem~\ref{lpce}
with $k=1$, $n = d+1$:

\begin{theorem}
Let $P$ be a second order hyperbolic operator with $C^2$ coefficients.
Let $\phi$ be a strongly pseudoconvex function with respect to $P$.
Then for compactly supported $u$ we have
\[
\| \etf u \|_{L^\frac{2(d+1)}{d-1} \cap \tau^{-\frac14}H^\frac12_\tau}
\lesssim \|\etf P u\|_{
L^\frac{2(d+1)}{d+3}+ \tau^{\frac14}H^{-\frac12}_\tau}, \quad \tau > \tau_0
\]
\end{theorem}

Applied to unique continuation problems, this gives

\begin{theorem}
Let $P$ be a second order hyperbolic operator with $C^2$ coefficients.
Let $\Gamma$ be a smooth surface which is strongly pseudoconvex with
respect to $P$. Then unique continuation
for $P+V$ across $\Gamma$ holds for all potentials $V \in L^\frac{d+1}2$.
\end{theorem}
This improves an earlier result in \cite{MR97i:35012}.  One can also
produce versions of this with potentials in mixed norm spaces.

\subsection{The heat equation}

Consider a second order parabolic operator in $\R \times \R^d$,
\[
P = \partial_t - \partial_j g^{jk}(t,x) \partial_k
\]
We denote by $\sigma$ and $\xi$ the time, respectively the space Fourier
variable. Then the symbol of $P$ is
\[
p(t,x,\sigma,\xi) = -i\sigma + g^{jk} \xi_j \xi_k
\]
This vanishes only at $\sigma=0,\ \xi=0$ so we should treat $P$ as an
elliptic operator. However, the results in Theorem~\ref{lpcee} cannot 
be applied directly due to the different scaling associated to the 
heat operator. Instead one needs to adapt that setup to the current
problem. This is discussed in what follows.

First we note that a time derivative is roughly equivalent with two
time derivatives. Hence the size of the frequency  is now
$(|\sigma|^2+|\xi|^4)^\frac14$ and the dyadic regions in frequency correspond
to $(|\sigma|^2+|\xi|^4)^\frac14 \approx 2^j$. The weighted Sobolev spaces
are redefined accordingly,
\[
\|u\|_{H^k_\tau} = \| (|D_t|^2+|D_x|^4+\tau^4)^\frac{k}4 \hat{u} \|_{L^2}
\]
and similarly for $H^{k,p}_\tau$.

The principal symbol of the conjugated operator $P_\phi$ now has
the form
\begin{eqnarray*}
p_\phi(t,x,\sigma,\xi,\tau) &=& p(t,x,\sigma,\xi+i\tau \nabla \phi)
\\ &=& g^{jk} \xi_j \xi_k - \tau^2 g^{jk} \partial_j
\phi \partial_k \phi + i (\sigma + 2 \tau g^{jk} \xi_j \partial_k \phi )
\end{eqnarray*}
Observe that the time derivatives of $\phi$ do not appear
in this formula, as the terms containing them are lower order
terms. By the same token, time derivatives are also excluded
from the definition of the Poisson bracket, namely
\[
\{ \Re p_\phi,\Im p_\phi\} = ( \Re p_\phi)_x ( \Im p_\phi)_\xi -
 (\Re p_\phi)_\xi  (\Im p_\phi)_x
\]

Another adjustment one needs to make concerns the scale of the
localization in the $L^2$ Carleman estimates, which now can be done on
parabolic balls of size $\tau^{-1}\times (\tau^{-\frac12})^d$. Because
of this less time regularity for the coefficients is needed.

Taking all these considerents into account the analysis proceeds very
much like in the case of Theorem~\ref{lpcee}.  The operator $P_\phi$
is elliptic at all frequencies larger than $\tau$, so the analysis
must concentrate on the frequency region
\[
|\sigma|+|\xi|^2 \lesssim \tau^2
\]

A short computation shows that any surface $\Gamma$ which is not
tangent to the time slices is strongly pseudoconvex with respect to
$P$.  If $\phi$ is a strongly pseudoconvex function with respect to
$P$ with $\nabla _x \phi \neq 0$ then for fixed $t$ and $x$ the
characteristic set of $\Re p_\phi$ is a cylinder on top of an
ellipsoid, while the characteristic set of $\Im p_\phi$ is an oblique
plane.  Then the characteristic set of the full operator is an $d-1$
dimensional ellipsoid which is the intersection of the cylinder with
the plane, and has $d-1$ nonvanishing curvatures.  Therefore, we use
(an adapted version of) Theorem~\ref{lpcee} in $d+1$ dimensions with
$k=0$ and $n=d+1$.

\begin{theorem}
  Let $P$ be a second order  parabolic operator whose
  coefficients are $C^1$ in $x$ and $C^{\frac12}$ in time.  Let $\phi$
  be a strongly pseudoconvex function with respect to $P$ with
  $\nabla_x \phi \neq 0$.  Then for compactly supported $u$ we have
\[
\| \etf u \|_{H^{\frac1{d+3},\frac{2(d+3)}{d+1}}_\tau\cap
  \tau^{\frac14}H^1_\tau}  \lesssim \|\etf P
u\|_{H^{-\frac1{d+3},\frac{2(d+3)}{d+5}}_\tau+
  \tau^{-\frac14}H^{-1}_\tau}, \quad
\tau > \tau_0
\]
\end{theorem}

Applied to unique continuation problems, this yields

\begin{theorem}
Let $P$ be a second order parabolic operator whose
  coefficients are $C^1$ in $x$ and $C^{\frac12}$ in time. 
Let $\Gamma$ be a smooth surface. Then unique continuation
for $P+V$ across $\Gamma$ holds for all potentials $V$ which
have the multiplicative mapping property
\[
V: H^{\frac1{d+3},\frac{2(d+3)}{d+1}} \to  H^{-\frac1{d+3},\frac{2(d+3)}{d+5}}
\]
\end{theorem}
This includes the case $V \in L^{\frac{d+2}{2}}$, for which the result
is new. This can be relaxed to a slightly larger Morrey space.
Applications of these ideas to strong unique continuation problems are
contained in a forthcoming paper of the authors.

\begin{remark}
  By using the mixed norm estimates in part (b) of Theorem~\ref{tpq1}
  we can also obtain versions of these results with potentials in mixed
  norm spaces. For simplicity we state a weaker form of the Carleman 
estimates,
\[
\| \etf u \|_{L^q L^r \cap
  \tau^{\frac14}H^1_\tau}  \lesssim \|\etf P
u\|_{L^{q'} L^{r'}+
  \tau^{-\frac14}H^{-1}_\tau}, \quad
\tau > \tau_0
\]
which holds whenever\footnote{The exponent $q=2$ is now allowed
because the pair $(q,r)$ is no longer the endpoint.}
\[
\frac{2}{q}+\frac{d}{r} = \frac{d}2, \qquad 2 \leq q \leq \infty
\]
This yields the unique continuation result  for $P+V$ provided 
that  $V \in L^{\tilde q}_t L^{\tilde r}_x$ 
where the exponents $\tilde q$ and $\tilde r$ satisfy the scaling
relation\footnote{For $d =1$ one needs the additional restriction
  $\tilde q
  \geq 2$}
\[
\frac{2}{\tilde q} + \frac{d}{\tilde r} = 2, \qquad 1 \leq \tilde q,
\tilde r \leq \infty
\]
\end{remark}

\subsection{The Schr\"odinger equation}

Here we consider the second order Schr\"odinger operator in $\R \times
\R^d$,
\[
P = i \partial_t - \partial_j g^{jk}(x) \partial_k
\]
For this one needs to use the same setup as in the case of parablic
equations.

  We say that a surface $\Gamma$ is noncharacteristic if it
is not tangent to the time slices. As in the case of the wave
equation, the condition (\ref{is}) is always satisfied, but the
condition~\eqref{rsf} may or may not hold.

If $\phi$ is a strongly pseudoconvex function with respect to $P$ with
$\nabla_x \phi \neq 0$ then
\[
p_\phi(t,x,\sigma, \xi,\tau) = \sigma - g^{jk} \xi_j \xi_k - \tau^2 g^{jk} \partial_j
\phi \partial_k \phi +  2i \tau g^{jk} \xi_j \partial_k \phi 
\]
The characteristic set of the real part is a paraboloid, which has $d$
nonvanishing curvatures. The characteristic set of the full operator
is an $d-1$ dimensional ellipsoid which is the intersection of the
paraboloid with a vertical plane, and has $d-1$ nonvanishing
curvatures. Then (a variant of) Theorem~\ref{lpce} gives

\begin{theorem}
$\!$ Let $P$ be a second order  Schr\"odinger operator whose
  coefficients are $C^2$ in $x$ and $C^1$ in time.  Let $\phi$
  be a strongly pseudoconvex function with respect to $P$ with
  $\nabla_x \phi \neq 0$.  Then for compactly supported $u$ we have
\[
\| \etf u \|_{L^{\frac{2(d+2)}{d}}\cap
  \tau^{-\frac14}H^\frac12_\tau}  \lesssim \|\etf P
u\|_{L^\frac{2(d+2)}{d+4}+
  \tau^{\frac14}H^{-\frac12}_\tau}, \quad
\tau > \tau_0
\]
\end{theorem}

Applied to unique continuation problems, this yields

\begin{theorem}
$\!$ Let $P$ be a second order  Schr\"odinger operator whose
  coefficients are $C^2$ in $x$ and $C^1$ in time.  
Let $\Gamma$ be a smooth surface. Then unique continuation
for $P+V$ across $\Gamma$ holds for all potentials $V \in
L^\frac{d+2}{2}$.
\end{theorem}
 One can also produce versions of this result involving
mixed norm spaces.

\begin{remark}
Using the mixed norm version of Theorems~\ref{tpq1},~\ref{tpq} 
one obtains the Carleman estimates
\[
\| \etf u \|_{L^q L^r \cap
  \tau^{\frac14}H^1_\tau}  \lesssim \|\etf P
u\|_{L^{q'} L^{r'}+
  \tau^{-\frac14}H^{-1}_\tau}, \quad
\tau > \tau_0
\]
which holds whenever\footnote{It would be interesting to see if the
  exponent $q=2$ is admissible.}
\[
\frac{2}{q}+\frac{d}{r} = \frac{d}2, \qquad 2 < q \leq \infty
\]
This yields the unique continuation result  for $P+V$ provided 
that  $V \in L^{\tilde q}_t L^{\tilde r}_x$ 
provided the exponents $\tilde q$ and $\tilde r$ satisfy the scaling
relation\footnote{For $d =1$ one needs the additional restriction
  $\tilde q
  \geq 2$}
\[
\frac{2}{\tilde q} + \frac{d}{\tilde r} = 2, \qquad 1 < \tilde
q,\tilde r \leq \infty
\]
\end{remark}

\end{document}